\documentclass{pspum-l} \usepackage{amsmath, amsthm, amscd, amsfonts,amssymb} 
\numberwithin{equation}{section} 

\newcommand{\teq}{\arabic{section}.\arabic{equation}}
\newcommand{\teql}{\Alph{section}.\arabic{equation}}

\newcommand{\sqr}[2]{{\vcenter{\vbox{\hrule height.#2pt\hbox{\vrule width.#2pt
height#1pt \kern#1pt\vrule width.#2pt}\hrule height.#2pt}}}}

\newcounter{eqcount}

\renewcommand{\labelenumi}{{{\rm (\teq \alph{enumi})}}} 
\newenvironment{edesc}{\refstepcounter{equation}\begin{enumerate}}%
{\end{enumerate}}
\newenvironment{triv}{\refstepcounter{equation}\begin{list}%
{{\hbox{\rm(\teq)\ }}} \item }{\end{list}}

\newcommand{\ring}[1]{{\mathbb #1}}
\newcommand\bZ{{\ring{Z}}}

\newcommand\bC{{\ring{C}}} \newcommand\bR{{\ring{R}}}
 \newcommand\bQ{{\ring{Q}}}

\newcommand{\csp}[1]{{\mathbb #1}}
\newcommand{\tsp}[1]{{\mathcal #1}}

\newcommand{\bM}{\csp{M}}
\newcommand{\prP}{\csp{P}}

\newcommand{\sQ}{\tsp{Q}}
\newcommand{\sP}{{\tsp {P}}} 
 
\newcommand{\sT}{{\tsp {T}}} \newcommand{\sH}{{\tsp {H}}}
\newcommand{\sX}{{\tsp {X}}} 
\newcommand{\sM}{{\tsp {M}}} \newcommand{\sD}{{\tsp {D}}}

 \newcommand{\bN}{{\csp {N}}}

\newcommand{\eql}[2]{{\rm (\ref{#1}\ref{#2})}} 

\newcommand{\vect}[1]{{\pmb #1}} 
 \newcommand{\bg}{\vect{g}}
 \newcommand{\be}{{\vect{e}}}
\newcommand{\bp}{{\vect{p}}} \newcommand{\bx}{{\vect{x}}}
 \newcommand{\bw}{{\vect{w}}}
 \newcommand{\bz}{{\vect{z}}}
\newcommand{\bm}{{\vect{m}}}  
\newcommand{\bh}{{\vect{h}}}  \newcommand{\bu}{{\vect{u}}}
\newcommand{\row}[2]{{#1_1,\ldots,#1_{#2}}}

\newcommand{\smatrix}[4]{{\big(\begin{array}{cc}
\!\lower2pt\hbox{$\scriptstyle#1$} &\lower2pt\hbox{$\scriptstyle#2$}\!
\\\! \raise2pt\hbox{$\scriptstyle#3$} &\raise2pt\hbox{$\scriptstyle#4$}
\!\end{array}\big)}}

\newcommand{\texto}[1]{{\textr{#1}}}
\newcommand{\GL}{\texto{GL}} \newcommand{\SL}{\texto{SL}}
 \newcommand{\ind}{\texto{ind}}
\newcommand{\PSL}{\texto{PSL}} \newcommand{\PGL}{\texto{PGL}}
 
 \renewcommand{\ni}{\texto{Ni}}
 \newcommand{\Pic}{\texto{Pic}}
\newcommand{\textr}[1]{{\text{\rm #1}}}
\newcommand{\tr}{\textr{tr}} 
\newcommand{\abs}{\textr{abs}}  \newcommand{\cyc}{\textr{cyc}}
 
\newcommand{\pC}{{\textr{C}}} \newcommand{\inn}{\textr{in}}
 \newcommand{\Aut}{\textr{Aut}}
\newcommand{\pr}{\textr{pr}}

\newcommand{\rd}{\texto{rd}}

\newcommand{\tG}[1]{{}_{#1}\tilde G}

\newcommand{\GAP}{{\bf GAP}}

\newcommand{\sph}{{\vphantom 1}}

\newcommand{\textb}[1]{{\text{\bf #1}}}
\newcommand{\bfC}{{\textb{C}}}
\newcommand{\longmapright}[2]{\smash{\mathop{\hbox to
#2pt{\rightarrowfill}}\limits^{#1}}}
\newcommand{\longmapleft}[2]{\smash{\mathop{\hbox to
#2pt{\leftarrowfill}}\limits^{#1}}}

\newcommand{\mapright}[1]{\smash{\mathop{\longrightarrow}\limits^{#1}}}
\newcommand{\mapleft}[1]{\smash{\mathop{\longleftarrow}\limits^{#1}}}

\newcommand{\np}{{+}}   \newcommand{\nm}{{-}}

\newcommand{\lrang}[1]{{\langle #1\rangle}}

\newcommand{\eqdef}{\stackrel{\text{\rm def}}{=}}

\newfont{\sevenrm}{cmr7}
\newfont{\bsevenrm}{cmbx7}
\newfont{\mathseven}{cmsy7}
\newfont{\bigmath}{cmsy10 scaled 1200}
\newfont{\fiverm}{cmr5}
\newfont{\bfiverm}{cmbx5}
\newfont{\hel}{cmbx10 scaled 1200}
\newfont{\eu}{eufb10}
\newfont{\sseu}{eufm5}
\newfont{\seu}{eufm7}
\newfont{\Cal}{cmmib10}
\newfont{\sCal}{cmmib7}
\newfont{\zch}{eusb10}

\theoremstyle{plain}
\newtheorem{thm}{Theorem}[section] 
\newtheorem{lem}[thm]{Lemma}
\newtheorem{princ}[thm]{Principle}

\newtheorem{prop}[thm]{Proposition}
\newtheorem{cor}[thm]{Corollary}

\theoremstyle{definition}
\newtheorem{defn}[thm]{Definition}
\newtheorem{exmp}[thm]{Example}
\newtheorem{guess}[thm]{Conjecture}
\newtheorem{quest}[thm]{Question}
\newtheorem{prob}[thm]{Problem}

\theoremstyle{remark}
\newtheorem{rem}[thm]{Remark}

\newcommand{\xs}{\times^s\!}
\newcommand{\wsp}{{$\,$---$\,$}} 

\newcommand{\twid}{\~{\hbox{$\,$\!\!}}}

\begin{document}

\newcommand{\sHt}[2]{{\sH(A_{{#1}},\bfC_{3^{#2}})}}
\newcommand{\sHtp}[2]{{\sH_+(A_{{#1}},\bfC_{3^{#2}})}}
\newcommand{\sHtm}[2]{{\sH_-(A_{{#1}},\bfC_{3^{#2}})}}
\newcommand{\sHtpm}[2]{{\sH_{\pm}(A_{{#1}},\bfC_{3^{#2}})}}

\newcommand{\pari}{{\textr{par}}}
\newcommand{\Ct}[1]{{\bfC_{{3^{#1}}}}}
\let\phi=\varphi
\newcommand{\balpha}{{\vect{\alpha}}} \newcommand{\bbeta}{{\vect{\beta}}}
\newcommand{\et}{{\text{\sl et}}}
\newcommand{\ram}{{\text{\rm ram}}}
\newcommand{\G}{{\ring{G}}}
\newcommand{\app}{{App}}
\newcommand{\sV}{{\tsp V}}
\newcommand{\HH}{{{}_+\sH_\infty}}
\newcommand{\NNi}{{{}_+\ni_\infty}} 
\newcommand{\CC}{{{}_+\pC}}
\newcommand{\bfCC}{{{}_+\bfC}}
\newcommand{\mq}{\textr{MQ}}
\newcommand{\one}{{\pmb 1}}
\newcommand{\C}{\text{\rm C}}
\newcommand{\cha}{{\text{\rm char}}}
\newcommand{\lcm}{{\text{lcm}}}
\newcommand{\Out}{{\text{\rm Out}}}
\newcommand{\Sp}{{\text{\rm Sp}}}
\newcommand{\Conj}{{\text{\rm Conj}}}
\newcommand{\red}{{\text{\rm rd}}}
\newcommand{\sym}{{\text{\rm Sym}}}

\newcommand{\hc}{{\text{\rm ${\frac1 2}$-canonical\ }}}
\newcommand{\N}[1]{{\{1,\dots,{#1}\}}} \newcommand{\mt}[1]{{\cite[#1]{FrMT}}}
\newcommand{\Spin}{{\text{\rm Spin}}}
\setcounter{tocdepth}{2}
\newcommand{\sh}{{\text{\bf sh}}}
\newcommand{\mb}[1]{{\mbox{\boldmath{#1}}}}
\font\eightrm=cmr8  \font\eightit=cmsl8 \let\it=\sl 
\newcommand{\fva}{\cite{FrV2}}
\newcommand{\fval}[1]{\cite[#1]{FrV2}}
\newcommand{\fvb}{\cite{FrV3}}
\newcommand{\fvbl}[1]{\cite[#1]{FrV3}}
\newcommand{\reg}{\text{\rm reg}}
\newcommand{\Gen}[2]{{\text{\rm Gen}(#1,#2)}}
\newcommand{\HMGen}[2]{\text{\rm H-MGen}(#1,#2)}
\newcommand{\HM}{\text{\rm H-M}}
\newcommand{\Pow}{\text{\rm Pow}}
\newcommand{\PowGen}[2]{{\text{\rm PowGen}(#1,#2)}}
\newcommand{\gu}{\text{g-u}}
\newcommand{\tpi}[5]{{\theta\Bigl[\!\!\begin{array}{c}\raise-
2pt\hbox{{$#1$}}\\\raise2pt\hbox{{$#2$}}
\end{array}\!\!\Bigr](#3,\Pi_{#4,#5})}}
\newcommand{\pdelta}{{\pmb \delta}}
\newcommand{\pepsilon}{{\pmb \epsilon}} 
\newcommand{\bn}{{\pmb n}}

\title[Lifting Invariants]{Alternating groups and \\ moduli space lifting Invariants}

\author[M.~D.~Fried]{Michael D.~Fried, Emeritus UC Irvine}

\date{\today}

\address{Math.~Dept., Emeritus, UC Irvine}
\email{mfried@math.uci.edu, mfri4@aol.com}

\begin{abstract} The genus of a curve discretely separates decidely different algebraic
relations in two variables to focus us on the connected moduli space $\sM_g$. Yet, modern applications      
also require a data variable (function) on the curve. The resulting 
spaces are versions, depending on our needs for this data variable, of {\sl Hurwitz spaces\/}. A {\sl Nielsen class}
(\S\ref{nz.1}) consists of $r\ge 3$ conjugacy classes $\bfC$ in the data variable monodromy $G$. It generalizes the
genus.  

Some Nielsen classes define  connected spaces. To detect, however, the components of others requires further subtler
invariants. We regard our Main Result  (MR) as level 0 
of Spin invariant information on moduli spaces. 

In the MR, $G=A_n$ (the alternating group), $r$ counts the data variable branch points and
$\bfC=\bfC_{3^r}$ is $r$ repetitions of the 3-cycle conjugacy class. This Nielsen class defines two spaces called
absolute and inner: 
$\sH(A_n,\bfC_{3^r})^\abs$ of degree $n$, genus 
$g=r-(n-1)>0
$ covers  and $\sH(A_n,\bfC_{3^r})^\inn$ parametrizing Galois closures of such covers. The parity of a spin invariant
precisely identifies the two components of each space. The inner result is the deeper.  

We examine the effect of combining the MR, \cite{AP}  and \hc\ classes  on $\sM_g$.  First: 
\S\ref{maxAltExt} considers an  analog of a
famous conjecture of Shafarevich: With $H$ the
composite group of all Galois extensions $K/Q$ with group some alternating group, does the canonical map  $G_\bQ\to H$ 
have pro-free kernel.  Second: Thm.~\ref{thetaNullRes} produces  nonzero  automorphic ($\theta$-null power)
functions on the reduced Hurwitz spaces 
$\sH_+(A_n,\bfC_{3^r})^{\abs,\red}$ (resp.~$\sH_-(A_n,\bfC_{3^r})^{\abs,\red}$) when $r$ is even (resp.~odd), for either $g=1$ or $n\ge 12g+4$ . 
\end{abstract}

\subjclass[2000]{Primary  11F32,  11G18, 11R58; Secondary 20B05, 
20C25, 20D25, 20E18, 20F34} 

\keywords{Moduli spaces of covers, Braid group, mapping class group, $A_n$, 
$\Spin_n$, Frattini covers, Serre's lifting invariant, $\Theta$-characteristics}

\thanks{Support from BiNational Science Foundation
Grant--Israel \#87-00038,  Alexander von
Humboldt Found., Institut f\"ur
Experimentalle Math., July 1996, from 
NSF \#DMS-99305590 and \#DMS-0202259 and the {\sl Workshop on algebraic curves and 
monodromy}, April 2004 in Milan.}  

\maketitle

\!\!\!\!\!\!\tableofcontents


\section{Introduction and notation} 

\subsection{Nielsen class notation} \label{nz.1} \!\!\S\ref{nz.1} reviews  our Main Result (MR) 
(\S\ref{quickReview}). This section's remainder reminds of Nielsen classes, Hurwitz spaces and
the small lifting invariant used in the MR.  We  repeatedly use that components of Hurwitz
spaces translate to braid (or Hurwitz monodromy) orbits on Nielsen classes. 

\subsubsection{A quick review of results} \label{quickReview} 
\S \ref{suppLems}, \S\ref{ns.2} and  \S\ref{nf.4}  treat  
the space of  projective line covers with $r$ ($\ge n-1$) 3-cycles as branch cycles
($n\ge 4$). For $n\ge 5$, $\Ct r$ denotes
$r$ repetitions of the conjugacy class of 3-cycles in $A_n$. 
When $n=4$ (or 3) there are two conjugacy classes of 3-cycles, so  $\bfC_{3^r}$ is
ambiguous. \S\ref{nt} states precise results in that hard case, so crucial to the
complete analysis. To simplify,  assume here  
$n\ge 5$. \S\ref{nt.2} shows  each allowable $\bfC_{3^r}$ replacement when $n=4$ has a
comparable result.  There are results for both degree
$n$ covers and Galois covers (degree
$n!/2$). 

We use both cases, denoted respectively  $\sH(A_n,\bfC_{3^r})^\abs$ and
$\sH(A_n,\bfC_{3^r})^\inn$. The MR says the following for either. A 
\hc  class (spin) invariant (noted by a
$\pm$ subscript) completely determines components of ${\sHt n r}\!^*$, $*=\inn$ or $\abs$.  
\begin{itemize}
\item Thm.~\ref{thmB}: For $r\!\ge n\!\ge 5$, 
\!${\sHt nr}\!^*$ has 2 components,
${\sHtpm nr}\!^*$. 
 \item
Thm.~\ref{thmA}:
${\sHt n{n\nm1}}\!^*$ has exactly one  \!(spin \!$(-1)^{n\nm1}$) component. \end{itemize}  
Cor.~\ref{nth.1} shows these spaces have useful moduli properties on their own.   
\begin{itemize} \item Each 
of \!${\sHtp nr}\!^*$ and ${\sHtm nr}\!^*$ has definition field $\bQ$. \item A dense subset of
$\sH_{\pm}(A_n,\Ct r)^\abs(\bar
\bQ)$    give $(A_n,S_n,\Ct r)$  {\sl realizations\/}.
\end{itemize}
Thm.~\ref{thetaNullRes} uses that  for some $(r,n)$,  $\sH_{\pm}(A_n,\Ct r)^\abs$ dominates  $\sM_g$, the
moduli of genus $g=r- (n-1)$ curves. That produces Hurwitz families with nonzero  
$\theta$-nulls.

\subsubsection{Nielsen class preliminaries} \label{gpNC} Let $G$ be a subgroup of $S_n$. For $g\in G$,  we say
$g$ contains $i$ (or $i$ is in the support of $g$) if a nontrivial disjoint
cycle of $g$ contains $i$. Also, for $\bg\in G^r$, $\Pi(\bg)=g_1\cdots g_r$
and    $\lrang{\bg}$ is the group the entries of $\bg$  generate. 
Permutations from $S_n$ act on the right of integers. 

Now consider  
 $r$ conjugacy
classes from $G$: $\bfC=(\row 
{\C} 
r)$, often one conjugacy class  
repeated many times. The main definitions don't depend on the order of their listing. For
example, for  $\bg\in G^r$,  $\bg\in \bfC$
means entries of $\bg$ are in the conjugacy classes of $\bfC$ in
{\sl some\/} order. This gives the {\sl
Nielsen class\/}, $\ni(G,\bfC)$ of $(G,\bfC)$:
$$\{\bg\in G^r\mid \Pi(\bg)=1,\ \bg\in \bfC, \text{ and }
\lrang{\bg}=G\}.$$

Denote a (ramified) cover of the sphere $\prP^1=\prP^1_z$ by a nonsingular connected curve $X$ by  
$\phi:X\to \prP_z^1$. Then $\phi$ has a degree, $\deg(\phi)=n$, and $G=G_\phi$, a transitive
subgroup \wsp the {\sl geometric monodromy\/} of $\phi$ \wsp of $S_n$ (symmetric group on $n$ letters).
Denote the set $\{g\in G\mid (1)g=1\}$ -- the stabilizer of 1 -- by $G(1)$. 

A {\sl
branch point\/} of $\phi$ is a $z'\in \prP^1_z$ for which the fiber $X_{z'}$ over $z'$ has fewer than 
$n$ points. So, 
$\phi$ has a {\sl branch point set\/} $\bz_\phi=\{\row z r\}\in (\prP^1)^r\setminus
\Delta_r/S_r =U_r$. 

Given a labeling of points of $X$ over $z_0$, the  following data attaches an $r$-tuple $\bg\in G^r$ to
$\phi$: 
 {\sl classical generators\/} $[P_1],\dots, [P_r]$ of the fundamental
group 
$\pi_1(\prP^1\setminus \bz, z_0)$ \cite[\S1.2]{BaFr}.  
Classical generators  present $\pi_1(\prP^1\setminus \bz, z_0)$ modulo one relation: 
$[P_1]  \dots[P_r]=1$. These topics have down-to-earth treatments at \cite{Fr08b}. 

\begin{defn}[branch cycles] \label{branchcycles} Any cover $\phi$ then corresponds to  a homomorphism 
$\psi_\phi: \pi_1(\prP^1_z\setminus\{\bz_\phi\},z_0)\to G$. 
The classical generators then assign  
$\bg\in
\ni(G,\bfC)$  by $\psi_\phi(P_i)=g_i$, $i=1,\dots,r$:  a {\sl branch cycle description\/} of $\phi$.
\end{defn}

Another set of classical generators 
will produce a different $\bg$ for $\phi$, yet in the same Nielsen
class. We say $\phi$  is in $\ni(G, \bfC)$. This with other background topics, studded
with classically motivated examples, is in
\cite[Chap.~4]{FrBook}. 

In a connected family of $r$-branched covers, we expect $\bz_\phi$ to move with $\phi$, while
$\ni(G_\phi,\bfC_\phi)$ is constant in $\phi$. We say $\phi$ is in the Nielsen class. A 
grasp on categories of covers requires useful equivalence relations. We use {\sl absolute, inner\/}, 
and their reduced versions in this paper (\S \ref{absEquiv} and \S \ref{rnz.4}).   To each equivalence of covers in a
given Nielsen class, there is a corresponding equivalence on the Nielsen class.  

Any of  
\cite[\S3.1]{BaFr}, \cite{Fr3}, \cite[\S I.F]{Fr4}, \cite[Chap.~I]{MM},
\cite[p.~60]{Se1} or \cite[Thm.~4.32]{VB}  explain Riemann's Existence
Theorem (RET):  how equivalence classes of $\phi\,$s branched over a fixed $\bz\in U_r$  correspond 
one-one to Nielsen class representatives modulo equivalence.

\subsubsection{Families of covers} \label{covs-reps} 
I explain absolute and inner equivalence for the families in Thm.~\ref{thmA} and \ref{thmB}. A family of 
$\ni(A_n,\bfC_{3^r})^\abs$ (resp.~$\ni(A_n,\bfC_{3^r})^\inn$) covers over a space $S$ is a degree $n$
(resp.~Galois, with group $A_n$) cover 
$\Phi: \sT\to S\times
\prP^1_z$, where the fiber of $\sT$ (resp.~$\sT/A_n(1)$) over $s\times \prP^1_z$, $s\in S$, is a cover
with 3-cycles as branch cycles. This determines 
$\Psi_\Phi: S
\to
\sH(A_n,\bfC_{3^r})^\abs$ (resp.~$\to \sH(A_n,\bfC_{3^r})^\inn$). Both these spaces have
{\sl fine moduli\/} (\S\ref{intSNorm} and Rem. \ref{stack}). 

That is, $\Psi_\Phi$ determines $\Phi$  up to equivalence on families.
Fine absolute (resp.~inner) moduli  holds because 
$A_n(1)$ is its own normalizer in
$A_n$ (resp.~$A_n$ has no center; \S \ref{absEquiv}). Analogous results hold for reduced
equivalence (Rem.~\ref{intSNorm}). 

In practice fine moduli means we know all families in a Nielsen class from knowing 
properties of one family over the corresponding space $\sH(G,\bfC)^*$. Further, properties of that family
come from an action of a {\sl mapping class group\/} on Nielsen classes. For inner and absolute equivalence
the mapping class group  is the {\sl Hurwitz monodromy\/} group $H_r$, a braid group quotient. 

With  $\ni(A_n,\bfC)^\inn=\ni(A_n,\bfC)/A_n$ and $\ni(A_n,\bfC)^\abs=\ni(A_n,\bfC)/S_n$: 
\begin{itemize}
\item \S \ref{BHMon} reminds that $H_r$ orbits on $\ni(A_n,\bfC)^\inn$ (resp.~$\ni(A_n,\bfC)^\abs$) 
$\Leftrightarrow$ 
$\sH(A_n,\bfC)^\inn$ (resp.~$\sH(A_n,\bfC)^\abs$) components .
\end{itemize} 
Lem.~\ref{lnf.5} says $H_r$ orbits on $\ni(G,\bfC)^\inn$ correspond one-one with 
orbits on $\ni(G,\bfC)$, but depending on $(G,\bfC)$ this may not hold with \lq\lq$\abs$\rq\rq\ replacing
\lq\lq$\inn$\rq\rq\ (Ex.~\ref{mod8exs}).

Our M(ain) R(esult) is
Thms.~\ref{thmA} and \ref{thmB}: Listing absolute/inner components for 
$(A_n,\bfC=\bfC_{3^r})$, $n\ge 4$, $r\ge n-1$.  The MR proof takes up
\S\ref{suppLems},
\S\ref{ns.2}, \S\ref{nt} and \S\ref{nfo}.  \S\ref{nth} and \S\ref{halfclasses} tie to \cite{Fr2},
\cite{Me}, \cite{Se3} and
\cite{AP} for a list of  applications. 

\subsubsection{Lifting invariants} \label{nz.2} A Frattini cover $G'\to G$ is a group cover (surjection)  
where restriction to any proper subgroup of $G'$ is not a cover. I now explain how any central Frattini
extension 
$\psi: R\to G$ gives a 
lifting invariant \cite[Part II]{Fr4}. 

A special case comes from 
alternating groups.    Let $\Spin^+_n$ be the (unique) nonsplit degree 2 cover of the connected
component $O_n^+$ (of the identity) of the 
orthogonal group (\cite[\S II.C]{Fr4} or 
\cite{Se2}). Regard $S_n$ as a
subgroup of the orthogonal group $O_n$. The alternating group $A_n$
is in  $O_n^+$, the kernel of the determinant map. Denote its pullback to $\Spin_n^+$ by
$\Spin_n$ and  identify
$\ker(\Spin_n\to A_n)$ with the multiplicative group $\{\pm 1\}$. Odd order elements of
$S_n$ are in
$A_n$.  Any odd order $g\in A_n$ has a unique odd order lift, $\hat g\in \Spin_n$. 

Let
$\bg\in A_n^r$, with
$g_1\cdots g_r=\Pi(\bg)=1$.  
If entries of $\bg$ have odd order, 
define the {\sl spin\/} lifting invariant of $\bg$ to be  
\begin{equation} \label{spinlift} s(\bg)=s_{\Spin_n}(\bg)=\hat g_1\cdots \hat
g_r\in \{\pm 1\}.\end{equation}  A degree $n$ (absolute) cover has the same lifting invariant as
its Galois closure. So, s(\bg) will not distinguish between absolute and inner
classes. 

\subsection{Main Result and its corollaries} \label{nz.4} \S\ref{An3r} states and outlines the proof of the
MR. Then, the three remaining subsections discuss the applications. 

\subsubsection{$\Spin_n\to A_n$ and 3-cycles} \label{An3r}
Strong Coalescing Lem.~\ref{lnf.2'} applies 
for $n\ge 5$. It says, for $\bg\in 
\ni(A_n,\bfC_{3^r})$, there is $q\in H_r$  so  $(\bg)q=\bg'$ has 
$(g'_2)^{-1}=g_1'$. We then   
induct on $(r,n)$ to describe all components of
$\sH(A_n, \bfC_{3^r})^*$, with $*=\abs$ or $\inn$. 

If a curve cover $\phi:X\to \prP^1_z$
in
$\ni(A_n,\bfC_{3^r})$  corresponds to a point of a component labeled 
$\oplus$  (resp.,  $\ominus$) in the {\sl Constellation diagram\/} (Fig.~1), then any branch cycle
description $\bg$ for $\phi$  has  
$s(\bg)=+1$ (resp.~-1; \S\ref{nz.2}). Lifting invariants are the same for covers
representing points on the same component. Fig.~1 labels each
component at   $(n,r)$ with a symbol $\oplus$ or $\ominus$ corresponding to the lifting invariant value.
Thm.~\ref{thmA}   generalizes \cite[Cor.~\ref{nf.2}]{Se2} in  
showing there is {\sl exactly\/} one component when the curves are genus 0: row-tag $\mapright{g=0}$ in
Fig.~1. 

\begin{thm} \label{thmA} For  $r=n-1$, $n\ge 5$,  ${\sHt n {n\nm1}}^\inn$
has exactly one connected component.  
Further, $\Psi^\inn_\abs: {\sHt n {n\nm1}}^\inn\to 
{\sHt n {n\nm1}}^\abs$ has
degree 2. \end{thm}

The row 
with tag $\mapright{g\ge 1}$ illustrates this theorem. 

\begin{thm} \label{thmB} For each $r\ge n\ge5$, ${\sHt 
nr}^\inn$ has exactly two connected components,
${\sHtp nr}^\inn$ (symbol $\oplus$) and ${\sHtm nr}^\inn$ (symbol 
$\ominus$). Denote their respective (connected) images in ${\sHt nr}^\abs$ by 
$\sH_{\pm}(A_n,\bfC_{3^r})^\abs$. The maps $$\Psi^{\inn,\pm}_\abs:
\sH_{\pm}(A_n,\bfC_{3^r})^\inn\to
\sH_{\pm}(A_n,\bfC_{3^r})^\abs \text{ \rm have degree 2}.$$
\end{thm} 

\newcommand{\icona}{\ominus}
\newcommand{\iconb}{\oplus}
\newcommand{\iconc}{}
\newcommand{\icon}{{\ominus\, \oplus}} 

\begin{figure}[h]  \label{An3rConst}
\caption{Constellation of spaces $\sH(A_n,\bfC_{3^r})^*$} 
\begin{tabular}{|c|c|c|c|c|c|c|} \hline 
&&&&&&\\ 
$\mapright{g\ge 1}$& $\icon$&$\icon$&
\dots & $\icon$& $\icon$ &$\mapleft{1\le g}$ 
\\ \hline &&&&&&\\ 
$\mapright{g=0}$ & 
$\icona$ & $\iconb$ & \dots & $\icona$ & $\iconb$ & $\mapleft{0=g}$ 
\\ \hline &&&&&&\\
$n\ge 4$ & $n=4$ & $n=5$ & \dots & $n$ even & $n$ odd & $4\le n$\\
\hline
\end{tabular} \end{figure}

We include a column for $n=4$, though 
$A_4$ has   two conjugacy classes of 3-cycles, with representatives $(1\,2\,3)$ and 
$(3\,2\,1)$. Denote the 1st by $\C_{+3}$, the 2nd by $\C_{-3}$. Then, 
$\bfC_{\pm 3^{s_1,s_2}}$ indicates $s_1$ (resp.~$s_2$) repetitions of $\C_{+3}$ (resp.~$\C_{-3}$); abbreviate 
to
$
\C_{\pm 3^{s_1}}$ if $s_1=s_2$. Expression \eqref{dagger}  must hold for $\ni(A_4,\bfC_{\pm 3^{s_1,s_2}})$ to
be nonempty. 
\begin{exmp} \label{pm3^2} As in \eql{ent.1}{ent.1b}, $\ni(A_4,\bfC_{\pm 3^2})$ has two braid orbits,
with  reps.:  
$$\bg_{4,+}=((1\,3\,4),(1\,4\,3),(1\,2\,3),(1\,3\,2)) \text{ and }
\bg_{4,-}=((1\,2\,3),(1\,3\,4),(1\,2\,4),(1\,2\,4)).$$
\end{exmp} 

Trickiest point: if $s_1\ne s_2$, then $\ni(A_4,\bfC_{\pm3^{s_1,s_2}})$ and
$\ni(A_4,\bfC_{\pm3^{s_2,s_1}})$ are distinct. Yet,  any  $\beta\in S_4\setminus A_4$ (as in
\S\ref{nt.2}) conjugates between them. So
$$\ni(A_4,\bfC_{\pm3^{s_1,s_2}})^\abs=\ni(A_4,\bfC_{\pm3^{s_2,s_1}})^\abs.$$  So, if for
$n=4$, you put any one allowed value of $(s_1,s_2)$ for each $r$, then Figure 1 still is valid. Examples: For
$r=3$,
$\{s_1,s_2\}=\{0,3\}$ give  one $\abs$ component; for $r=4$, $\{s_1,s_2\}=\{2,2\}$ and the two
components (*=abs or in) are both over $\bQ$; for $r=5$, $\{s_1,s_2\}=\{1,4\}$ there are 
two (*=abs or in)  components ($\abs$ over $\bQ$, $\inn$ over $\bQ(\sqrt{-2})$). For
$r\ge 6$ there are several  values of 
$\{s_1,s_2\}$. 

\begin{exmp}[$\inn$ is stronger than $\abs$] \label{mod8exs} Let $\C_d$ be the class of $d$-cycles
in
$A_n$. For
$n\equiv 1\mod 4$ and $d(n)=\frac{n\np1}2$ consider the Nielsen class $\ni(A_n,\bfC_{d(n)^4})$ of 4
reps.~of $d(n)$-cycles. A special case of \cite[Thm.~5.5]{osserman} gives one braid orbit on
$\ni(A_n,\bfC_{d(n)^4})^\abs$, suggesting there are no distinguishing properties of these spaces as $n$
varies. Yet,
\cite[Prop.~5.15]{twoorbit} starts their many differences with this: 
\begin{edesc} \label{n14} \item \label{n14a} for $n\equiv 1\mod 8$, $\ni(A_n,\bfC_{d(n)^4})^\inn$ has two
braid orbits, with the corresponding spaces conjugate over a quadratic extension of $\bQ$; while 
\item \label{n14b} for $n\equiv 5 \mod 8$, $\ni(A_n,\bfC_{d(n)^4})^\inn$ has but one braid orbit. \end{edesc} 
\end{exmp}   
     
\subsubsection{Corollaries on $A_n$ realizations}    Let $\sH'$ be 
an (irreducible) component of  $\sH(A_n, \bfC_{3^r})^\inn$. 
\!\!Cor.~\ref{corAr1} says   
$\sH'$ (and its map to $U_r$) has definition field  $\bQ$. 

Any $\bp'\in \sH'$ produces
a field extension $\hat L_{\bp'}/\bQ(\bp')(z)$ regular and Galois over $\bQ(\bp')$, with group $A_n$. 
Let $\sH$ be the image of $\sH'$ in 
$\sH(G,\bfC)^{\abs}$. 
Then, any $\bp\in \sH$ produces
a regular degree $n$ field extension $L_\bp/\bQ(\bp)(z)$. 
This is the natural extension of function fields for a cover
$\phi_\bp: X_\bp\to \prP^1_z$ representing $\bp$. The geometric (resp.~arithmetic)
Galois closure of $L_\bp/\bQ(\bp)(z)$ has group $G_\bp=A_n$ (resp.~$\hat G_\bp$ between $A_n$ and $S_n$).
So, 
$\hat G_\bp$ is either $A_n$ or $S_n$: respectively,  each $\bp\in \sH$ produces  either
an $(A_n,A_n)$  or an $(A_n,S_n)$ realization over $\bQ(\bp)$.

Cor.~\ref{corAr1} produces  a dense
set of $\bp\in \sH(\bar \bQ)$ with $(G_\bp,\hat G_\bp)$ equal to $(A_n,S_n)$ (resp.~equal to $(A_n,A_n)$).
It is subtler to ask if either conclusion holds restricting to $\bp\in \sH(K)$ for  $[K:\bQ]< \infty$. 
Combining Cor.~\ref{corAr2} and
\cite{Me} shows the $(A_n,A_n)$ conclusion if $r=n-1$, even for $K=\bQ$. 

\!Let $\bfC$ be {\sl any\/} odd order classes in $A_n$. As \cite{twoorbit} shows, the algorithm behind the
MR is applicable to much more than the case of 3-cycles. 

Example: It makes sense to speak of -- and to identify
-- the
$+$-components, $\sH_+(A_n,\bfC)^*$ (as in \eql{n14}{n14a}  with $*=\inn$, there may be more than one), among
all components of
$\sH(A_n,\bfC)^*$. 
\cite[Prop.~6.8]{BaFr} then interprets as follows.
\begin{lem} \label{ratPtSpin}  Assume $\bp\in
\sH_+(A_n,\bfC)^\abs(K)$, corresponding to  $\phi_\bp:
X_{\bp}\to
\prP^1_z$. With $\bz_\bp$ its branch points, assume there is $z_0\in
\prP^1_z\setminus
\bz$ with all points of $\phi_\bp^{-1}(z_0)$ in $K$. Then, 
$\bp$ gives an
$(A_n,A_n)$ realization and   
$\bp'\in \sH(\Spin_n,\bfC)^\inn=\sH_+(A_n,\bfC))^{\inn}$ over it gives a $\Spin_n$ regular $K$ realization.
(Rem.~\ref{stack} for why this is nonobvious.) 
\end{lem} 

\subsubsection{$G_\bQ$ and canonical fields from alternating groups}
\label{QAltPAC} \newcommand{\alt}{\text{\rm alt}}

Denote the absolute Galois group of a field $K$ by $G_K$.  
One goal of arithmetic geometry is to present $G_\bQ$ as a
known group quotient
$N=G(F/\bQ)$ (so $F$ is Galois over $\bQ$) by a known subgroup $M=G_F$. \fvb\ produced such
presentations with $N$ a product of symmetric groups and $M=\tilde F_\omega$, the profree group on a
countable number of generators. Fields $F$ with known arithmetic properties enhance 
applications. 

The archetype is  {\sl Shafarevic's conjecture\/}: For $F=\bQ^\cyc$, $N=\hat \bZ^*$ (profinite 
invertible integers) and $M=\tilde F_\omega$. \S\ref{hpconj} explains the mystery of whether 
$\bQ^\alt$, the composite of all $A_n$ extensions of $\bQ$, $n\ge 5$,  should have $G_{\bQ^\alt}=\tilde
F_\omega$.  Our MR makes this plausible
using
\cite[Thm.~A]{FrV3}: For
$F\subset
\bar\bQ$  P(seudo)A(lgebraically)C(losed), 
$G_F$  is Hilbertian if and only if it is $\tilde F_\omega$. In particular, if  
$\bQ^\alt$ is  PAC,  then $G_{\bQ^\alt}=\tilde
F_\omega$.

The following sequence makes a case for  $\bQ^\alt$ being PAC.  First, 
it is PAC  if  each $\bQ$ curve $X$ has a $\bQ$ cover  
$X\to
\prP^1_z$ of degree $n$ giving  an $(A_n,A_n)$ realization  (\S\ref{secrigp-aigp}) for some $n\ge 5$
($n$ allowed to vary with $X$). Every curve of genus
$g$ appears as a geometric $A_n$ cover with odd order branching, for many possible degrees, 
from  
\cite{AP} (\S\ref{maxAltExt}). Further, most curves give a corresponding geometric point on $\infty$-ly many
of the spaces
$\sH(A_n,\bfC_{3^r})^{\inn}$. If this implies each curve over
$\bQ$ gives a corresponding $\bQ$ point on one of these spaces, then we have the result. 

That 
these spaces have definition field  $\bQ$ (Thm.~\ref{thmB})  further encourages. 
Still, Prop.~\ref{getAnAn} shows 
there are  many $X\,$s over $\bQ$, of each positive even genus, with  {\sl no\/} odd branched $\bQ$  cover
$X\to\prP^1_z$ of any degree, much less any 
$(A_n,A_n)$ realization. 

A serious issue around  Hilbert's Irreducibility Theorem (HIT) arises in 
\S\ref{hitprop}.  
\begin{triv} \label{Anissue}  Extending  \cite[p.~36--37]{Mu76} to ask
if
$\sM_g$ is an {\sl Hilbertian variety}; why showing  HIT for the \hc\ covers $\sM_{g,{\pm}}\to
\sM_g$ is nontrivial. 
\end{triv}

\subsubsection{Spaces supporting $\theta$-nulls} \label{thetaNull} For each $(n,r)$, $r\ge n$, 
there is a map from $\sH_{\pm} (A_n,\bfC_{3^r})^\abs$ to $\sM_{g,\pm}$, the space of
genus $g=r-(n-1)$ curves with an attached half-canonical class (\S\ref{hcfamilies}).  From 
Thm.~\ref{thetaNullRes} the reduced versions (\S\ref{rnz.4}) of $\sH_{+} (A_n,\bfC_{3^r})^\inn$  support a canonical 
even $\theta$-null (with 2-division characteristic) $\theta_{n,r}[0]$. For the absolute spaces, there is such a $\theta$-null on $\sH_{+} (A_n,\bfC_{3^r})^\abs$ (resp. $\sH_{-} (A_n,\bfC_{3^r})^\abs$) if $r$ is even (resp.~odd). A  power of this
is the {\sl Hurwitz-Torelli\/} analog of an automorphic function.  
We only, however, know it is non-zero for absolute
spaces and, for given $g$, infinitely many explicit $(n,r)$ (including $g=1$). 

\section{Coalescing and supporting lemmas} \label{suppLems} \!\!\!This section shows how to braid $\bg\in
\ni(A_n,\bfC_{3^r})$ to where its first 2 or 3 entries are in a list of precise possibities. 
\cite{Se2} is a quick corollary. 

\subsection{Braid and Hurwitz monodromy groups} \label{BHMon} Generators 
$\row q {r-1}$ of the {\sl Hurwitz monodromy group}, $H_r$, a quotient of the braid group $B_r$, act
as permutations on   the right of $\ni(G,\bfC)$.   For $\bg\in \ni(G,\bfC)$: 
$$ (\bg)q_i=(g_1^\sph,\dots,g_{i-1}^\sph,g_i^\sph g_{i+1}^\sph
g_i^{-1}, g_i^\sph, g_{i+1}^\sph,\dots, g_r).$$ Generators $\row Q {r-1}$ of $B_r$ 
generate it freely modulo these relations: 
\begin{triv} \label{Brrel} 
$Q_iQ_j=Q_jQ_i$ for $|i-j|> 1$ and 
$Q_{i\np1}Q_iQ_{i\np1}=Q_iQ_{i\np1}Q_i$. 
\end{triv}
Add to \eqref{Brrel} one further  relation for $H_r$: \begin{triv}  \label{Hurwitz}
$q_1\cdots q_{r-1}q_{r-1}\cdots q_1=1$.
\end{triv} 
Also,  $H_r$ is the fundamental group of 
projective $r$-space 
minus the discriminant locus: $\prP^r\setminus D_r$; that is, the space of monic polynomials of degree $r$
with no repeated roots. Another description of $\prP^r\setminus D_r$  is as the quotient of $(\prP^1)^r
\setminus \Delta_r/S_r$. 

\noindent The word $Q_1\cdots Q_{r-1}Q_{r-1}\cdots Q_1\eqdef Q^{(r-1)}\in B_r$ conjugates $\bg=(\row g r)$  
by $g_1$: 
$$\bg\mapsto (\bg)Q^{(r-1)}=g_1^\sph\bg g_1^{-1}=(\dots, g_1^\sph g_ig_1^{-1},\dots).$$
So, to have $H_r$ acting on $\ni(G,\bfC)$ requires quotienting by $G$: $g\in G$ has
the effect 
$$\bg\mapsto  g^{-1}\bg g=(g^{-1}g_1g,\dots, g^{-1}g_rg)\in \ni(G,\bfC).$$ The resulting
set $\ni(G,\bfC)/G=\ni(G,\bfC)^\inn$ we call {\sl inner Nielsen classes}. Also, the element of
$q_{(r-1)}\eqdef q_1\cdots q_{r-1}=\sh$ acts as a {\sl shift operator\/} on
$\bg\in\ni(G,\bfC)$: 
\begin{equation} \label{shift} (\bg)q_{(r-1)}=(g_2,g_3,\dots,g_{r},g_1).\end{equation} The word of 
\eqref{Hurwitz} acts trivially on  $\ni(G,\bfC)^\inn$. 

Computations in the first four sections
are for $B_r$ acting on  Nielsen classes. Sometimes (as in the proof of Lem.~\ref{3-3braids})  we extend that
action to  a generalization of Nielsen classes. These  are Nielsen {\sl sets\/} $\ni(G,\bfC)_{g'}$, defined
by $(G,\bfC,g')$ where we replace the product-one condition by 
$\Pi(\bg)=g'$ with $g'\in G$.  Only 
elements of
$G$ centralizing
$g'$ can act by conjugation on $\ni(G,\bfC)_{g'}$. Corollaries, however, then pass to 
$H_r$ acting on Nielsen classes. App.~\ref{hurSpaces} translates 
between Nielsen classes and these spaces. This dictionary reduces \S2-\S4 to combinatorics
and group theory. 

There is a homomorphism 
$\alpha: B_r\to S_r$ by  $Q_i\mapsto \alpha(Q_i)=(i\, i\np1)$. For $\bh\in 
G^u$, denote 
juxtaposition of $k$ copies of $\bh$ by $\bh^{(k)}$. For example, 
$$((1\,2\,3),(1\,3\,2), (1\,2\,3),(1\,3\,2))=((1\,2\,3),(1\,3\,2))^{(2)}.$$

Suppose $Q\in B_r$. Call $Q$
{\sl local\/} to a subset $I$ of
$\{1,\dots,r\}$ if $Q$ is a product of braids  affecting only the
positions in $I$. Further, suppose $\bg\in G^r$ and  $i$ and $j$ are integers, 
with $1\le i< j\le r$. Then,  there exists $Q\in B_r$ local to the integers 
between 
$i$ and $j$ (inclusive) with $(\bg)Q=\bg'$ and $g_i'=g_j$.

\subsection{Coalescing} \label{nf.1}
The first part of Lem.~\ref{lnf.1} simplifies working with alternating groups. 
(Regard $A_1=\{1\}$ as the degree 1 alternating group.) \S\ref{nz.2} has the definition of the
lifting invariant 
$s(\bg)$ for $\bg$ having odd order entries.

\subsubsection{A starter lemma} Lem.~\ref{lnf.2} uses 
specific braids we regard as standard. 

\begin{lem}[3-cycle Lemma] \label{lnf.1} Let $\bg\in \Ct r$.
Let $G=\lrang {\bg}$ act on $\{1,\dots,n\}$.
Then, $G$ is a product of alternating groups, one copy for each orbit of
$G$.  Up to conjugacy in $S_n$, here are all 3-cycle pairs  with  
product 
a power of a 3-cycle.
\begin{edesc} \label{A3list1}
\item \label{(i)} $(g,g^{-1})\mapsto 1$.
\item \label{(ii)} $((i\,j\,k),(i\,k\,t))\mapsto (i\,j\,t)$, for $k\ne t$.
\item \label{(iii)} $(g,g)\mapsto g^{-1}$.\end{edesc}
Up to conjugacy in $S_n$, here are all 3-cycle pairs with product   not a 3-cycle power. 
\begin{edesc}  \label{A3list2}
\item \label{(iv)} $(g,g')\mapsto (g)(g')$ where $g$ and $g'$ have no common support.
\item \label{(v)} $((i\,j\,k),(i\,j\,t))\mapsto (j\,k)(i\,t)$, for $k\ne t$.
\item \label{(vi)} $((i\,j\,k),(i\,l\,m))\mapsto (i\,j\,k\,l\,m)$ with
$\{i,j,k,l,m\}$ distinct integers.\end{edesc} 
Let $(\row{\bg'} t)=\bg$ with $\Pi(\bg_i')=1$, $i=1,\dots,t$.
Then, $\prod_{i=1}^ts(\bg'_i)=s(\bg)$. If $g$ is a 3-cycle, then
$s(g^{(3)})=s(g,g^{-1})=1$. Finally, $s((i\,j\,k),(i\,k\,t),(i\,t\,j))=-1$.
\end{lem}

\begin{proof} 
Assume $\bg$ generates a transitive group. Then, the 
first statement says $\lrang{\bg}=A_n$. This is well-known 
from the following chain of deductions: $\lrang{\bg}$ is primitive, and a 
primitive group containing a 3-cycle is the alternating or symmetric group. 
If $\bg$ isn't transitive, then each 3-cycle has support on one 
of the orbits. Thus, you can apply the first argument to the 3-cycles 
supported on each orbit of $\lrang{\bg}$. 

Everything else is elementary. Example: Let $\hat g$ be the (unique) order 3 lift to $\Spin_n$
 of the 3-cycle $g\in A_n$. Then,
$s(g^{(3)})=\hat g^3=1$ and $s(g,g^{-1})=\hat g\hat g^{-1}=1$. Note that $s(h\bg h^{-1})=
\hat h s(\bg) \hat h^{-1}$ for any $h\in 
A_n$ and any lift of $h$ to $\hat h\in \Spin_n$.
Thus,
assume $((i\,j\,k),(i\,k,t),(i\,t\,j))=\bg$ is a 3-tuple in $A_5$: 
$A_5$ acts on the first five integers from $\{1,\dots,n\}$. Then,  $s(\bg)$
doesn't depend on whether we see $\bg$ as elements in $A_5$ or in $A_n$,
$n\ge5$. Identify $A_5$ with
$\PSL_2(\bZ/5)$ and $\Spin_5$ with $\SL_2(\bZ/5)$. Thus, the final
calculation is an explicit computation with $2\times2$ matrices. This
appears in Part C of the proof of \cite[Ex.~3.13]{Fr4} or in \cite{Se2}.
\end{proof}

\subsubsection{Disappearing sequences} The cases of \eqref{A3list1} separate according to the
conjugacy class in $S_n$ of the product of the three pairs. The phrase {\sl coalescing types\/} 
refers to this. Below we add to  
coalescing types \eql{A3list1}{(ii)} and \eql{A3list1}{(iii)} the possibility of $\bg$ having as its first 3
or 4 entries these tuples  (up to conjugation)  having product-one:  
\begin{edesc} \label{A3list3} \item \label{(ii')} $((1\,2\,3),(1\,3\,4),(1\,4\,2))$; 
\item \label{(ii'')} $((1\,2\,3),(1\,3\,4),(1\,2\,4),(1\,2\,4))$; and 
\item \label{(iii')} $((1\,2\,3)^{(3)})$. 
\end{edesc} Then, \eql{A3list3}{(ii')} and \eql{A3list3}{(ii'')} (resp.~\eql{A3list3}{(iii')}) correspond to
\eql{A3list1}{(ii)} (resp.~\eql{A3list1}{(iii)}).  Only when 
$n=3$ or 4 are we forced to use \eql{A3list3}{(ii'')}
as a braiding target (see \S\ref{nt.2}).  
  
 Recall the  homomorphism 
$\alpha$ of \S\ref{nz.1}. 
Denote the subgroup of $Q\in B_r$ with $\alpha(Q)$ permuting 
$\{1\,\dots,k\}$ 
by $B_r^{(k)}$. 

For any $\bg\in \ni(G,\bfC)$, use $g_1\cdots g_r=1$ to draw the following conclusion.
For each $i\in \N n$ there
exists $1\le j_1<\cdots< j_k\le r$ with these properties.
\begin{edesc} \label{enf.1}
\item \label{enf.1a} The support of $g_{j_1}$ contains $i$.
\item \label{enf.1b} $g_{j_1} g_{j_2}\cdots g_{j_k}$ fixes $i$.
\end{edesc}
Call the sequence $\row j k$ a {\sl disappearing\/} sequence for $i$. For $G$ transitive, a  braid to
a conjugate of the $r$-tuple (say, by Lem.~\ref{lnf.5}),  replaces $i$ by any desired integer. 

\begin{lem}[Coalescing]\label{lnf.2} For  
$\bg\in \ni(A_n,\Ct r)$, $n\ge 3, r\ge n-1$, there is a $Q\in B_r$ with
$(\bg)Q=(\row {g'} r)$ and $(g'_1,g'_2)$ is a disappearing sequence of length 2, (coalescing type
\eql{A3list1}{(i)} or \eql{A3list1}{(ii)}) or $(\bg)Q$ has coalescing type \eql{A3list1}{(iii)}.   

Stronger still, if the first 3 terms of 
$\bh\in
\ni(A_n,\Ct r)$ are a type \eql{A3list1}{(iii)} disappearing sequence,  then either
$(h_1,h_2,h_3)=(h_1^{(3)})$ (type \eql{A3list3}{(iii')}), or the first two terms of $(\bh)Q_2^2$ are a
disappearing sequence. 
\end{lem}

\begin{proof} Suppose we find $Q\in B_r$ with $(\bg)Q=\bg'$ and $(g'_i,g'_j)$
($i< j$) a pair of coalescing type from \eqref{A3list1}. Let $Q'=(Q_1\cdots Q_{i-1})^{-1}(Q_2\cdots
Q_{j-1})^{-1}$. Then, $((\bg)Q)Q'=\bg''$ has  
$g_1''=g'_i$ and $g_2''=g'_j$.

Apply an element of $B_r$ to $\bg$ to assume a given
disappearing sequence for $i$ is $1,\dots, l$. For example, to
put $g_{j_1}$ in the first position, apply
$Q_1^{-1}\cdots Q_{j_1\nm1}^{-1}$. To simplify notation, assume
$i=1$. Such a braiding moves $g_{j_1},\dots,g_{j_l}$. Still, it 
leaves them  in the same order they originally appear 
(reading left to right). 

If $l$ is two, Lem.~\ref{lnf.1} lets us take
$(g'_1,g'_2)$ to be one of \eql{A3list1}{(i)} or (with $k=1$) \eql{A3list1}{(ii)}. So, assume $l>2$. One
further  assumption: For all $Q\in B_r^{(k)}$, \begin{triv} \label{enf.1c}  $l$ is the  
shortest length of a disappearing sequence for 1 in $(\bg)Q$. \end{triv} \noindent This 
assumption lets us
prove $l=3$ and we may assume $(g_1',g_2',g_3')$ satisfies \eql{A3list3}{(iii')}. 

A disappearing sequence corresponds to integers in a chain: $$1\mapsto 
i_1\mapsto 
i_2\mapsto \cdots\mapsto i_{l\nm1}\mapsto 1,$$ where the 1st 3-cycle maps $1\mapsto i_1$ and the last
($l$th) maps $i_{l-1}\mapsto 1$.  Suppose this disappearing  sequence for 1 has a 3-cycle, 
say 
$g_u$,  not contributing 
to this chain. This violates \eqref{enf.1c} with $Q=1$. So, too,  none of $\row i {k\nm1}$ is 1. So, the
disappearing sequence has  the 
form \begin{equation} \label{enf.2} ((1\,i_1\,t_1),(i_1\,i_2\,t_2),\dots, 
(i_{l\nm2}\,i_{l\nm1}\,t_{l\nm1}),(i_{l-1}\,1\,t_l)).\end{equation} 

Now suppose in \eqref{enf.2}, $t_1\ne i_2$. Apply $Q_1$ to \eqref{enf.2} to get this 
replacement for the first two positions: 
$$((1\,i_1\,t_1)(i_1\,i_2\,t_2)(t_1\,i_1\,1),(1\,i_1\,t_1))= 
((1\,i_2\,\cdot),(1\,i_1\,t_1)).$$ Thus, the 1st  and   the 3rd through $k$th  
positions give a shorter disappearing 
sequence for 1. This 
violates \eqref{enf.2} with $Q=Q_1$.  Conclude $t_1=i_2$. Similarly, applying $Q_j$ 
to \eqref{enf.2} forces 
$t_j=i_{j\np1}$, $j=1,\dots,l\nm1$. Note: The last of these gives 
$t_{l\nm1}=1$. 
If $i_3=1$, then the first three positions contain the 3-cycles  
\begin{equation} \label{**} ((1\,i_1\,i_2),(i_1\,i_2\,1),(i_2\,1\,i_4)).\end{equation} This is 
\eql{A3list1}{(iii)},  a disappearing sequence of length 3 for 1. From \eqref{enf.1c} we are done, unless 
$i_3\ne 1$. In this case, apply $Q_1^{-1}$ to \eqref{enf.2}. Now the first two positions are 
 $$((i_1\,i_2\,i_3),(i_3\,i_2\,i_1)(1\,i_1\,i_2)(i_1\,i_2\,i_3))= 
((i_1\,i_2\,i_3), (1\,i_2\,i_3)).$$ The 3-cycles in the 2nd--$l$th 
positions 
give a length $l-1$ disappearing sequence for 1. This contradicts 
\eqref{enf.1c}. Conclude the first paragraph by inducting on $l$. 

Consider the hypotheses of the 2nd paragraph statement. In \eqref{**}, if $i_4= i_1$, then
$(h_1,h_2,h_3)=(h_1^{(3)})$. Otherwise,   applying $Q_2^2$ gives  the desired conclusion: 
$((1\,i_1\,i_2),(i_2\,i_4\,i_1),(i_1\,i_2\,1))Q_2=
((1\,i_1\,i_2),(i_1\,1\,i_4) ,(i_2\,i_4\,i_1))$.   
\end{proof}

\subsection{Invariance Corollary} \label{nf.2}
Cor.~\ref{cnf.3} reproves \cite{Se2}. For $g\in A_n$ of odd order, let $w(g)$ by the sum of
$\frac{\ell^2-1}8\mod 2$ over the length $\ell$ of 
  disjoint cycles in $g$. 

\begin{cor}[Invariance] \label{cnf.3} Let $n\ge 3$. 
If $\phi: X\to\prP^1$ is in the Nielsen class $\ni(A_n,\Ct {n\nm1})$, then   
$\deg(\phi)=n$,  $X$ has genus 0, and $s(\phi)=(-1)^{n\nm1}$.  

Generally, for any genus 0 Nielsen class of odd order elements, and representing $\bg=(\row g r)$, $s(\bg)$ is
constant, equal to $(-1)^{\sum_{i=1}^r w(g_i)}$. 
\end{cor}

\begin{proof} Induct on n. Apply 
Lem.~\ref{lnf.2} to  $\bg\in\ni(A_n,\Ct {n\nm1})$. Coalesce $g_1$ and $g_2$ in any of the cases (2.4a), (2.4b)
and (2.4c) to get $\bg'=(g_1g_2,g_3,\dots,g_{n-1})$. In each case of \eqref{A3list1}, $\bg'$ has 
fewer than $n\nm 1$ 3-cycles as entries (still with product 1). If $\lrang{\bg'}$ is transitive,  apply
RET to produce  $X'\to\prP^1$, a (connected) cover having $\bg'$ as its branch cycles. R-H
(in \eqref{RH}) implies $2(n+g(X')-1)=2r'$ with $r'=n\nm 2$ (case \eql{A3list1}{(ii)} or
\eql{A3list1}{(iii)}) or
$n\nm 3$ (case
\eql{A3list1}{(i)}). 
 This is a 
contradiction: The genus of  $X'$ would be 
negative. Conclude, $\lrang{\bg'}$ has more than one orbit  in each case. 

In case \eql{A3list1}{(iii)}, $g_1g_2=g_1^{- 1}=g_2^{-1}$: $\lrang{\bg'}$ has just one 
orbit. So, we can assume \eql{A3list1}{(i)} or \eql{A3list1}{(ii)}.  The formula is clear for $n=3$. Now do
an induction. 

Case \eql{A3list1}{(i)}: $\bg'$ has $n-3$ branch cycles, spread on
2 or 3 orbits. First assume $\lrang{\bg'}$ has orbits of
length $n_1$, $n_2$ and
$n_3$ ($n_1+n_2+n_3= n$). Thus, $\bg'$ has $n_i-1$ entries
supported on the $i$th orbit, $i=1, 2, 3$. Write $\bg'$ as
$(\bg_1,\bg_2,\bg_3)$, with the
3-cycles of $\bg_i$ having support on the $i$th orbit. According to Lem.~\ref{lnf.1},   
$s(\bg')=s(\bg)=\prod_{i=1}^3 s(\bg_i)$.
Apply the induction
assumption to conclude
$$s(\bg)= (-1)^{n_1-1}(-1)^{n_2-1}(-1)^{n_3-1}=(-1)^{n-1}.$$

Now we show  there cannot be
just 2 orbits. The orbit of length $n_i$ supports at least $n_i-1$
3-cycles. Thus, there must be at least $n_1-1+n_2-1=n-2$ of these 3-cycles.
There are, however, only $n-3$ of them.

Case \eql{A3list1}{(ii)}: 
Here $\bg'=(g_1g_2,g_3,\dots,g_r)$. Let $\hat g_1$,
$\hat g_2$ and $\widehat{g_1g_2}$ be respective lifts of
$g_1$, $g_2$ and $g_1g_2$ to $\Spin_n$. 
Lem.~\ref{lnf.1} gives $\hat g_1\hat g_2=-\widehat{g_1g_2}$.
Conclude: $s(\bg)=-s(\bg')$. In the product $g_1g_2$, exactly one integer
from the union of the support of $g_1$ and $g_2$ disappears. So,
$\bg'$ must have exactly two orbits of respective lengths
$n_1$ and $n_2$ with $n_1+n_2=n$. Apply the induction assumption exactly
as for \eql{A3list1}{(i)}. Thus:
$$s(\bg)= -s(\bg')=(-1)^{n_1-1}(-1)^{n_2-1}=(-1)^{n-1}.$$

Now for the general case where  $\bg$ has odd  order entries, but maybe not 3-cycles. Write 
$g_i$ as a product of disjoint cycles, $(g_{i,1},\dots, g_{i,k_i})=\bg_i$. Then,  juxtaposed $\bg_i\,$s
give $\bg^*=(\row {\bg} r)$ in a new Nielsen class, still  of genus 0. From
Lem.~\ref{lnf.1}, $s(\bg^*)=s(\bg)$. This reduces us to where all entries of $\bg$ are cycles.

To conclude, replace each $g_i$ (conjugate to $(1\,\dots\,k)$) by $(h_{i,1},\dots,
h_{i,k_i})=\bh_i$, conjugate to $((1\,2\,3),(1\,4\,5),\dots,(1\,k-1\,k))$. Call the juxtaposed branch cycle
$\bh$. The changes are canonical and only depend on the lengths of the disjoint cycles in $\bg$.
Apply Lem.~\ref{nf.1} to see  $s(\bg)/s(\bh)$ is $\prod_{i=1}^r u_i$ with  
$$ u_i=s((1\,\dots\,k)^{-1}, (1\,2\,3),(1\,4\,5),\dots,(1\,k-1\,k)).$$ Conclude easily from \cite[Lem.~2]{Se2}
: $u_i$ is  
$(-1)^{\frac{k^2\nm1}8+\frac{k-1}2}$ (see Rem.~\ref{cliffUse}). 
\end{proof}

\begin{rem}[Clifford algebra] \label{cliffUse}   The proof of
Cor.~\ref{cnf.3} uses the Clifford algebra only in computing $u_i$ from \cite[Lem.~2]{Se2}.  An induction
reduces this to computing directly $s((1\,\dots\,k)^{-1}, (1\,2\,3),(1\,4\,5\,\dots\,k))$, the
lifting invariant for a polynomial map one can write down by hand. Is there a simple proof this has value 
$(-1)^{(k^2+(k-2)^2-2)/8+1}$ without using the Clifford algebra?\end{rem}

\def\sbgp#1{{[\bg_{{#1},\bullet}}]} \def\sbgm#1{{\bg_{{#1},-}}}
\def\Sbgp#1 #2{{\bg_{{#1},+},((1\,2\,3),(3\,2\,1))^{({#2})}}}
\def\Sbgm#1 #2{{\bg_{{#1},-},((1\,2\,3),(3\,2\,1))^{({#2})}}}

\subsection{Product-one and H-M reps.} \label{nf.3} This section generalizes \cite[\S 3.F]{Fr4}.

\begin{defn}[H-M-Nielsen 
class generators] \label{dnf.4} For $r=2s$, 
let $\bfC$ be a collection of conjugacy classes from a group $G\le S_n$.
We don't assume $G$ is transitive.
Suppose $\bg\in\bfC$ has this form:
$(g_1^\sph,g_1^{-1},\dots, g_s^\sph,g_s^{-1})$.  
We say it is an H-M representative (H-M  rep.)  
of $\lrang{\bg}$. Also, $\bg$ is an {\sl
H-M rep.\/}
of $\ni(G,\bfC)$ if $\lrang{\bg}=G$. \end{defn}  \!The following is from  \cite[Lemma 3.8]{BFr}.

\begin{lem}[Product-one] \label{lnf.5}
 Let $\bg\in \bfC$ with $\Pi(\bg)=1$.
Let   $i,i\np1,\dots,j
\bmod r$ be consecutive integers with $g_i g_{i+1}\cdots g_j=1$ (including
$i,i\np1,\dots,r\nm1,r, 1, 2,\dots, j$). Let
$\gamma\in \lrang{ g_i,\dots,g_j}$. There is $Q\in B_r$ with  
\begin{equation} \label{enf.3} (\bg)Q=(g_1,\dots,g_{i-1}, \gamma g_i\gamma^{-
1}, \gamma g_{i\np1}\gamma^{-
1}\dots, \gamma g_j\gamma^{-1},g_{j+1},\dots,g_r).\end{equation} \end{lem}

Easily (as in the next lemmas) find
$Q\in B_r$ that takes an H-M rep.~to
\begin{equation} \label{enf.4} (g_1^\sph,\dots,g_{u}^\sph, g_{u}^{-1},\dots, g_1^{-1})=[g_1,\dots,g_{u}]
\eqdef [\bg].\end{equation} 

\begin{lem} \label{lnf.6} Take $\bg$ as in \eqref{enf.4}. For any
$\pi\in S_u$, there exists $Q\in B_r$ with
$$[\row g u]Q=[g_{(1)\pi},\dots,g_{(u)\pi}].$$
\end{lem}

\begin{proof} Transpositions generate $S_u$. It suffices to show
this when $\pi=(1\,2)$ with $Q\in B_r$ local to  the first four
entries.
Take $Q_{1,2}=
Q_1^{-1} Q_3^\sph $. Then:
\begin{equation} \label{enf.5} \begin{array}{rl}
(g_1^\sph,g_2^\sph,g_2^{-1}, g_1^{-1})Q_{1,2}&=(g_2^\sph,
g_2^{-1}g_1^\sph g_2^\sph,
g_2^{-1}, g_1^{-1})Q_3\\
&=(g_2^\sph, g_2^{-1}g_1^\sph g_2^\sph, g_2^{-1}g_1^{-1}g_2^\sph,
g_2^{-1}).\end{array} \end{equation} 
Product-one Lemma \ref{lnf.5} gives $Q'\in B_r$ conjugating \eqref{enf.5} by $g_2$ 
(fixing
the coordinates beyond the first four).
Conclude by taking $Q=Q_{1,2}Q'$. \end{proof}

\begin{lem}[Generator] \label{lnf.7} Assume the following  for 
$\bg=(\bg',\bg'')\in\bfC$:
\begin{triv} \label{enf.6}  $\bg'= [\row {g'} {u'}]$, and $\Pi(\bg'')=1$.
\end{triv}
\noindent Then, for any $h\in \lrang{\bg'}$, there is a $Q\in B_r$ with
$Q(\bg)=(\bg',h\bg''h^{-1})$. 

Suppose the following
holds for each $\bg=[\row g u]\in \ni(G,\bfC)$.
\begin{triv} \label{enf.7} For $\bg(i)=[g_1,\dots,g_{i-1},g_{i+1},\dots,g_u]$,
$\lrang{\bg(i)}=G$, $i=1,\dots,u$.
\end{triv}
\noindent Then, all H-M representatives of $\ni(G,\bfC)$ fall in one $B_r$ orbit.
\end{lem}

\begin{proof} We show the statement after
\eqref{enf.6}. Induct on the number of entries from $\bg'$  
to get the product $h$. So, it suffices to take $h=g_j'$ for
some $j$ between 1 and $u'$. Apply Lem.~\ref{lnf.6} to assume with no
loss $j=u'$. Then, apply Lem.~\ref{lnf.7} to
$([g_{u'}],\bg')$. This gives the conclusion to \eqref{enf.6}. The conclusion
following \eqref{enf.7} comes from repeated application of the above to
$(\bg(i)^\sph,g_i^\sph,g_i^{-1})$. \end{proof}

\!\begin{lem}[Blocks] \label{lnf.8} Suppose $\bg= (\row {\bg} u)$
with $\Pi(\bg_i)=1$, for all but possibly one  $i_0\in \{1,\dots, u\}$. 
For any $\pi\in S_u$, and $\tau_i\in \lrang{\bg_i}$, $i=1,\dots,u$, there exists $Q\in B_r$ with
$(\bg)Q=(\bg_{(1)\pi},\dots,\bg_{(u)\pi})$. Also,
for any $i$ and $j$ (excluding $j=i_0$ if $i_0$ exists), there exists $Q\in B_r$ with
$$(\bg)Q=(\bg_1,\dots,\bg_{j-1},\tau_i^\sph\bg_j\tau_i^{-1},\bg_{j+1},
\dots,\bg_u).$$
\end{lem}

\begin{proof} The case $u=2$ suffices to show we can permute the appearance of 
the
$\bg_i\,$s. For this, assume $\Pi(\bg_1)=1$, and braid
every entry of $\bg_1$, in order from left to right,
past every entry of $\bg_2$. This gives the effect of 
\begin{equation} \label{enf.8} (\bg_1,\bg_2)Q=(\alpha\bg_2\alpha^{-1},\bg_1), 
\alpha=\Pi(\bg_1).\end{equation} Done, since $\alpha=1$. The last sentence reduces to  cases
$i=j=1$ and $i=1$, $j=2$. 

For the 1st,  apply Lem.~\ref{lnf.5} to $\bg_1$. For the 2nd,  with 
$g$ an entry  of $\bg_2$ ($\Pi(\bg_2)$ may not be 1),  braid to $(g\bg_1g^{-1},\bg_2)$ 
with 
$\bg_2$ written 
$(\bh, g, \bh')$. Braid to $(\bh,\bg_1,g,\bh')$ as above. Then, braid the
sequence $\mapsto (\bh, g, g^{-1}\bg_1 g, \bh')\mapsto (g^{-1}\bg_1 g, \bg_2)$. 
\end{proof} 

\section{Coalescing targets} \label{ns.2} The induction goal, for given $n$ and
$n-1=r\ge 4$ (resp.~$r\ge n\ge 5$), is to apply a  $Q\in B_r$ to any $\bg\in\ni(A_n,\bfC_{3^r})$ so
 $(\bg)Q$ is an (resp.~one of two)  exemplar(s). \S\ref{nge5} lists 
coalescing targets for $n\ge 5$. Yet, these  require  
the intricate case $n=4$ (\S\ref{3-4targets}).   

\subsection{Coalescing Targets, $n\ge 5$} \label{nge5} We use  
$\sbgp u =[\row g u]$~(Def.~\ref{dnf.4})
with $g_1=(1\,2\,3),\ g_2=(1\,4\,5), \dots ,g_u=(1\,2u\,2u\np1)$ to list braid targets.  

\subsubsection{Normal forms for $\bfC_{3^r}$} \label{nformC3r} The pairs of Nielsen class 
representatives $\bg$ below have respective lifting invariants $s(\bg)=+1$, 
and $s(\bg)=-1$ (\S\ref{nz.2}). 

\begin{equation} \label{ens.1a}\begin{array}{rl}&r\ge n \text{ both odd}, (\sbgp {\frac{n-3}
{2}},(1\,n\nm1\,n)^{(3)}, ((1\, 2\, 3), (3\, 2\, 1))^{(\frac {r-n}{2})}),\\
&(\sbgp {\frac{n-3}{2}},
(1\,n\nm2\,n\nm1),(1\,n\nm1\,n),(1\,n\, n\nm2),
((1\, 2\, 3), (3\, 2\, 1))^{(\frac{r-n}{2})}).\end{array}\end{equation}
\begin{equation}\label{ens.1b} \begin{array} {rl}
&\text{odd }r \ge n \text{ even, }(\sbgp {\frac{n-2}{2}},(1\,n\nm1\,n)^{(3)},
((1\, 2\, 3), (3\, 2\, 1))^{(\frac{r-n-1}{2})}),\\
&(\sbgp {\frac{n-4}{2}},
(1\,n\nm2\,n\nm1),(1\,n\nm1\,n),(1\,n\, n\nm2),
((1\, 2\, 3), (3\, 2\, 1))^{(\frac{r-n+1}{2})}).\end{array} \end{equation}
\begin{equation} \label{ens.2a} \begin{array} {rl}&\text{even }r \ge n \text{ odd, }
(\sbgp {\frac{n-1}{2}},((1\, 2\, 3), (3\, 2\, 1))^{(\frac{r\nm n\np1}{2})}),\\ 
&(\sbgp {\frac{n-3}{2}},(1\,n\nm2\,n\nm1)^{(2)},  (1\,n\nm2\,n),(1\,n\,n\nm 1),
((1\, 2\, 3), (3\, 2\, 1))^{(\frac{r\nm n\nm 1}{2})}).\end{array}\end{equation}
\begin{equation} \label{ens.2b} \!\begin{array}{rl}&r\ge n \text{ both even}, (\sbgp {\frac{n-2}
{2}},(1\,n\nm1\,n), (1\,n\nm1\,n)^{\nm1}\!,((1\, 2\,  3), 
\!(3\, 2\, 1))^{{(\frac{r-n}{2})}}),\\ 
&(\sbgp {\frac{n-4}{2}},(1\,n\nm1\,n\nm2)^{(2)},(1\,n\nm1\,n),(1\,n\,n\nm 2),
((1\, 2\, 3), (3\, 2\, 1))^{(\frac{r\nm n}{2})}).\end{array}\end{equation}

\newcommand{\nb}{\text{\rm nb}}
\subsubsection{More on $\bfC_{3^r}$ normal forms} \label{nfComm} Each  $\bg$ in \S\ref{nformC3r}
ends with 
$((1\, 2\, 3), (3\, 2\,  1))^{(t)}$ for some  $t$. Denote this end part  $\bg_e$, 
and the beginning part 
$\bg_b$: $\bg=(\bg_b,\bg_e)$. We chose $\bg_b$  to be transitive  on 
$\{1,\dots,n\}$, with $s(\bg_b)=s(\bg)$. Refer to the spin lifting value by a subscript: as in
\eqref{ens.2b}$_{\pm}$ indicating the two \eqref{ens.2b} listings. 

\begin{defn} \label{nub} Each $\bg_b$ in \S\ref{nformC3r} starts with part of an (element braid equivalent to an) 
H-M rep.~$\sbgp
{\frac{n-u}{2}}$ (Def.~\ref{dnf.4}). The quirky part is $\bg_{\nb}$, what is left after the H-M in $\bg_b$. 
This is the {\sl nub\/}. \end{defn}
 
For example, respective $+$ and $-$ nubs  of
\eqref{ens.1a} are $$(1\,n\nm1\,n)^{(3)}\text{ and }(1\,n\nm2\,n\nm1),(1\,n\nm1\,n),(1\,n\, n\nm2).$$ We see
the value of $s(\bg)$ from the nub alone. 

The -1 rep.~(resp.~+1 rep.) of \eqref{ens.1b} (resp.~ \eqref{ens.2a}) also
works for $r=n-1$.

Strong Coalescing Lemma \ref{lnf.2'} gives the tools for braiding any $\bg\in \ni(A_n,\bfC_r)$ to where it has
the correct $\bg_e$. So, in the induction of 
\S \ref{nfo}, the significant braidings are where there is no $((1\, 2\, 3), (3\, 2\, 1))^{(t)}$ tail. 

For example,  the $\bg_b$ part of the -1 rep.~of \eqref{ens.2a} braids to 
\begin{equation} \label{ens.3a} (\sbgp
{\frac{n-5}{2}},(1\,n\nm3\,n\nm2)^{(3)},(1\,n\nm2\,n\nm1),(1\,n\nm1\,n),(1\,n\, n\nm2)).\end{equation}  Also,
the $\bg_b$ part of the 2nd element of \eqref{ens.2b} braids to   
\begin{equation} \label{ens.3b} (\sbgp {\frac{n-6}{2}},
(1\,n\nm4,n\nm3)^{(3)},(1\,n\nm2\,n\nm1),(1\,n\nm1\,n),(1\,n\, n\nm2)).\end{equation}
Respectively, these  are \eql{enf.11}{enf.11b} and \eql{enf.11}{enf.11a} in the proof of Lem.~\ref{lnf.2'}.
 
\newcommand{\bep}{{\pmb \epsilon}}

\subsection{Coalescing targets for  $n=3,4$} \label{3-4targets} While $n=3$ is easy,  
 $n=4$ is not. 
\subsubsection{Two conjugacy classes of 3-cycles} \label{2-3-cycles} 

The Klein 4-group
$K$ is a normal subgroup of
$A_4$. A 3-cycle in $A_4$ determines its conjugacy class by whether it maps to  $(1\,2\,3)$ or
$(1\,3\,2)$ in  $A_4/K=\bZ/3=A_3$.  

\begin{lem} \label{allow} With $G=A_3$ or $A_4$, $\ni(G,\C_{\pm 3^{s_1,s_2}})$ is nonempty if and only if  
\begin{equation} \label{dagger} s_1-s_2\equiv 0
\bmod\ 3 \ (s_1+s_2=r).\end{equation}   

Subject to \eqref{dagger}, $(s_1,s_2)$ (resp.~unordered pairs $\{s_1,s_2\}$) label nonempty inner
(resp. absolute) Nielsen classes of 3-cycle conjugacy classes in either $A_3$ or
$A_4$. \end{lem} 

It is
convenient to select  $(2\,3)=\beta$ to  conjugate a 3-cycle  in $A_4$ 
to the conjugacy class of its inverse.
For any 3-cycle $\alpha\in A_4$, denote its conjugate 
$\beta\alpha\beta^{-1}$ by ${}^\beta\alpha$. Similarly, if $\bg$ is an
$r$-tuple of elements of $A_4$, its conjugate 
by $\beta$ is ${}^\beta\bg$. 

Let  $\bg_i$ be an $r_i$-tuple of $A_4$ 3-cycles, with $\Pi(\bg_i)=1$; $r_i$ varies with 
$i=1,\dots,t$. For $\bep\in (\bZ/2)^t$, denote $({}^{\beta^{\epsilon_1}}\bg_1,
\dots, {}^{\beta^{\epsilon_t}}\bg_t)$ by ${}^{\bep}(\row {\bg} t)$. When no other notation  suggests the 
division between  $\bg_1,\dots,\bg_t$,
replace the comma separators by semi-colons:  ${}^\bep(\bg_1;\dots;\bg_t)$ to unambiguously shows the action
of $\bep\in(\bZ/2)^t$.

\subsubsection{The 3-Lemma} \label{nt.1} We need a precise result for $G=A_3$. Assume \eqref{dagger}. 

\begin{lem}[$3$-Lemma] \label{lnt.1} $B_r$ applied to 
$((1\,2\,3)^{(s_1)}\!,\!(3\,2\,1)^{(s_2)})$ is $\ni(A_3,\bfC_{\pm3^{s_1,s_2}})$. 

If $\bg^*=(\bg,\bg')\in \ni(A_n,\bfC_{3^{r'}})$, $n\ge 5$,  and $\bg\in
\ni(A_3,\bfC_{\pm3^{s_1,s_2}})$, $r=s_1+s_2$, then there is 
$Q\in B_{r'}$ with
$(\bg^*)Q=(\bar\bg,\bg')$ where $\bar\bg$ is 
$$\begin{array}{rl}&((1\,2\,3),(3\,2\,1))^{({\frac r 2})} \text{ if $r$ is even}, and \\
&((1\,2\,3)^{(3)},((1\,2\,3),(3\,2\,1))^{(\frac{r-3}{2})}) \text{ if $r$ is odd}.\end{array}$$ \end{lem}

\begin{proof} Since $A_3$ is cyclic of order 3, the first statement is obvious. 

Since $n\ge 5$,  apply Blocks Lem.~\ref{lnf.8} to conjugate $\bg$ by $\gamma=(2\,3)(k\,j)$ with $k$ and $j$
any  integers distinct from  1, 2 and 3. So, with no loss assume $s_1\ge s_2$.  Braid $\bg$ to 
$((1\,2\,3)^{(s_1-s_2)}, ((1\,2\,3),(3\,2\,1))^{(s_2)})$. With no loss, take $s_2=0$.
Thus,  \eqref{dagger} implies 3 divides $s_1$. We take $s_1$
even; the other case is similar. So, 
$\bg$ is $(\row {\bg} {s_1/3})$ with each $\bg_i$ equal
$(1\,2\,3)^{(3)}$. By assumption $\lrang{\bg_i,\bg'}=A_n$, $n\ge 5$. 
Several applications of the the Blocks Lemma, using $\gamma$ above, produces $Q'\in B_r$ with 
$$(\bg,\bg')Q'=(\bg_1,\gamma\bg_2\gamma^{-
1},\bg_3,\gamma\bg_4\gamma^{-1},\dots,\gamma\bg_{s_1/3}\gamma^{-
1},\bg').$$ This the desired target with $r$ even. 
\end{proof}

\subsection{The case $n=4$} \label{nt}
Most difficulties 
are in this induction on $r$ for $n=4$.

\subsubsection{$A_4$ targets} \label{nt.2} Conjugating 
by $\beta=(2\,3)$ switches $s_1$ and $s_2$ in list \eqref{ent.1}. 
\begin{edesc} \label{ent.1} \item \label{ent.1a} $r=3, s_1=3, s_2=0:
\qquad \qquad \bg_{3,-}=((1\,2\,3),(1\,3\,4),(1\,4\,2))$.
\item\label{ent.1b} $r=4, s_1=s_2=2: 
\begin{array}{rl}\qquad \ \ \ \ \bg_{4,+}=&((1\,3\,4),(1\,4\,3),(1\,2\,3),(1\,3\,2)),\cr
\qquad \ \ \ \ \bg_{4,-}=&((1\,2\,3), (1\,3\,4), (1\,2\,4),
(1\,2\,4)).\end{array}$ 
\item \label{ent.1c} $r\ge 5, s_1=3+s_1',s_2=s_2':
\begin{array}{rl} \bg_{r,+}=& ((1\,3\,4)^{(3)},
(1\,2\,3)^{(s_1')},(3\,2\,1)^{(s_2')}),\\ 
\bg_{r,-}=&(\bg_{3,-},
(1\,2\,3)^{(s_1')},(3\,2\,1)^{(s_2')}).\end{array}$\end{edesc}
 
\begin{lem}[4-Lemma] \label{lnt.2}  Assume \eqref{dagger} holds for $(s_1,s_2)$. Then,     
any \\ $\bg\in \ni(A_4,\bfC_{\pm \C_{{3^{s_1,s_2}}}})$ braids either to an element in
\eql{ent.1}{ent.1a},
\eql{ent.1}{ent.1b} or \eql{ent.1}{ent.1c}, or its conjugate by $\beta$. If $s_1=s_2$, some braid achieves 
conjugation by $\beta$.
\end{lem}

We divide the proof into four subsections.  The first on $r=3$ and $r=4$ showing how to use the
$\sh$-incidence matrix. The next two treat separately when 
\eql{A3list1}{(i)} and 
\eql{A3list1}{(ii)} hold, inducting on
$r$ using Coalescing Lemma \ref{lnf.2}. The last considers the case $s_1=s_2$ to show conjugation by $\beta$
is braidable.

\subsubsection{$r=3$ and 4, and the $\sh$-incidence matrix} \label{shinc} Modulo conjugation by $S_4$ and
action of
$B_3$ here are the strings of two  or three 3-cycles with product 1. 
\begin{equation} \label{ent.2} ((1\,2\,3),(3\,2\,1)),\  (1\,2\,3)^{(3)}, \ 
((1\,2\,3),(1\,3\,4),(1\,4\,2)).\end{equation}
Only the 3rd is transitive on $\{1,2,3,4\}$. 
This finishes the case $r=3$.  

\cite[\S2.10]{BaFr} has the rubric for the $\sh$-incidence matrix. It works for
all values of $r$, though  for $r=4$ it is usually possible to do it by hand. 
The result is information on natural $j$-line (curve) covers (reduced Hurwitz spaces as in \S\ref{rnz.4}) one
can  compare with modular curves (which are a special case). We do 
$\ni(A_4,\bfC_{\pm 3^2})$ here. \cite[\S9]{BaFr} and \cite[\S6, \S 7.2]{thomp-gen0} interpret these 
computations. 

The computation works by using these three important groups:
\begin{edesc} \label{cuspgp} \item \label{cuspgpa} $\sQ''=\lrang{\sh^2, q_1q_3^{-1}}$; \item
\label{cuspgpb} the {\sl cusp group\/}
$\text{Cu}_4=\lrang{q_2,\sQ''}/\sQ''$; and
\item \label{cuspgpc} the {\sl mapping class group\/} $\bar M_4=\lrang{\gamma_0,\gamma_1}=H_4/\sQ''$ generated
freely by
$\gamma_0=q_1q_2$,  
$\gamma_1=q_1q_2q_3=\sh$ (\S\ref{BHMon}) of respective orders 3 and 2.  
\end{edesc} 

This induces an action of $\bar M_4$ on $\ni(G,\bfC)^*/\sQ''=\ni(G,\bfC)^{*,\rd}$ ($*=\inn$ or $\abs$), {\sl
reduced\/} Nielsen classes. That is,   
$\gamma_0,\gamma_1,\gamma_\infty$ in \eqref{cuspgp} are names for $H_4$ elements on reduced
classes.   The orders of $\gamma_0$ and $\gamma_1$ in \eql{cuspgp}{cuspgpc} come easily from the Hurwitz
relation \eqref{Hurwitz} mod $\sQ''$. So, too, does the relation $\gamma_0\gamma_1\gamma_\infty=1$, with 
$\gamma_\infty=q_2$. 

The cover $\bar \sH(G,\bfC)^{*,\rd}
\to \prP^1_j$ (as in Prop.~\S\ref{pnz.5}) has as branch cycles (Def.~\ref{branchcycles})
$(\gamma_0,\gamma_1,\gamma_\infty)$ on  
$\ni(G,\bfC)^{*,\rd}$. This gives us the genus of its components.

A pairing on  $\gamma_\infty$ orbits
$(\bg)\text{Cu}_4=O=O_\bg$ gives $\sh$-incidence matrix  entries:  $(O,O')\mapsto $ the
cardinality of $O$ intersected with the shift on $O'$:   
$$\smatrix {\vdots} {|O\cap (O')\gamma_1|} {\cdot} {\cdots} = \smatrix {\vdots}  {|O\cap (O')\gamma_0|}
{\cdot} {\cdots}.$$ 
\begin{itemize}\item  Blocks $\Leftrightarrow$ components of $\sH(G,\bfC)^{*,\rd}$, or of $\sH(G,\bfC)^*$.
\item Fixed points of $\gamma_0$ or $\gamma_1$ appear on the diagonal. 
\end{itemize} 

Consider $g_{1,4}=\bg_{4,-}\in \ni(A_4,\bfC_{\pm 3^2}) =((1\,2\,3), (1\,3\,4), (1\,2\,4), (1\,2\,4))$. from
Cor.~\ref{nf.2}, $s(g_{1,4})=s((1\,2\,3), (1\,3\,4), (1\,4\,2))=-1$. 

Subdivide $\mapsto \ni(A_3,\bfC_{\pm 3^2})^{\inn,\rd}$ according to the 
sequences of conjugacy classes $\C_{\pm 3}$;   $q_1q_3^{-1}$ and $\sh$ switch these rows: 
$$\begin{tabular}{cccc} &[1] +\,-\,+\,- &[2] +\,+\,-\,- &[3] +\,-\,-\,+ \\
&[4] -\,+\,-\,+ &[5] -\,-\,+\,+ &[6] -\,+\,+\,-
\end{tabular} $$
 
Here is the notation in the charts below where $O_{i,j}^k$ appears: $k$ is the cusp width, and $i,j$
corresponds to a labeling of orbit representatives. The diagonal entries for 
$O_{1,1}^4$ and $O_{1,4}^4$ are nonzero. In detail, however, $\gamma_1$
(resp.~$\gamma_0$) fixes 1 (resp.~no) element of $O_{1,1}$, and neither of $\gamma_i$, $i=0,1$, fix any
element of $O_{1,4}$. 
 
 $$\begin{array}{rl} \text{H-M rep.} \mapsto \bg_{1,1}=& ((1\,2\,3),
(1\,3\,2), (1\,3\,4), (1\,4\,3))\\
\bg_{1,3}=& ((1\,2\,3), (1\,2\,4), (1\,4\,2), (1\,3\,2))\\
\text{H-M rep.} \mapsto  \bg_{3,1}=& ((1\,2\,3), (1\,3\,2), (1\,4\,3), (1\,3\,4))
\end{array} $$ 

\begin{table}[h] \label{sh-incA4}
\begin{tabular}{|c|ccc|}  \hline $\ni_0^+$ Orbit & $O_{1,1}^4$\ \vrule  &
$O_{1,3}^2$\ \vrule & 
$O_{3,1}^3$ \\ \hline $O_{1,1}^4$ 
&1&1&2\\ 
$O_{1,3}^2$ &1 &0&1  \\ $O_{3,1}^3$ &2&1&0 \\  \hline \end{tabular} 
\ 
\begin{tabular}{|c|ccc|} \hline $\ni_0^-$ Orbit & $O_{1,4}^4$\ \vrule  &
$O_{3,4}^1$\ \vrule & 
$O_{3,5}^1$ \\ \hline $O_{1,4}^4$ 
&2&1&1\\ 
$O_{3,4}^1$ &1 &0&0  \\ $O_{3,5}^1$ &1&0&0 \\  \hline \end{tabular} \end{table}

\begin{prop} \label{A43-2} On $\ni(\Spin_4,\bfC_{\pm3^2})^{\inn,\rd}$
(resp.~$\ni(A_4,\bfC_{\pm3^2})^{\inn,\rd}$)
 $\bar M_4$ has one (resp.~two) orbit(s). So,  $\sH(\Spin_4,\bfC_{\pm3^2})^{\inn,\rd}$ 
 (resp.~\!$\sH(A_4,\bfC_{\pm3^2})^{\inn,\rd}$) has one (resp.~two)  component(s),  
 $\sH_{0,+}$ (resp.~$\sH_{0,+}$ and $\sH_{0,-}$). 

Then, $\sH(\Spin_4,\bfC_{\pm3^2})^{\inn,\rd}$ maps one-one to 
$\sH_{0,+}$ (though changing $A_4$ to $\Spin_4$ give different moduli).
The compactifications of $\sH_{0,\pm}$ both have genus 0 from Riemann-Hurwitz applied to 
$(\gamma_0,\gamma_1,\gamma_\infty)$ on reduced Nielsen classes (Ex.~\ref{exA43-2}).\end{prop} 

\begin{rem}[Complication in  \eqref{ens.2a}]
\label{comp4Min} Prop.~\ref{A43-2} says the lifting invariant separates the two $B_4$ orbits on
$\ni(A_4,\bfC_{\pm 3^2})$. We use this for convenient substitution, when we have all but the
first 4 entries matching an item in list 
\eqref{ens.1a}. For example: We can substitute  
$\bg'=((1\,n\nm1\,n\nm2)^{(2)},(1\,n\nm1\,n),(1\,n\,n\nm 2))$ for any 3-cycle  4-tuple $\bg$, with
product-one and $s(\bg)=-1$,  on
$\{1,n\nm 2,n\nm1,n\}$. 
\end{rem}

\subsubsection{Case \eql{A3list1}{(i)}} \label{liftone} Now assume $r\ge 5$. 
With no loss, $\bg$ braids to $(g,g^{-1},\bg')$ 
 $\Pi(\bg')=1$. 
Suppose $\bg'$ is intransitive on $A_4$.  Then,  $\lrang{\bg'}=A_3$. Apply \ 
Product-one Lemma \ref{lnf.5} to find $Q\in B_r$, $\gamma\in A_4$ and $j\in\{2,3\}$ 
with $$\bg''=(\bg)Q=\gamma\bg\gamma^{-1}=((1\,j\,4),(4\,j\,1);(1\,2\,3)^{(s_1')},(3\,2\,1)^{(s_2')}),$$  
with
$s_1'\equiv s_2'\mod 3$. Since $r\ge 5$, some braid of $\bg''$ puts one of these at its 
head: $$((1\,j\,4),(4\,j\,1),\!(1\,2\,3)^{(3)}), \!((1\,j\,4),(4\,j\,1),\!(3\,2\,1)^{(3)}),
\!((1\,j\,4),(4\,j\,1),\!(1\,2\,3)^{(2)},\!(3\,2\,1)^{(2)}).$$ It is easy to braid the first two to
\eql{ent.1}{ent.1c}. Example: For the 1st, if $j=2$,  conjugate by
$(2\,3\,4)$ (Prod.~One Lem.) and then  slide $(1\,3\,4)^{(3)}$ to the
front (Blocks Lem.~\ref{lnf.8}).  

Now suppose $\bg'$ is transitive on $A_4$. Blocks Lemma \ref{lnf.8} gives $Q\in B_r$ 
with $(\bg)Q=((1\,2\,3),(3\,2\,1),\bg')$. The induction assumption gives   a 
$Q'$ putting $\bg'$ in a preferred form depending on $s_1'$ and 
$s_2'$. If $r-2\ne 4$, or if $r-2=4$ and $s(\bg)=+1$, the Blocks Lemma allows combining the head
$((1\,2\,3),(3\,2\,1))$ with the tail   
$((1\,2\,3)^{(s_1'')},(3\,2\,1)^{(s_2'')})$ for a 
preferred form  in \eqref{ent.1}.  

That leaves only deciding how to braid $((1\,2\,3),(3\,2\,1), \bg_{4,-})$ to a normal form.  The Blocks Lemma
braids this to $(\bg_{4,-},(1\,2\,4),(4\,2\,1))$, with $(1\,2\,4)^{(3)}$ in the 3rd through 5th entries.
Apply it again to get $\bg_{6,-}$ (in \eql{ent.1}{ent.1c}). 

\subsubsection{Case \eql{A3list1}{(ii)}}   Without loss, 
assume we can't braid to case \eql{A3list1}{(i)} and   $\bg$ has the form  
$((1\,3\,4),(1\,4\,2),\bg')$. Apply the induction assumption to 
$((1\,3\,2),\bg')=\bg^*$. 
If $\bg^*$ is intransitive on $A_4$, then 
$\bg^*=((1\,2\,3)^{(s_1')},(3\,2\,1)^{(s_2')})$. We may assume 
$$\bg=((1\,3\,4),(1\,4\,2),(1\,2\,3)^{(s_1')}, (3\,2\,1)^{(s_2'\nm1)}).$$ Braid 
this 
to $((1\,2\,3),(1\,3\,4),(1\,4\,2), (1\,2\,3)^{(s_1'\nm1)}, 
(3\,2\,1)^{(s_2'\nm1)})$ ($\bg_{r,-}$) by 
braiding \\ $(1\,2\,3)$ past $((1\,3\,4),(1\,4\,2))$. 

Now assume  
$\bg^*$ is 
transitive on $A_4$. Recall  $\alpha: B_r\to S_r$ (\S\ref{BHMon}).  Apply the
induction  with $r{-}1$ 
replacing $r$. Denote the $i$th  entry of $(\bg)Q'$ by 
$(\bg)Q'[i]$. Let $i=(1)\alpha(Q)$. Thus, there is $Q\in B_{r-1}$ and $Q'\in B_r$ with the 
following properties. 
\begin{edesc} \label{ent.4} 
\item \label{ent.4a} $(\bg^*)Q$ is in the list \eql{ent.1}{ent.1a}--\eql{ent.1}{ent.1c} with $r\nm1$ 
replacing $r$.
\item \label{ent.4b} $(\bg)Q'[j]=(\bg^*)Q[j]$, $j=1,\dots,i-1$, 
$(\bg)Q'[i](\bg)Q'[i\np 1]=(\bg^*)Q[i]$ \hfil \break and
$(\bg)Q'[j\np 1]=(\bg^*)Q[j]$, 
$j=i+1,\dots,r-1$. \end{edesc}

Consider possibilities for  reexpanding $(\bg^*)Q$ to give $(\bg)Q'$ by putting two
3-cycles (using  
\eql{ent.1}{ent.1a}) with product $(\bg^*)Q[i]$ in its place.  We can dismiss all but
one by showing, contrary to assumption,  case \eql{A3list1}{(ii)} holds. The case $r{-}1=3$ 
illustrates with  $(\bg^*)Q=\bg_{3,-}$. With $((1\,2\,4),(1\,4\,3))$ 
in place of 
$(1\,2\,3)$, the result is $((1\,2\,4),(1\,4\,3),(1\,3\,4),(1\,4\,2))$. Thus, $(\bg)Q'$ 
is in case \eql{A3list1}{(i)}. By braiding it is clear this one substitution 
suffices for 
$r-1=3$. 

Now consider  $r-1\ge 4$. For $\bg_{4,+}$, no matter the substitution, you are 
in case \eql{A3list1}{(i)}. By braiding, assume substitutions in $\bg_{4,-}$ are for 
$(1\,3\,2)$ or for $(1\,4\,2)$. For the first, this produces 
$((1\,3\,4),(1\,4\,2),(1\,3\,2),(1\,3\,4),(1\,4\,2))$. Apply $Q_2^{-1}Q_3^{-1}$ 
to 
get $((1\,3\,4),(1\,3\,4),(1\,2\,3),(2\,3\,4),(1\,4\,2))$. Now apply $Q_3$ to 
get
$$((1\,3\,4),(1\,3\,4),(1\,2\,4),(1\,2\,3),(1\,4\,2)).$$ With $(1\,2\,4)$ (resp.~$(1\,4\,2)$) in the 3rd
(resp.~5th) position, this braids to 
\eql{A3list1}{(i)}.  
Similarly, substitution for $(1\,4\,2)$ gives  
$((1\,3\,2),(1\,3\,2),(1\,3\,4),(1\,4\,3),(1\,3\,2))$ in \eql{A3list1}{(i)}. 

Assume $r\ge 5$ is odd. Substitution for $(1\,3\,4)$ in $\bg_{r,+}$ gives 
$$((1\,3\,4)^{(2)},(1\,3\,2),(1\,2\,4),(1\,2\,3)^{(s_1')},(3\,2\,1)^{(s_2')}).$$ 
Since 
$s_1'>0$, this is 
\eql{A3list1}{(i)}. Substitution in $\bg_{r,+}$ for $(1\,2\,3)$ gives 
\!\begin{equation} \label{ent.5a} ((1\,3\,4)^{(3)},
(1\,2\,3)^{(s_1'\nm1)},(1\,2\,4),(1\,4\,3),(3\,2\,1)^{(s_2')}).\end{equation} As  
$(1\,3\,4)$ and $(1\,4\,3)$ appear, this 
braids to  
\eql{A3list1}{(i)}. Finally,  
substitute for $(3\,2\,1)$: 
\!\begin{equation} \label{ent.5b} ((1\,3\,4)^{(3)},
(1\,2\,3)^{(s_1')},(1\,3\,4),(1\,4\,2 ),(3\,2\,1)^{(s_2'\nm 1)}).\end{equation} 
This is \eql{A3list1}{(i)} if $s_2'>1$. So, let $s_2'=1$ and $s_1'\equiv 1 \bmod 3$. Rewrite 
\eqref{ent.5b} as 
\begin{equation} \label{ent.5c} ((1\,3\,4)^{(3)},
(1\,2\,3)^{(s_1'\nm 1)},(1\,2\,3),(1\,3\,4),(1\,4\,2 )).\end{equation} 
Since $(2\,3\,4)\in \lrang{(1\,2\,3),(1\,3\,4),(1\,4\,2 )}$, Blocks Lem.~\ref{lnf.8} gives a braid that
conjugates 
$(1\,3\,4)^{(3)}$ to
$(1\,2\,3)^{(3)}$. Apply it again to braid    
\eqref{ent.5c} to  $$((1\,2\,3),(1\,3\,4),(1\,4\,2 ), (1\,2\,3)^{(s_1'\np
2)}) \text{: $s_2'=0$ in  $\bg_{r,-}$,  
\eql{ent.1}{ent.1c}}.$$ The 
$\bg_{r,-}$ substitutions are easier for they imitate previous substitutions. 
Example: substituting in one of the beginning three entries  
duplicates the case $r-1=3$. 

\subsubsection{Braiding $\beta$ when $s_1=s_2$} 
It suffices to  braid $\beta$ for any
braid orbit rep. If $\bg$ has lifting invariant +1, and $s_1=s_2$, then it braids to $\bg''$ in
\S\ref{liftone} (so
$s_1'=s_2'$). Example: If
$j=2$, then
$\beta\bg''\beta^{-1}=((1\,3\,4),(1\,4\,3);(3\,2\,1)^{(s_1')},(1\,2\,3)^{(s_1')})$. Apply the Prod.~One
Lem.~to conjugate by $(1\,2\,3)$.  Then braid the commuting pieces $(3\,2\,1)^{(s_1')}$ and
$(1\,2\,3)^{(s_1')}$ past each other to return to $\bg''$. 

Prop.~\ref{A43-2} showed, for $r=4$, any $\bg$ with lifting invariant -1 braids to $\bg_{4,-}$.
As above, it now suffices to show we can braid $\beta$ when the lifting invariant is -1 to  
$\bg_{r,-}$,
$r\ge 5$. Similar to braids above, braid $\bg_{r,-}$ to 
$$\bg^\dagger=((1\,3\,2)^{(2)},(1\,3\,4),(1\,4\,2); (1\,2\,3)^{(s_1'+1)},(3\,2\,1)^{(s_2'-2)}), \text{
with $s_2'=s_1'\np3$}.$$ Again, ~Prop.~\ref{A43-2} braids the conjugate by $\beta$ on the
4-entry head of $\bg^{\dagger}$, and the paragraph above just braided it on its tail.
That completes Lem.~\ref{lnt.2}.  

\section{Improving the Coalescing Lemma and full induction} \label{nf.4} 
We improve Coalescing Lemma \ref{lnf.2} to show for $n\ge 5$ we can braid to    
\eql{A3list1}{(i)}. 

\subsection{Set up for Strong Coalescing} Consider how a  
{\sl disappearing sequence\/} in Lem.~\ref{lnf.2}  produces   \eql{A3list1}{(i)}, \eql{A3list1}{(ii)} or
\eql{A3list1}{(iii)}.  Use $l$ for the length
of the sequence. 

\subsubsection{Some tough braidings} 
Coalescing types \eql{A3list1}{(i)} or \eql{A3list1}{(ii)} have $l=2$. 
Condition \eqref{enf.1c} on the shortest disappearing sequence is an induction assumption: no element
of $B_r^{(l)}$ takes our choice of disappearing sequence to a smaller disappearing subsequence. From this,
Lem.~\ref{lnf.2}  concludes $l=3$ and coalescing type \eql{A3list1}{(iii)} occurs. Though $n=3$ needs type
\eql{A3list1}{(iii)}, we will see  that case is special. 

Here are hard cases with $n=4$ and 5 ($r=6$) for  
finding coalescing type \eql{A3list1}{(i)}:  
\begin{equation} \label{enf.9a}  \bg_0=((1\,2\,3)^{(3)},(1\,4\,5)^{(3)}); \end{equation}
\begin{equation} \label{enf.9b} \bg_1=((1\,2\,3)^{(3)},(2\,1\,4)^{(3)})\text{; or }
\bg_2=((1\,2\,3)^{(3)},(1\,2\,4)^{(3)}).\end{equation} We introduce an extra notation now for later use. If
in addition to dropping the tail of the \eqref{ens.1a}${}_+$ term (as in \S\ref{nfComm}), we drop the head
$\sbgp {\frac{n-3} 2}$, the nub that remains is conjugate to $(1\,2\,3)^{(3)}$. All
6-tuples in \eqref{enf.9a} and \eqref{enf.9b}  are juxtapositions of two of these. We refer to the
resulting 6-tuples as (having type) 
\eqref{ens.1a}${}_+$+\eqref{ens.1a}${}_+$ 

In \eqref{enf.9a}, each integer, $i=1,\dots,5$, 
has a disappearing sequence of length three. Also, there is no 
$1\le i< j\le 6$ with $(g_i,g_j)$ a disappearing sequence for any integer. 
Still, Lem.~\ref{3-3braids} shows we can braid \eqref{enf.9a} to type \eql{A3list1}{(i)}. 

\subsubsection{Coalescing tricks} 
Lem.~\ref{3-3braids} conceptually finds explicit braids forced on
us. Later we will be less complete; all are similar. Just as 
$A_4$ has  $ \C_{\pm3}$, $A_5$ has the two 5-cycle classes $\C_{\pm 5}$.  
Similarly, as with the  two $B_4$ orbits on $\ni(A_4,\C_{\pm3^2})^\inn$ (\S\ref{2-3-cycles}) separated
by  lift invariants,  there is a similar result for $(A_5,\bfC_{\pm 53^2})$. 
\begin{lem} \label{3-3braids} There are two $B_4$ orbits on each of 
\begin{equation} \label{5-cycleBraids} \ni(A_5,\bfC_{\pm 53^2})^\inn, \ni(A_5,\bfC_{+ 5^23^2})^\inn\text{ and
}
\ni(A_5,\bfC_{\pm 5^2})^\inn,\end{equation} separated by lift invariants. There is one $B_4$ orbit on
$\ni(A_5,\bfC_{3^4})^\inn$. 

An H-M rep.~braids to 
$\bg_0$ in
\eqref{enf.9a} (so, to a 6-tuple of type
\eql{A3list1}{(i)}). Some braid takes  $\bg_1$ and $\bg_2$ in \eqref{enf.9b} to  an H-M rep. or to a
juxtaposition of two conjugates to $\bg_{3,-}$ (in \eql{ent.1}{ent.1a}; see Rem.~\ref{twoTypes}).  
 \end{lem} 

\begin{proof} 4-Lemma \ref{lnt.2}  shows how the $\sh$-incidence matrix gives  
two $B_4$ orbits on 
$\ni(A_4,\C_{\pm3^2})^\inn$ (with reps.~$\bg_{4,\pm}$; \cite[\S2.10.1]{BaFr} uses
$\ni(A_5,\bfC_{3^4})^\inn$ for similar purposes). \cite[Prop.~5.11]{BaFr} showed the orbit results for
\eqref{5-cycleBraids}. 

Now consider $\bg_1$. We easily braid it to
$$\bg_1'=((2\,1\,4),(1\,2\,3),(1\,2\,3),(1\,2\,3),(2\,1\,4),(2\,1\,4)).$$ Note: The first two entries, and the
4th and 5th entries have coalescing type \eql{A3list1}{(ii)}. Invariance Cor.~\ref{cnf.3} says  $s(\bg_1)$
is
$(-1)^2s(\bh_1)$ with
$$\bh_1=((1\,4\,3)^*,(1\,2\,3),(2\,3\,4)^*,(2\,1\,4))\in \ni(A_4,\C_{\pm3^2})^\inn.$$ 
The * superscripts on $\bh_1$ entries remind those positions are coalescings. 
Braids we now apply to $\bh_1$ track the *$\,s$. Also,
$(\bh_1)Q_2^{-1}=((1\,4\,3)^*,(2\,3\,4)^*,(3\,4\,1),(2\,1\,4))$. 

Prop.~\ref{A43-2} implies $\bh_1$ braids to $\bg_{4,+}=((1\,3\,4),(1\,4\,3),(1\,2\,3),(1\,3\,2))$. Now we ask,
where did the *$\,$s end up? If in the 1st and 2nd (resp.~3rd and 4th) positions, reexpand so 
$((1\,2\,3),(1\,3\,2))$ are in the 5th and 6th position; we braided
$\bg_1$ to have type \eql{A3list1}{(i)}. If, rather, the *$\,$s fall one in the 1st
or 2nd, the other in the 3rd or 4th, then reexpanding gives two juxtaposed 
\eql{A3list3}{(ii')} types.  

For
$\bg_2$, follow part of the plan for $\bg_1$, braiding to
$$((1\,2\,4),(1\,2\,3),(1\,2\,3),(1\,2\,3),(1\,2\,4),(1\,2\,4)).$$ Apply
$Q_1^{-1}Q_4$ to get $((1\,4\,3),(1\,2\,3))$ in the 4th and 5th entries, and $((1\,2\,3),(2\,3\,4))$ in the
1st and 2nd positions. Now use the shift to braid to  where the first four positions are
$((1\,2\,4),(2\,3\,4),(1\,4\,3))$.  This is coalescing type \eql{A3list3}{(ii')}. 

Now consider $\bg$. Its lifting invariant is clearly one. As with $\bg_1$, braid to 
$$\bg'=((1\,4\,5),(1\,2\,3),(1\,2\,3),(1\,2\,3),(1\,4\,5),(1\,4\,5)).$$ The pair $((1\,4\,5),(1\,2\,3))$ has
product $(1\,4\,5\,2\,3)$ and the pair $((1\,2\,3),(1\,4\,5))$ has product $(1\,2\,3\,4\,5)$. Apply the
2-orbit outcome for $\ni(A_5,\bfC_{+ 5^23^2})$ to see that
$$\bh=((1\,4\,5\,2\,3)^*,(1\,2\,3),(1\,2\,3\,4\,5)^*,(1\,4\,5))$$ (with the *$\,$s as above) braids
to $((1\,2\,3\,4\,5)^*,(5\,4\,3\,2\,1)^*,(1\,2\,3),(3\,2\,1))$, and we know where the *$\,$s  
must end up. So, upon their expansion we have a 6-tuple \\ 
$\bg''=(\bg^\dag,(1\,2\,3),(3\,2\,1))$ of type
\eql{A3list1}{(i)}. Go further, still: $\bg^\dag\in \ni(A_5,\bfC_{3^4})^\inn$ and it has lifting invariant 1.
By the first statement of the lemma, it braids to an H-M. 
This concludes showing 
\eqref{ens.1a}${}_+$+\eqref{ens.1a}${}_+$ for $n=5$ braids to an H-M rep.
 \end{proof} 

\begin{rem}[Which coalescing type] \label{twoTypes} \S\ref{part2}  shows  
$\bg_1$ and
$\bg_2$ braid to both coalescing types \eql{A3list1}{(i)} and
\eql{A3list3}{(ii')}, though Lem.~\ref{3-3braids} left this ambiguous. \end{rem}

\begin{rem}[Using Braid packages] \label{braidPack} \cite{MSV} describes recent Braid package 
applicable to our by-hand calculation in Lem.~\ref{3-3braids}. 
Our approach conceptually shows why the braidings give
components separated by lifting invariants. Still,  \GAP\ is additional corroboration. The package documented
by
\cite{MSV} does not have the
$\sh$-incidence matrix of \cite{BaFr} used in \S\ref{shinc}. 
\end{rem} 
 
\subsection{Strong Coalescing} Suppose  
$\bg\in \ni(A_4,\Ct r)$ with $r\ge 3$. Excluding the case where $r=4$ and $s(\bg)=-1$, 4-Lemma
\ref{lnt.2} gives a braid of $\bg$    
to $\bg'$ with either $(g'_1,g'_2)$
of type \eql{A3list1}{(i)} (conjugate to $((1\,2\,3),(3\,2\,1))$) or
\eql{A3list3}{(ii')} ($(g_1',g_2',g_3')$ conjugate in $S_4$ to $((1\,2\,3),(1\,3\,4),(1\,4\,2))$. 

Lem.~\ref{lnf.2'} reduces proving Thm.~\ref{thmA} and Thm.~\ref{thmB} to cases 
like those of Lem.~\ref{3-3braids}. 

\begin{lem}[Strong Coalescing] \label{lnf.2'}
If $n\ge 5$, $r\ge 4$, then  $\bg\in \ni(A_n,\Ct r)$ braids to $\bg'$ with 
$(g'_1,g'_2)$ of coalescing type \eql{A3list1}{(i)}. \end{lem} 

\S\ref{part1} proves, for $n\ge 5$, we get either coalescing type \eql{A3list1}{(i)}, or a sum of types, $T_1+T_2$
where the $T_i\,$s are either  \eql{A3list3}{(ii')} or  \eql{A3list3}{(iii')}  (see \eqref{enf.10}). Then, \S\ref{part2} produces from this type \eql{A3list1}{(i)} and the proof
of Lem.~\ref{lnf.2'}. 

\subsubsection{Proof of type \eql{A3list1}{(i)},  or type $T_1+T_2$ as above}  \label{part1} 
Suppose no braid of
$\bg$ has coalescing type either
\eql{A3list1}{(i)} or
\eql{A3list1}{(ii)}. Then Coalescing Lemma \ref{lnf.2} shows each $i\in \{1,\dots,n\}$ occurs
in a disappearing sequence of type $(i\,i_1\,i_2)^{(3)}$ in $(\bg)Q$
for every $Q\in B_r$. We show this is impossible. To simplify,
take $(i\,i_1\,i_2)=(1\,2\,3)$. Then, there is $Q\in B_r$ with $(\bg)Q=((1\,2\,3)^{(3)}, \bg^*)$. Now
apply the same argument to $\bg^*$, and reduce the case to one of the 6-typles $\bg_0$, $\bg_1$ or 
$\bg_2$ in either \eqref{enf.9a} or \eqref{enf.9b}. 

Now we may assume $\bg$ braids to $\bg'$ with $(g_1',g_2')$ of coalescing type \eql{A3list1}{(ii)}.  Coalesce
the first two positions of $\bg'$ to $\bh$ with $g_1'g_2'=h_1^*$, its remaining entries are in order
from $\bg'$. As in the proof of Lem.~\ref{3-3braids}, *  tracks where the coalesced  entry ends up in the
braiding. Suppose
$\lrang{\bh}$ is transitive on a set containing at least five integers. Then, our induction assumption
shows the only nontrivial case to be when the transitive set of 3-cycles includes $h_1^*$. With no loss,
either:  
\begin{edesc} \label{h*exp} \item \label{h*expa} the first two entries of $\bh$ are $(h_1^*,(4\,2\,1))$ of
type \eql{A3list1}{(i)}; or   
\item \label{h*expb} the first three entries are $(h_1^*,h_2,h_3)$ of type \eql{A3list3}{(ii')}.
\end{edesc}  
Expand $h_1^*$. For \eql{h*exp}{h*expa} we can braid from $\bg$ to one of type
\eql{A3list3}{(ii')}, our desired conclusion. For \eql{h*exp}{h*expb}, Lem.~\ref{3-3braids} braids to where
the 1st and 2nd (resp.~3rd and 4th) entries have type \eql{A3list1}{(i)}. 
Now assume the orbits of $\bh$ have four or three integers. 

If all orbits of $\lrang{\bh}$ have cardinality 3, then restrict to one in  
the support of $g_1'g_2'$ (take this to be $(1\,2\,4)$ for simplicity). Exclude the previous
case. Then, 3-Lemma~\ref{lnt.1} reduces us to  $\bh=((1\,2\,4)^{(3)}, g^{(3)})$  with
$g=(3\,5\,6)$.  Expand $\bh$ to show  $\bg_{4,-}$  (in
\eql{ent.1}{ent.1b}) and $g^{(3)}$ juxtaposing. 
Braid \!$\bg_{4,-}$ to get \!$((1\,2\,3),(1\,2\,4)^{(2)},(3\,2\,4),g^{(3)})$. 
Then, coalesce $(g,(1\,2\,3))$ from the 1st and last 3-cycles. Braid the result, $(3\,5\,6\,1\,2)$, past $(1\,2\,4)^{(2)}$ to replace that by 
$(2\,3\,4)^{2}$.  A braid now taking $((2\,3\,4),(3\,2\,4))$ to the 1st and 2nd position shows we can braid $\bg$ to type \eql{A3list1}{(i)}: our ultimate goal. 

The final case reduces to where $\lrang{\bh}$ is transitive on 4 integers, and
the support of 
$(g_1',g_2')$ includes a 5th integer. Now apply 4-Lemma \ref{lnt.2}.  This is similar to the case we just
finished, except for one possibility: $\bh=\bg_{4,-}$. Cor.~\ref{cnf.3} shows $s(\bg)=-s(\bh)=+1$. In the
expansion of $\bh$, coalesce the last two terms from $\bg_{4,-}$ to get $\bh^\dag\in \ni(A_5,\bfC_{3^4})$. 
From Lem.~\ref{3-3braids} this braids to an H-M rep., so the expanded $\bh^\dag$ braids to type
\eql{A3list1}{(i)} for the conclusion.

\subsubsection{Finish of producing type \eql{A3list1}{(i)}} \label{part2} Assume, contrary to our
goal,  no braid has type \eql{A3list1}{(i)}. Apply \S\ref{part1} to braid $\bg$ to have the shape   
$(\bg_1',\bg_2',\bg_3')$ as follows. 
\begin{edesc} \label{enf.10}
\item \label{enf.10a} $\bg_i'$, $i=1,2$, is conjugate to  $(1\,2\,3)^{(3)}$,   
$((1\,2\,3),(1\,3\,4),(1\,4\,2))$.
\item \label{enf.10b} As $n\ge 5$, the Blocks Lemma allows conjugating $\bg_1'$ and $\bg_2'$ so
each has some orbit support outside the other, though some common support. \end{edesc} 

Since $\lrang{(1\,2\,3),(1\,3\,4),(1\,4\,2)}=A_4$,  
the Blocks Lemma allows anything from $\{1,2,3,4\}$ to 
be common to all  3-cycles in $(\bg_1',\bg_2')$. So assume these  have 1 in their support. 
With $\bg_1'$ of type
\eql{A3list1}{(iii)}  and $\bg_2'$ of type \eql{A3list3}{(ii')}, with no loss, assume one of these  for $(\bg_1,\bg_2)$:
\begin{edesc} \label{enf.11} \item \label{enf.11a}
$\bh_1=((1\,2\,3)^{(3)},(1\,4\,5),(1\,5\,6),(1\,6\,4))$ (common support 1), or
 \item \label{enf.11b} $\bh_2=((1\,2\,3)^{(3)},(1\,3\,4),(1\,4\,5),(1\,5\,3))$ (common support 1 and 3).
\end{edesc}
If both $\bg_1'$ and $\bg_2'$ have type \eql{A3list3}{(ii')} we may assume one of the following:
\begin{edesc} \label{enf.12} \item \label{enf.12a} $\bh_3=((1\,2\,3),(1\,3\,4),(1\,4\,2), 
(1\,5\,6),(1\,6\,7),(1\,7\,5))$; 
or
\item \label{enf.12b} $\bh_4=((1\,2\,3),(1\,3\,4),(1\,4\,2), 
(1\,4\,5),(1\,5\,6),(1\,6\,4))$; or
\item \label{enf.12c} $\bh_5=((1\,2\,3),(1\,3\,4),(1\,4\,2), 
(1\,3\,4),(1\,4\,5),(1\,5\,3))$.\end{edesc} 

Both $\bh_2$ and $\bh_5$  have the  3-cycle 
subsequence  
$((1\,2\,3),(1\,3\,4),(1\,4\,5))$ with product  $(1\,2\,5)$.  Also, 
$\bh_4$ has the subsequence 
$((1\,3\,4),(1\,4\,5),(1\,5\,6))$ with product  $(1\,3\,6)$.  Lem.~\ref{lnf.9} 
shows  each braids to  
type \eql{A3list1}{(i)}.  Conclude these cases with Lem.~\ref{lnf.10}. 
   
\begin{lem} \label{lnf.9} Assume $\bg\in \ni(A_n,\bfC_{3^r})$ and $i<j<k$ have these properties:  
\begin{edesc} \label{enf.13}
\item \label{enf.13a} $\lrang{g_i,g_j,g_k}$ is transitive on a five integer subset 
from 
$\{1,\dots,n\}$; and 
\item \label{enf.13b} $g_ig_jg_k$ is a 3-cycle.\end{edesc}
Then, $\bg$ braids to $\bg'$ of coalescing type \eql{A3list1}{(i)}.\end{lem}

\begin{proof} A braid from Blocks Lem.~\ref{lnf.8}, and a conjugation from Prod-one Lem.~\ref{lnf.5}   
allows assuming 
$i=1$, 
$j=2$ and $k=3$ and the integers 
of transitivity in \eql{enf.13}{enf.13a} are $\{1,2,3,4,5\}$. 
Lemma \ref{3-3braids} gives  transitivity of $B_4$ on $\ni(A_5,\bfC_{3^4})$. 
Thus, there exists $Q\in B_4$ with 
$$(g_1,g_2,g_3,(g_3g_2g_1)^{-1})Q=((1\,2\,3),(3\,2\,1),(1\,4\,5),(5\,4\,1)).$$
This gives $Q'\in B_r$ with 
$(\bg)Q'=((1\,2\,3),(3\,2\,1),(1\,4\,5),g'_4,\dots,g'_r)$:  type  
\eql{A3list1}{(i)}. \end{proof}

\begin{lem} \label{lnf.10} The position of the
5-cycle determines a 3-typle in $\ni(A_5,\bfC_{+53^2})^\inn$. Conclude, there are braids of
$\bh_1$ and
$\bh_3$ to type
\eql{A3list1}{(i)}. 
\end{lem}

\begin{proof} The first sentence says $\bg'=((1\,2\,3),(1\,4\,5),(5\,4\,3\,2\,1))$ has $B_3$ orbit 
$\ni(A_5,\bfC_{+53^2})^\inn$. Check: If
$(g_1',g_2',(5\,4\,3\,2\,1))\in \ni(A_5,\bfC_{+53^2})^\inn$, then conjugating by a power of
$(5\,4\,3\,2\,1)$ gives
$\bg'$. This completes the first sentence.  

The two cases $\bh_1$ and $\bh_3$ are similar. So we do just the harder,
$\bh_1$. Braid to
$((1\,2\,3),(1\,2\,3),(1\,4\,5),(1\,5\,6),(1\,6\,4),(1\,2\,3))$, then coalesce the 2nd and 3rd (resp.~5th and
6th) entries to get $\bw=((1\,2\,3),w_2^*, (1\,5\,6),w_4^*)$. Cor.~\ref{nf.2} shows $s(\bw)=-1$. Another
application of it shows $s((1\,2\,3\,4\,5),(1\,5\,6\,3\,2),(3\,6\,4))=-1$. 

According to Lem.~\ref{3-3braids}, this implies there is a braid of $\bw$ to 
$$\bu=((1\,2\,3\,4\,5)^*,(1\,5\,6\,3\,2)^*,(3\,4\,6),(3\,4\,6)).$$  Reexpand $\bu$ to a 6-tuple $\bh_1'$.
From  the first sentence conclude, for some $m$ and $t$, $(\bh_1')Q_1^mQ_2^t$ has 
$((1\,2\,3),(4\,5\,6))$ (resp.~$((1\,5\,6),(1\,3\,2))$) as its 1st and 2nd (resp.~3rd and 4th) entries. So,
from the 1st and 4th entries, 
$\bh_1'$ braids to type \eql{A3list1}{(i)}. 
\end{proof}

\subsection{Starting induction and $S_n$ conjugation} \label{ns.1} We prove Thms.~\ref{thmA} and \ref{thmB} by inducting
in lexicographic order on pairs
$(n,r)$, $r\ge n-1\ge 4$. Recall: $(n,r)\ge (n',r')$ if $n>n'$ or if $n=n'$ and $r>r'$.

Strong Coalescing Lem.~\ref{lnf.2'} braids anything in    
$ \ni(A_n,\bfC_{3^r})$ to $\bg=(\row g r)$ with $g_1g_2=1$. Rewrite 
$(g_3,\dots,g_r)$ as 
$(\bg'_1,\dots,\bg'_t)$ with this property.
There are $t$ disjoint orbits $\row I t$, $t\le 3$,
of $\bg'$  on $\N n$. For any orbit $I_j$ of length at least 3, $\bg_j'$
consists of the 3-cycles with support in $I_j$. So,
each $\bg'_j$ generates a transitive group on $I_j$.
It may happen $I_j$ has length 1: $\bg_j'$ would be
empty, so assuring  
with $|I_j|=n_j$ that $\sum_j n_j=n$.
Induct by assuming (on $I_j$) $\bg'_j$
has the form \eqref{ens.1a}, \eqref{ens.1b}, \eqref{ens.2a} or \eqref{ens.2b}.  The
induction relies on the nontrivial case
$n=4$.

\subsubsection{\eql{A3list1}{(i)} coalescing} 
We extend $(g,g^{-1})\mapsto 1$ in list  \eql{A3list1}{(i)} of
Lem.~\ref{lnf.2'}. 
\begin{edesc} \label{genList} 
\item \label{(i.a)}  $(n,r)\mapsto (n,r-2)$, the $r-2$ elements
$g_3,\dots,g_r$ remaining
from coalescing $(g,g^{-1})$ at the beginning of
$Q(\bg)$, are transitive.
\item \label{(i.b)}  $(n,r)\mapsto ((n_1,n_2),r-2)$ or $((n_1,n_2,n_3),r-2)$:
$\lrang{g_3,\dots,g_r}$ has orbits of cardinality $n_1,n_2$ (resp., $n_1,n_2,n_3$).
\end{edesc}

\subsubsection{$S_n$ conjugation lemma} \label{ns.3}
Suppose $\bg\in \ni( A_n,\bfC_{3^r})$ and $\beta=(2\,3)$. If you can braid conjugation by
$\beta$, then Product-one Lem.~\ref{lnf.5} produces 
$S_n$ conjugation from  braiding. For $n\ge 5$, Lem.~\ref{lns.1} extends 4-Lemma
\ref{lnt.2}. 

\begin{lem}[$S_n$-Conjugation]\label{lns.1} Assume $\bg\in \ni( A_n,\bfC_{3^r})$, 
$n\ge  5$, has the form $((g,g^{-1})^{(u)},\bg')$ where $g=(i\,j\,k)$ and neither $j$ 
nor $k$ appear in the supports of the entries of $\bg'$. Then, for any $\gamma\in S_n$, 
there is a  braid $Q$ with $(\bg)Q=\gamma\bg\gamma^{-1}$. The conclusion applies 
to all 
$r$-tuples appearing in \S\ref{nformC3r}. \end{lem}

\begin{proof} Let $\beta_0=(j\,k)$. Since $S_n=A_n \cup \beta_0 A_n$, the Product-one Lemma
gives the  conclusion if 
we show it for $\beta_0$. By hypothesis, conjugation of $\bg$ by $\beta_0$ 
produces $(g^{-1},g,\bg')$. Apply $Q_1$ to braid this back to $\bg$. 

For most $r$-tuples in \S\ref{nformC3r}, $\{2,3\}=\{j,k\}$ works, so long as the nub
(Def.~\ref{nub}) has no support in
$\{2,3\}$.   These are the only exceptions.  
\begin{itemize} \item  $((1\,2\,3),(1\,2\,3)^{-1}, 
(1\,3\,4),(1\,4\,5),(1\,5\,3))$ from 
\eqref{ens.1a}.
\item  $((1\,2\,3),(1\,2\,3)^{-1}, 
(1\,3\,4)^{(2)},(1\,3\,5),(1\,5\,4))$ from 
\eqref{ens.2a}.\end{itemize} In each case, see easily that the desired  $Q$ is the
composition of
$Q_1$ and the Blocks Lemma braid that effects conjugation by $(2\,3)(4\,5)$. 
\end{proof}

\subsection{The general induction} \label{nfo}  \S \ref{nt} treats cases $n=3$ 
and 4. Now take 
$n\ge 5$ to complete the induction for Theorems \ref{thmA} and \ref{thmB}.  Prop.~\ref{A43-2} has the initial case, $(n,r)=(5,4)$. 
\S\ref{5.1.a} does case \eql{genList}{(i.a)}. The remaining subsections handle   \eql{genList}{(i.b)}. 

\subsubsection{Case \eql{genList}{(i.a)}} \label{5.1.a} The induction
gives $Q\in B_{r-2}$ with $(g_3,\dots,g_r)Q$
a standard form for the given $r-2$ and $n$. The
Blocks Lemma then gives $Q'$ with
$$(g_1,g_1^{-1},(g_3,\dots,g_r)Q)Q'=((g_3,\dots,g_r)Q,(1\,2\,3),(3\,2\,1)),$$
standard form for $r$ and $n$.

\subsubsection{Setup for two orbit case of \eql{genList}{(i.b)}} Denote $(g_3,\dots,g_r)$ by $\bg'$  
and assume $\lrang{\bg'}$ has two orbits. Blocks Lem.~\ref{lnf.8} conjugates $\bg$ by  
$\alpha\in A_n$ so $\alpha\bg'\alpha^{-1}$ has orbits 
$\{1,\dots,n_1\}=O_1$ and $\{n_1+1,\dots,n\}=O_2$ for 
$n_1$ with $\alpha g_1\alpha^{-1}=(1\,n_1\np 1\,n_1\np 2)$.

\begin{lem}[Orbits] \label{lnfo.1}  We may
assume the following.
\begin{edesc} \item $(g_3,\dots,g_s)$ has support in $O_1$, with one integer in common
to the support of
$g_1$  ($\Pi(g_3,\dots,g_s)=1$). 
\item  $(g_{s\np1},\dots,g_r)$ has
support in $O_2$  ($\Pi(g_{s\np1},\dots,g_r)=1$) with two integers in common to the
support of $g_1$.
\end{edesc} 
\end{lem}

\begin{proof} Let $O_1$ and $O_2$ be the orbits of
$\lrang{\bg'}$. As the supports of $g_i$, $i\ge 3$, are in either $O_1$ or $O_2$,
we may braid the former to the left of the latter. This gives 
$(g_3,\dots,g_s)$ and
$(g_{s\np1},\dots,g_r)$ as in the lemma statement. Two integers in the support of
$g_1$ lie in one orbit, one in the other. 

Conjugate $\bg$ by
some $\beta'\in S_n$ to get $\bg'$ satisfying the lemma's 
conclusion. If $\beta'$ is in $A_n$, we are
done.  If not, assume some 2-cycle $\gamma'$  fixing the support
of 
$g_1$  switches two integers either in $O_1$ or in $O_2$. Then take the product of
$\beta'$ and $\gamma'$ to finish the proof. We chose 
$O_1$ to have only one integer in common with the support
of $g_1$. As $O_1$ contains at least two other integers,
choose $\gamma'$ to switch these.  
  \end{proof}

Now apply induction to
$(g_3,\dots,g_s)=\bar\bg_1$ and
$(g_1^\sph,g_1^{-1},g_{s\np1},\dots,g_r)=\bar\bg_2$ coming from Orbits 
Lem.~\ref{lnfo.1}. Put each  in 
a normal form  from \S\ref{nformC3r}. This requires special attention to the case
$|O_i|=3$ or 4: 3-Lemma \ref{lnt.1} allows braiding the endings (as in  \S \ref{ns.2})
appropriately. So assume 
$\bg=(\bar\bg_1,\bar\bg_2)$ with $\bar\bg_2$ in a normal form with 
$\{1,n_1+1,\dots,n\}$ replacing $\{1,2,\dots,n\nm n_1\np1\}$. As in   
\S \ref{ns.2}, divide $\bar\bg_2=(\bg_{2,b},\bg_{2,e})$
into beginning and end parts:  
$\bar\bg_{2,e}$ is   
$((1\,n_1\np1\,n_1\np2),(1\,n_1\np2\,n_1\np1))^{(t)}$ for some $t$. 
The next lemma braids to replace $\bg_{2,e}$ by $((1\,2\,3),(3\,2\,1))^{(t)}$. 

\begin{lem}[Tail Lemma] \label{tailLem} Assume $n\ge 5$. Let $\bh_1$ have 3-cycle
entries  transitive on $\{1,\dots,n_1\}$ and $\Pi(\bh_1)=1$. Similarly, 
$\bh_2$ has 3-cycle entries transitive on 
$\{1,n_1\np1,\dots,n\}$ and $\Pi(\bh_2)=1$. Then, there is a $Q\in B_r$ with 
$$(\bh_1,\bh_2, (1\,n_1\np1\,n_1\np2),(1\,n_1\np2\,n_1\np1))Q=(\bh_1,\bh_2, 
(1\,2\,3),(3\,2\,1)).$$\end{lem}

\begin{proof} This is an application of the Blocks Lemma. Since 
$\lrang{\bh_1,\bh_2}=A_n$, this group contains $\gamma$ conjugating 
$(1\,n_1\np1\,n_1\np2)$ to $(1\,2\,3)$. The Blocks Lemma says you can achieve 
this conjugation on the block $ ((1\,n_1\np1\,n_1\np2),(1\,n_1\np2\,n_1\np1))$ 
by an element of $B_r$ leaving the other blocks fixed.
\end{proof}

Apply the Tail Lemma to  assume $\bg$ is 
$(\bar\bg_1,\bar\bg_2, ((1\,2\,3),(3\,2\,1))^{(t)})$ with  $\bar\bg_{1,e}$ 
and $\bar\bg_{2,e}$  empty and $\bar \bg_1$ (resp.~$\bar \bg_2$) a \S\ref{nformC3r} normal form
 on its orbit $O_1$ (resp.~$O_2\cup\{1\}$). 
This gives the following principle, using the nub (Def.~\ref{nub})
of a normal form. 

\begin{princ} Let $\bar\bg_{i,\nb}$ be the nub of $\bg_i$, $i=1,2$. If we can braid
$(\bar\bg_{1,\nb},\bar\bg_{2,\nb})$ to one of the normal form nubs juxtaposed with an
H-M rep. \wsp call that a {\sl stable\/} nub \wsp then we  are  done. \end{princ} So,
to complete the proofs of Thm.~\ref{thmA} and \ref{thmB} requires two things.
\begin{edesc} \label{enfo.1}
\item \label{enfo.1a} Listing juxtapositions $(\bar\bg_{1,\nb},\bar\bg_{2,\nb})$ of
normal form nubs. 
\item \label{enfo.1b} Braiding each of the tuples in list \eql{enfo.1}{enfo.1a} to a
stable nub. 
\end{edesc} 

\subsubsection{List \ref{nfo}.A: Repeats from list \ref{nformC3r}} \label{ListA} We comment
on our naming conventions when both $(\bar\bg_{1,\nb},\bar\bg_{2,\nb})$ have the same
type. If
$n\ge 5$, then
\eqref{ens.1a}${}_+$ falls outside \eql{enfo.1}{enfo.1a}. Still, we must consider 
$|O_1|=|O_2\cup\{1\}|=3$. This does fit \eqref{ens.1a}${}_+$, though with
$n=3$. So we use that to label this case in \eqref{ens.1a}${}_+$+\eqref{ens.1a}${}_+$
below. Also, \eqref{ens.1b}${}_{\pm}$+\eqref{ens.1b}${}_{\pm}$ give the same entries as
\eqref{ens.1a}${}_{\pm}$+\eqref{ens.1a}${}_{\pm}$, so we leave them out. 

With a  natural  renaming of  integers, here are
those 
$\bar\bg_1$ and $\bar\bg_2$ with both from the same place 
in the list  \S\ref{nformC3r}. Case \eqref{ens.1a}${}_-$+\eqref{ens.1a}${}_-$ looks
odd, using $n=4$ (even) for simplicity, though theoretically
we only allowed $n$ odd. Finally, in
\eqref{ens.2b}${}_{\pm}$+\eqref{ens.2b}${}_{\pm}$, we don't use the nub (which
has 4 entries), but rather a stable nub, and then we  substitute the braid of 
Lem.~\ref{lnf.10} for that.  
 
\begin{description} 
\item[\eqref{ens.1a}${}_+$+\eqref{ens.1a}${}_+$]
$((1\,2\,3)^{(3)},(1\,4\,5)^{(3)})$
\item[\eqref{ens.1a}${}_-$+\eqref{ens.1a}${}_-$] 
$((1\,2\,3),(1\,3\,4),(1\,4\,2),(1\,5\,6),(1\,6\,7),(1\,7\,5))$
\item[\eqref{ens.2a}${}_-$+\eqref{ens.2a}${}_-$]  
$((1\,2\,3)^{(3)},(1\,3\,4),(1\,4\,5),(1\,5\,4),(1\,6\,7)^{(3)}$,\hfil\break 
$(1\,7\,8),(1\,8\,9),(1\,9\,7))$
\item[\eqref{ens.2b}${}_+$+\eqref{ens.2b}${}_+$] $((1\,2\,3),(1\,3\,2), 
(1\,3\,4),(1\,4\,3), 
(1\,5\,6),(1\,6\,5),(1\,6\,7),(1\,7\,6))$
\item[\eqref{ens.2b}${}_-$+\eqref{ens.2b}${}_-$] 
$((1\,2\,3)^{(3)},(1\,4\,5),(1\,5\,6),(1\,6\,4), 
(1\,7\,8)^{(3)}$, \hfil\break $(1\,9\,10),(1\,10\,11),(1\,11\,10))$
\end{description}

\begin{lem} Each juxtaposition in the list above braids to a stable nub.
\end{lem}

\begin{proof} 
The Blocks Lemma braids each of \eqref{ens.2a}${}_-$+\eqref{ens.2a}${}_-$ and 
\eqref{ens.2b}${}_{\pm}$+\eqref{ens.2b}${}_{\pm}$ to a juxtaposing of  types 
\eqref{ens.1a}${}_{\pm}$+\eqref{ens.1a}${}_{\pm}$. Lems.~\ref{lnf.9} and 
\ref{lnf.10} braids these respectively to 
H-M reps. We are done. 
\end{proof} 

\subsubsection{List \ref{nfo}.B: Distinct pairs from 
\S\ref{nformC3r} and concluding the proof} Use the same principles as in
\S\ref{ListA}. Example:  \eqref{ens.1a}${}_+$+\eqref{ens.1a}${}_-$ has the shape 
$$((1\,2\,3)^{(3)},(1\,4\,5),(1\,5\,6),(1\,6\,4)).$$  Lem.~\ref{lnf.10}
braids this to a stable nub from \eqref{ens.2b}. 

The situation from
\eqref{ens.1a}${}_+$+\eqref{ens.2b}${}_-$ is slightly different: 
$$((1\,2\,3)^{(3)},(1\,4\,5)^{(3)},(1\,5\,6),(1\,6\,7),(1\,7\,5)).$$ The juxtaposition
$((1\,2\,3)^{(3)},(1\,4\,5)^{(3)})$ from
\S\ref{ListA} braids to an H-M rep., so we are done. In a like manner, we find 
there are no serious new cases.   

Finally, the Orbits Lemma has a simple variant when 
$\bg'=(g_3,\dots,g_r)$ has three orbits. This supports the arguments following 
it to produce the same kind of lists. 

\section{Applications to $G_\bQ$} \label{nth}
Our applications of  Theorems \ref{thmA} and \ref{thmB} --  
\S\ref{nth.1} and \S\ref{maxAltExt} -- are to the Inverse Galois Problem. We use the \hc spaces $\sM_{g,\pm}$ of \S\ref{hcfamilies} 
in \S\ref{maxAltExt}. Given a field $K$, 
$\bar K$ indicates an algebraic closure, and $G_K$ its absolute Galois group.  

\subsection{$(A_n,A_n)$ and $(A_n,S_n)$ realizations} \label{nth.1} \S~\ref{secrigp-aigp} explains  $(G,\hat G)$ regular
realization: $(A_n,A_n)$ and $(A_n,S_n)$ realizations are a special case. The {\sl no centralizer condition\/} holds for the
standard representation of
$A_n$. So Prop.~\ref{rigp-aigp} says $K$ points on corresponding Hurwitz spaces correspond to finding $K$ realizations
by covers in the Nielsen class
$\ni(A_n,\Ct r)$. Thm.~\ref{thmA} shows  the spin lifting invariant \eqref{spinlift} determines the components, so \cite[Thm.~3.16]{Fr4} 
gives part one of Cor.~\ref{corAr1}. The short proof for part two  is from Hilbert's Irreducibility Theorem (HIT).  

\begin{cor} \label{corAr1} For $n\ge 5$, each component of $\sH(A_n,\bfC_{3^r})^\inn$ or
$\sH(A_n,\bfC_{3^r})^\abs$ (with its map to $U_r$) has definition field $\bQ$. Further,
let $\sH^*$ be a  component of $\sH(A_n,\bfC_{3^r})^\abs$. Then, for a dense set
of  $\bp\in
\sH^*(\bar \bQ)$, the corresponding cover 
 gives an $(A_n,S_n)$ realization over $\bQ(\bp)$. \end{cor} 

\begin{proof} We need only show the last sentence. Let $\sH^{**}$ be the (unique according to Thm.~\ref{thmA}) component of
$\sH(A_n,\bfC_{3^r})^\inn$ lying over $\sH^*$. So $\sH^{**}\to
\sH^*\to U_r$ is a sequence of absolutely irreducible covers over $\bQ$, with $\sH^{**}\to
\sH^*$ Galois with group $\bZ/2$. (The cover $\sH^*\to U_r$ is far from Galois.) According to HIT (say, \cite[Chap.~10]{FrJ}, or
\S\ref{maxAltExt}), there is a dense set $\bz\in U_r(\bQ)$ so that for any point $\hat \bp\in \sH^{**}$ over $\bz$, $[\bQ(\hat \bp):\bQ]$ is
the degree of
$\sH^{**}/U_r$. Conclude: If $\bp$ is the image of $\hat \bp$ in $\sH^*$, then $[\bQ(\hat \bp):\bQ(\bp)]=2$, precisely what we need for 
an $(A_n,S_n)$ realization. \end{proof}

When 
$r=n-1$,
\cite{Me} shows there are an infinity of  $(A_n,A_n)$ realizations
from $\bQ$ points in $\sH(A_n, \bfC_{3^{n-1}})$ whose images are dense in $U_{n-1}$. Thm.~\ref{thmA} 
shows there is only one component of this moduli space. Conclude the following. 
\begin{cor} \label{corAr2}  When $r=n-1$,  $\bQ$ points giving  $(A_n,A_n)$ realizations
are analytically dense in $\sH(A_n, \bfC_{3^{n-1}})$. \end{cor}

\begin{rem}[Comments on Cors.~\ref{corAr1} and \ref{corAr2}] \label{Ar1Ar2comm} There is more detail on these corollaries in
\cite[p.~163--167]{FrK} to help a reader use the Hurwitz space  interpretation. This includes the effect of special assumptions on 
Mestre's parametrization. 
\end{rem}
 
\subsection{The maximal alternating extension of $\bQ$} \label{maxAltExt} We say $K$ is {\sl Hilbertian\/} if it satisfies {\sl HIT}:
Any irreducible $f\in K[x,y]$ of degree at least 2 in $y$ remains irreducible (over $K$) for $\infty$-ly many specializations of $x\mapsto
x_0\in K$.  

Call $K$ {\sl projective\/} if $G_K$ is projective: Surjective homomorphisms to
$G_K$ split. \S\ref{hpconj} puts a particular field $\bQ^\alt$ in the context of \cite{FrV3}. 

\newcommand{\hyp}{\text{\bf Hyp}}
\subsubsection{Hilbertian$+$Projective $\implies$ Pro-free Conjecture} \label{hpconj}  
Pro-free groups are projective, though most projective (even finitely generated) groups are not pro-free (see Ex.~\ref{proj-not-free}). 
All  results from \cite{FrJ} are in both editions (we use the 2nd). 

Let $I\le 
\bN^+$ be any infinite set of integers. Define  $\bQ^{\alt,I}$ to be the composite of all Galois extensions
$L/\bQ$ with group $A_n$ for some integer $n\in I$. \cite[p.~476]{FrV3} asks if $\bQ^{\alt,\bN^+}\eqdef \bQ^\alt$
is P(seudo)A(lgebraically)C(losed). From \cite[Thm.~10.4]{FrJ} this is equivalent to every absolutely irreducible curve over $\bQ$ having  a
rational point in $\bQ^{\alt}$. You can as well ask if any of the 
$\bQ^{\alt,I}$ are PAC. 

Recall: $\tilde F_\omega$ is the profree group on countably
many generators.

\begin{thm} If $\bQ^{\alt,I}$ is PAC, then there is a natural short exact sequence \begin{equation} \label{knownpres} 1\to \tilde F_\omega\to
G_\bQ\to G(\bQ^{\alt,I}/\bQ)\to 1.\end{equation}  \end{thm} 

\begin{proof}  Apply \cite[Thm.~A]{FrV3}: A  PAC, Hilbertian subfield of
$\bar \bQ$ has pro-free absolute Galois group. As $G(\bQ^{\alt,I}/\bQ)$ is  the product over
$n\in I$ of an infinite number of $A_n\,$s, it is a nontrivial finite extension of a Galois extension of $\bQ$. It is
automatic from
\cite[Prop.~13.9.4]{FrJ} that such a field   must
be Hilbertian. \end{proof} 

\newcommand{\FV}{\text{\rm FV}}

\cite{FrV3} presents $G_\bQ$ \wsp as \eqref{knownpres} would
\wsp  
 by known groups: Products of $S_n\,$s (instead of $A_n\,$s) on the right, $\tilde F_\omega$ on the left. It does this by producing PAC
fields $\bQ^\FV$ that are a composite of disjoint $S_n$ extensions.  

\begin{guess} \label{FrV} \cite{FrV3}  conjectures for $K\le \bar \bQ$ that Hilbertian $+$ Projective $\implies$ $G_K$
is profree. \end{guess} 

\begin{exmp}[Projective, but not pro-free]  \label{proj-not-free} For $G$  a finite group, consider its minimal
projective ({\sl universal Frattini}) cover  $\tilde \phi: \tilde G\to G$. It has pronilpotent kernel that is a 
direct product of  pro-free pro-$p$ groups (one for each prime $p$ dividing  $|G|$) \cite[\S 22.11]{FrJ}. Write
$\tilde G$ as the fiber product of group covers  $\tG p\to G$, with pro-$p$ kernel, $p||G|$. 

Assume $G$ is not cyclic, so neither is $\tilde G$. Suppose at least two primes divide $|G|$. Then, if $\tilde G$ is pro-free it
has rank at least 2. Yet, it has a finite index subgroup that is a {\sl product\/} of two nontrivial proper closed subgroups.
Therefore that subgroup is not pro-free, contradicting Schreier's theorem \cite[Prop.~17.6.2]{FrJ}.  
\end{exmp}

\subsubsection{Cohomological observation on $\bQ^\alt$ being projective}  For 
\eqref{knownpres} to hold,  
$\bQ^{\alt,I}$ must be projective. \cite[Cor.~p.~81]{SeGalCoh} says any totally imaginary $K/\bQ$, with each of its non-archimedian
completions of infinite decomposition degree over the completion of $\bQ$, must be projective. To see projectivity of $G_K$, consider each 
finite extension  
$L/\bQ$  and any Brauer-Severi variety $X$ over $L$. 
Inflating the corresponding Brauer group element from $L$ to $L\cdot K$ produces (torsion) generators of the Brauer group of all finite
extensions of $K$. From Class Field Theory, you determine
$[X]$ from its $\bQ/\bZ$ values at the completions of $L$. Inflating
$[X]$  from each completion $\tilde L$ of
$L$ to   
$K\cdot \tilde L$   kills $[X]$ (multiplication
by the degree gives the inflation). So, the inflation of $[X]$ to $K\cdot L$ is trivial $\Leftrightarrow$ $[X]$ has a $K\cdot L$ point. This 
statement on {\sl certain\/} varieties having rational points for each such $L$ is is equivalent to projectivity. 
It is, however, much weaker than $K$ is PAC.

This projectivity criterion should work to show $\bQ^\alt$ is projective, because  there are many known
Galois extensions of
$\bQ$ with group an alternating group. 
Yet, \cite{FrV3} skipped showing projectivity for the fields I labeled above as $\bQ^\FV$. Rather, it directly used  
\cite{FrJ78} to give each $\bQ$ curve a simple-branched $\bQ$ covering of $\prP^1_z$. 

\begin{prop}  This alternating group analog implies  $\bQ^{\alt,I}$ is PAC:  
\begin{triv} \label{Ancond} Each projective nonsingular curve $X/\bQ$ appears in a cover \\ $\phi:X\to\prP^1_z$ 
that gives an $(A_n,A_n)$ realization over $\bQ$ for some $n\in I$. \end{triv} \end{prop} 

\cite{Ne} has a well-known  result: Two
number fields with isomorphic absolute Galois group are conjugate. This result uses class field theory to conclude from 
valuation theory that  abelian extensions determine the correspondence between primes. Let $K$ be a number field, and denote  
 the composite of all Galois extensions of $K$ with group a Frattini cover of an alternating group extension by $\tilde K^\alt$.  

\!\begin{quest} Does $G(\tilde K^\alt/K)$ determine $K$ up to conjugacy. \end{quest} 

\!\!\subsubsection{Restricting condition \eqref{Ancond} to odd order branching} Restrict to $g\ge 2$. \cite[p.~36]{Mu76} discusses that $\sM_g$
is not {\sl unirational\/} (that is, the image of a map from some projective space)  for large
$g$. Yet, a    unirationity conclusion holds: \lq\lq [$\sM_g$] has lots of rational curves:\rq\rq copies of affine subsets of
$\prP^1$. These come from any algebraic surface 
$Z$ and a meromorphic function $f: Z\to \prP^1_z$ with general fiber having genus $g$. So, possibly $\sM_g$ still has sufficiently many 
 rational curves. 

\begin{quest} \label{hitprop} Let $U\to W$ by any irreducible $\bQ$ cover (finite and flat; so surjective \cite[Chap.~2 \S 7,
Prop.~4]{Mu66}) with $W$ open in $\sM_g$. Is   there a $\bQ$ rational curve $X\subset W$ where restriction of $U$  over $X$ remains
irreducible? 
\end{quest} 
  
 We don't know if \eqref{Ancond} is true. Prop.~\ref{getAnAn} shows there are curves over $\bQ$ for with no $(A_n,A_n)$
or an $(A_n,S_n)$ realization (over $\bQ$), for any $n$, from  odd order branched covers (as in 
\S\ref{useDiff}).  The proof would  simplify if    
Quest.~\ref{hitprop} had a yes answer; even restricting  $U$ to be  one of the covers 
$\sM_{g,\pm}$ (Prop.~\ref{hcspaces}). 

\begin{lem} \label{findirred} Assume $V\subset \sM_g$ is a $\bQ$ subvariety of dimension at least 1
satisfying the following conditions. 
\begin{edesc} \label{hccond} \item \label{hcconda}  There is a generically surjective $\bQ$ morphism $W\to V$ with the function field
$\bQ(V)$ of $V$ algebraically closed in
$\bQ(W)$. 
\item \label{hccondb}  $W$ is birational to an open subset of projective space $\prP^N$ (for some $N$). 
\item \label{hccondc}  Restricting 
$\sM_{g,\pm}$ to
$V$ has no $\bQ$ components of degree 1 over $V$. \end{edesc} Conclude: A set of curves $X/\bQ$ of genus $g$, corresponding to
$\bQ$ points dense in
$V(\bQ)$, have no $(A_n,A_n)$ or $(A_n,S_n)$ realizations (over $\bQ$) with odd order branching. 
\end{lem}  

\begin{proof}  Suppose $X$, over $\bQ$,  of genus $g$, has $\phi:X\to \prP^1_z$, over $\bQ$,  with odd
order branching. Lem.~\ref{lnth.2} gives  a 
$\bQ$ \hc\ class on $X$.  To conclude we find a dense set of $v\in V(\bQ)$ so the corresponding  $X_v\,$s have no
$\bQ$ \hc\ class.  

Prop.~\ref{hcspaces} says the connected spaces $\sM_{g,\pm}$ have respective degrees $2^{2g-1}\pm
2^{g-1}$  over
$\sM_g$. Denote by $V^{\pm}$ the restriction of each of these covers over $V$. 

By condition \eql{hccond}{hcconda}, over a Zariski open subset $V^*\le V$ the pullback map
$\pr_W: V^{\pm}\times_V W \to  W$ is a cover whose $\bQ$  irreducible components $W'$ have the same degrees over $ W$ 
as the corresponding components of the covers $V^{\pm} \to V$. Further, from \eql{hccond}{hccondc}, none of those degrees is 1. 

From
\eql{hccond}{hccondb}, we can apply Hilbert's Irreducibility Theorem (as in Cor.~\ref{corAr1}) to the collection of covers $\pr_W: W'\to
W$, above. So, there is  a dense set of
$\bp^*\in W(\bQ)$ with no
$\bQ$ point above them in {\sl any\/} of the $W'\,$s. Conclude: The image $\bp\in V(\bQ)$ of $\bp^*$ gives a curve
over
$\bQ$ with no \hc\ class  over $\bQ$. This is contrary to the above,  finishing the proof of the lemma. \end{proof}

We apply Lem.~\ref{findirred} to the hyperelliptic locus $\hyp_g$: genus $g$ curves with a degree 2 map to
$\prP^1$. For $g\ge 2$, a hyperelliptic curve is determined by the branch points of its canonical map to $\prP^1_z$, up to the action of
$\PGL_2(\bC)$ (M\"obius transformations) on these unordered branch points. See \S\ref{hurSpaces} for the $\PGL_2(\bC)$ action. So, with
$U_r$ as in
\S\ref{gpNC},
$\hyp_g$ is 
$U_r/\PGL_2(\bC)$  where
$g=r/2-1$. 

\begin{prop} \label{getAnAn} For each even $g\ge 2$, the space $\hyp_g=V$  satisfies the hypotheses  in \eqref{hccond}. Dense
in
$\hyp_g$ is a set of
$\bQ$ curves $X$ where $X$ has no  presentation as a 
$\bQ$ cover $\phi:X\to\prP^1_z$ with odd order branching. Thus, such a curve has no odd order branched cover over $\bQ$ fulfilling
\eqref{Ancond}.  
\end{prop} 

\begin{proof} 
The components $\hyp_{g,{\pm}}$ from restricting $\sM_{g,\pm}$ over $\hyp_g$ correspond to
the orbits of the fundamental group of $\hyp_g$ restricted to the action on \hc\ classes. Prop.~\ref{specreps} reveals 
this action by hand. Only when $g$ is odd is there in this restriction a component of degree 1. 

Consider the $\bQ$ map $\mu_r: U_r\to U_r/\PGL_2(\bC)$, $r\ge 4$. \!From \!\cite[Prop.~6.10]{BaFr}, over a Zariski open subset $U'\le
U_r/\PGL_2(\bC)$ the fibers of $\mu_r$ consist of copies of $\prP^3$, a variety with pure transcendental function field. This
gives  the pullback maps
$\hyp_{g,{\pm}}\times_{U_r/\PGL_2(\bC)} U_r \to  U_r$ the following property. Over a Zariski open subset $U^*\le U_r$, 
the absolutely irreducible component covers of   $$\hyp_{g,{\pm}}\times_{U_r/\PGL_2(\bC)} U^*\eqdef \hyp{}^*_{g,{\pm}}\to
 U^*$$  have respective degrees listed in 
Prop.~\ref{specreps}.  In all cases,  the covering degrees are distinct, each component is over $\bQ$, and when
the genus is even the degrees are at least $r$. Then the hypotheses of Lem.~\ref{findirred} apply and for a
dense set of hyperelliptic curves over $\bQ	$,  the proposition conclusion  holds.  
\end{proof}

\newcommand{\bl}{{\pmb \ell}} 
\section{$\frac 1 2$-canonical divisors and $\theta$-functions} \label{halfclasses} 
\S\ref{useDiff} explains how the irreducible
components 
$\sH_{\pm}(A_n,\bfC_{3^r})^{*,\red}$, \lq\lq$*=\abs$\rq\rq\ or \lq\lq$\inn$,\rq\rq\  support \hc\ classes,  and  how these then support the
analytic continuation of close-to-canonical $\theta$-functions.  The difference between the two cases $\pm$: When $r$ is even (resp.~odd) the $\theta\,$s 
for + are even (resp.~odd), for -  odd (resp.~even) in the
$\theta$ variables. \S\ref{oddord} then discusses the even natural $\theta$-nulls on the appropriate components.  The key issue is that these be
non-zero. At present we can only prove this for a given value of $g=r-(n-1)\ge 1$ for infinitely many $(r,n)$. \S\ref{hyplocus} 
computes components of 
$\sM_{\pm}$ over the hyperelliptic locus.

\subsection{Well-defined \hc\ classes} \label{useDiff} Let  
$\Phi: 
\sT\to \sH\times \prP^1_z$ be a family of covers with odd order branching. That is, for $\bp\in \sH$,  $\Phi$ restricts over the fiber
$X_\bp$ of
$\sT\to \sH$ over $\bp$  to $\phi_\bp: X_\bp \to \bp\times \prP^1_z$ with odd order branching.   
Lem.~\ref{def-hc} shows this defines a unique \hc\ divisor (from $\Phi$) on $X_\bp$. Then, Lem.~\ref{lnth.2} says, if $\Phi$ is a
 Hurwitz family, its {\sl reduced\/} Hurwitz family (\S\ref{rnz.4}) supports a \hc\ divisor class, and so a well-defined $\theta$ divisor
$\Theta_\bp$ at each $\bp\in \sH$. Then, for a fixed $\bp_0$, Prop.~\ref{thetaChoices} gives an expression for the effect on a $\theta$
function attached to
$\Theta_{\bp_0}$  after it has been analytically continued around a closed path based at $\bp_0$.
 
\subsubsection{Using differentials} \label{halfdiffs} Consider 
$\phi: X\to\prP^1_z$ branched over $\bz=\{\row z r\}$ and having $\bg=(\row g 
r)\in G^r$ as branch cycles. 

On $X$ there is a divisor class $\kappa$ that is completely 
canonical, being the divisor class of all global meromorphic differentials on $X$. 
\begin{edesc} \label{enth.1}
\item \label{enth.1a}  Any automorphism of $X$, in its extension to the collection of degree $2g-2$ divisor classes on $X$,  fixes $\kappa$.
\item \label{enth.1b} If  $X$ is a fiber in a family $\sX\to \sP$, then
coordinates  for 
$\sP$ 
suffice as coordinates locally describing $\kappa_\bp$ for $\bp\in \sP$. 
\end{edesc} 

Still, there is not usually a way to pick one representative {\sl divisor\/} for 
$\kappa$  explicitly. This is one difference between a general family 
of  curves and a family of $\prP^1_z$ covers. A member $\phi_\bp:X_\bp\to\prP^1_z$ does give such a 
canonical class 
divisor as the divisor of the differential $d\phi_\bp$. We accept that as part of the given data. Also, when all the branch 
cycles $\bg$ (\S\ref{gpNC}) have odd order \wsp $\phi_\bp$ has {\sl odd order\/} branching, this canonically produces a \hc\ divisor.
Here is how. 

In a neighborhood $N_{x_0}$ of  $x_0\in X$, there is a one-one function $x: N_{x_0}\to \prP^1_{x^*}$ and an integer $e$ so the following
holds. With $x_0^*$ the image of $x_0$ under $x$, 
$\phi$ composed with $x^{-1}$ (functional inverse) looks like $x^*\mapsto (x^*-x_0^*)^e + z_0=z$ or
$1/(x^*)^e$ (corresponding to $z_0=\infty$) in the image of $N_{x_0}$ under $x$. Here $e$ is the ramification index of $x_0$ over $z_0$: 
the length of a corresponding disjoint cycle in $g_i$ if $z_0=z_i$. Thus, $dz$ has order $e-1$ (resp.~$e+1$) at $x_0$ between the two cases
$z_0\in \bC$ (resp.~$z_0=\infty$). So,  each $e$ being odd, implies $e\pm 1$ is even.

\begin{lem} \label{def-hc} If $\phi=\phi_\bp$ has odd order branching, then the divisor $(d\phi)$ of the meromorphic differential $d\phi$ is 
$2D_\phi$ with $D_\phi$ a well-defined divisor on $X$. If $\phi$ has definition field a perfect field $K$, then $D_\phi$ does, too,
and so does its divisor class. \end{lem} 

The divisor $D_\phi$ in Lem.~\ref{def-hc} is a well-defined \hc divisor. Any divisor {\sl class\/} $\iota$ on $X$ with $2\cdot \iota=\kappa$
is called a \hc (divisor) class. There are $2^{2g}$ of these, differing pairwise by some 2-division point on the Jacobian $\Pic^{(0)}(X)$
(identified with divisor classes of degree 0 by Abel's Theorem,  \S\ref{tor1}). 

Denote the set of \hc classes by
$S_{\kappa/2}(X)$. They 
 canonically sit in $\Pic^{(g-1)}(X)$, the degree $g-1$ divisor classes. \cite{Se3} quotes \cite{atiyah} and \cite{Mu71} for basics
on $S_{\kappa/2}(X)$. Closest to our start is \cite{Fay}, for that works with moduli spaces of curves as do we. Still, we
switch to
\cite{Sh98} of necessity, for the production of automorphic functions, for that works with global moduli as do we, though our spaces are
reduced Hurwitz spaces, not (say) Siegel upper half-spaces. 

There are two types of  $y\in
S_{\kappa/2}(X)$. Assume the divisor $D$ represents $y$. Let $L(D)$  be the linear system of
meromorphic functions $f$ on $X$ satisfying $(f)+D\ge 0$. Call $y$ {\sl even\/} (resp.~{\sl odd}) if $\dim_\bC L(D)=\dim(y)$ is even (resp.
odd).

\subsubsection{Odd order branching and reduced equivalence} Continue with the family $\Phi$ of \S\ref{halfdiffs}. 
Assume some Nielsen class $\ni(G,\bfC)$ with odd order conjugacy classes defines $\Phi$. Then Lem.~\ref{def-hc} 
smoothly assigns a well defined \hc divisor on the
Riemann surfaces attached to points of $\sH(G,\bfC)^\abs$ (or $\sH(G,\bfC)^\inn$). 

Such 
representing divisors, however, disappear if we use {\sl reduced\/} equivalence of covers (\S\ref{rnz.4}; equivalencing by a $\PGL_2(\bC)$
action). Still, even the reduced Hurwitz spaces $\sH(G,\bfC)^{\abs,\rd}$ (and $\sH(G,\bfC)^{\inn,\rd}$)  carry well defined  
\hc classes. 

\begin{lem} \label{lnth.2}  Assume $\sH$ parametrizes a family of covers in a reduced Nielsen class
$\ni(G,\bfC)^{*,\rd}$ (as above) with odd order classes.  This canonically defines \\ $\bp\in \sH \mapsto$
a  \hc\ class on the associated curve $X_\bp$.  
\end{lem}

\begin{proof} Lem.~\ref{def-hc} produces a unique \hc divisor $D_\phi$ on  $\phi: X\to \prP^1_z$ representing the
Nielsen class before it is reduced. We have only to show the divisor attached to $\alpha\circ \phi$ for $\alpha\in
\PGL_2(\bC)$ is linearly equivalent to $D_\phi$. 

Replacing $\phi$ by $\alpha\circ\phi$ produces  $d(\alpha\circ\phi)$ 
as  the differential. With no loss assume an element in $\SL_2(\bC)$ represents $\alpha$. If $\alpha(z)=\frac{az+b}{cz+d}$, then 
$$d(\alpha\circ\phi)=d\phi/(c\phi +d)^2.$$ 
Thus, this has the same divisor as does $d\phi$ with the subtraction of two 
times the 
divisor of the function $c\phi+d$. Therefore, the \hc\ class is well-defined. 
 \end{proof}

\begin{exmp}[non-odd order branching and \hc\ classes] \label{hypbranch} For $g>1$ a \hc class, as in Lem.~\ref{lnth.2}, on the family of
a Nielsen class might be rare, unless the branch cycles have odd order. Still, Prop.~\ref{specreps} shows for 
$\ni(\bZ/2,\bfC_{2^{2s}})$ with  $g=\frac{2s-2}2$ odd, and $s\ge 2$, there is a globally defined  \hc\ class on the Hurwitz spaces. If,
however, $s\ge 4$, this class gives a degenerate $\theta$ (Rem.~\ref{hypcomps}).   
\end{exmp} 

\subsubsection{\hc classes in families of covers} \label{hcfamilies} 

Consider any smooth family $\Psi: \sT \to \sH$ of genus $g\ge 1$ curves.  There is a natural map $\Psi_{\sH,\sM_g}: \sH \to
\sM_g$, by $\bp \mapsto [\sT_\bp]$. More generally, we can define $\Psi_{\sH,\sM_g}$ if  $\sH$ is just a  
{\sl stack\/} of compact Riemann surfaces in the Grothendieck topology of finite covers (see Rem.~\ref{stack}).

Also, for any
integer
$k$, there is a canonical family
$\Psi_k:
\sP^{(k)}\to 
\sH$ with  the fiber 
$\sP^{(k)}_\bp$ over $\bp$  the variety $\Pic^{(k)}(\sT_\bp)$ of degree $k$ divisor classes on $\sT_\bp$. For $\sH=\sM_g$ (the moduli
space of projective non-singular curves of genus $g$ as in \S\ref{quickReview}), this defines a cover $\sM_{g,\pm} \to \sM_g$, with the fiber
of
$\sM_{g,\pm}$ over $m\in \sM_g$ consisting of the  $2^{2g}$ points 
$S_{\kappa/2}(X_m)\subset  \sP^{(g\nm 1)}$ (\S\ref{halfdiffs}). 

This has 
disjoint irreducible components $\sM_{g,+}$ and $\sM_{g,-}$. Prop.~\ref{hcspaces}  distinguishes, on the  fiber product
$\sH\times_{\sM_g}
\sM_{g,\pm}\eqdef
\sH_{\pm}$, the points on the two components. We introduce the divisors $\Theta$  on  $\Pic^{(0)}(\sT_\bp)$. These pull  back to its
universal covering space, $\widetilde {\Pic^{(0)}}(\sT_\bp)$, where Riemann's $\theta$  functions live.

\begin{prop} \label{hcspaces} Each  $\bp^+\in \sH_{+}$ ($\bp^-\in \sH_-$) lying over $\bp\in \sH$ corresponds (uniquely) to the divisor
$\Theta_{\bp^+}$ (resp.~$\Theta_{\bp^-}$) of  zeros of one of the
$2^{2g-1}+ 2^{g-1}$ (resp.~$2^{2g-1}- 2^{g-1}$) even (resp.~odd) $\theta$ functions (as defined by Riemann, up to an exponential factor) on
$\sT_\bp$.  Given an  even (resp.~odd) $\theta_{\bp_0}$ at 
$\bp_0\in \sH$,  there are normalizations that give a unique analytic  continuation of it to even (resp.~odd) $\theta\,$s along any path
in 
$\sH$ based at $\bp_0$.  
\end{prop}  

The next two subsections do the proof. \S\ref{tor1} describes the $2^{2g}$ $\theta\,$s (even
and odd) attached to a Riemann surface. Our  aim to get $\theta\,$s to depend on
just the coordinates describing those families. That isn't possible, as Riemann knew, though for 
families of say, Thm.~\ref{thmB}, the coordinates for the $\theta\,$s are especially good. 

\S\ref{tor1} introduces coordinates from the infinite degree 
Torelli space cover of $\sM_g$. That cover has the $2g\times 2g$ symplectic group over $\bZ$, $\Sp_{2g}(\bZ)$, as its
monodromy group. Then, \S\ref{tor2} uses the finite cover  
$\sM_{g,\pm}\to
\sM_g$.  I tie this not-easily found classical result to the telegraphic discussion in
\cite{Fay}. 

\newcommand{\pomega}{\pmb \omega}

\subsubsection{Torelli space coordinates} \label{tor1} I now explain Torelli space
$\sT_g$, an unramified covering of
$\sM_g$. For $m\in \sM_g$, the fiber $\sT_{g,m}$ over $m$ consists of all possible canonical
(first) homology bases for
$X_m$. Typical notation has such a basis  as $\bl=(\row {\alpha} g, \row \beta g)$ with the matrix of cup-product
intersections looking like
$J_{2g}=\smatrix {0_g} {I_g} {-I_g} {0_g}$ using the $g\times g$ zero, $0_g$, and identity, $1_g$,  matrices. With
$\bM_g(\bZ)$ ordinary $g\times g$ matrices in $\bZ$, that leaves possible choices (see the precise notation of
\S\ref{famCovers}) of basis as a homogenous space for  
\begin{equation} \label{SP2g} \Sp_{2g}(\bZ)=\bigl\{U=\smatrix A B C D\mid A,B,C,D\in \bM_g(\bZ), U J_{2g} U^\tr=J_{2g}\bigr\}.\end{equation} 
Denote  the
$\bZ$ module  that
$\bl$ generates by $\lrang{\bl}$.  

Giving $\bl$ fixes an identification of $\Pic^{(0)}(X)$ with $\bC^g/\lrang{\Pi(\bl)}$ with $\Pi(\bl)$ a lattice in $\bC^g$, thereby starting
Riemann's generalization of Abel's Theorem. I explain $\Pi(\bl)$. 

Choose a basis
$\pomega=(\row
\omega g)$ of holomorphic differentials on
$X_m$, using one of the normalizations also typical in the literature. For example, \cite[p.~3]{Fay} takes the integral of the $g$-tuple of
differentials along the paths $\row
\alpha g$  to  be the $g\times g$ matrix $2\pi i I_g$, while others, respectively, replace $\row \alpha g$ and $2\pi i I_g$ with $\row \beta
g$ and
$I_g$. Riemann showed such a choice determines 
$\pomega$ and the resulting $g\times {2g}$ matrix has columns  listing the transpose of $(\pmb \omega)$ integrated
in order along $\row {\alpha} g, \row \beta g$. Those columns form the matrix  $\Pi(\bl)=\bigl(2\pi i I_g|\tau_\bl\bigr)$, with $\tau_\bl$ 
symmetric. Then, (Riemann showed) 
$\lrang{\Pi(\bl)}$ is the  lattice those columns span.  

Now we  require (noncanonical) choices: Choose a
set of
$g$ points
$(\row {x'} g)$ on
$X$, so the complex vector
\begin{equation} \label{jacInv} (\row {\sum_{i=1}^g 
\int_{x_i'}^{x_i}\omega} g)=\sum_{i=1}^g(\row {
\int_{x_i'}^{x_i}\omega} g)=\bw\in \bC^g\end{equation} represents a degree 0 divisor class $[D=\sum_{i=1}^g (x_i-x_i')]$  in 
$\bC^g/\lrang{\Pi(\bl)}$. Modding out by $\lrang{\Pi(\bl)}$ assures  an integral independent of the path choices from $x_i'$ to
$x_i$. 

Some discussions choose $\row {x'} g$ to be the same point repeated $g$ times; still not canonical. (As 
$\bz$ appears in this paper from coordinates on $\prP^1_z$, we use
$\bw$ rather than the traditional $\bz$.)  With $\tr$ meaning transpose, Riemann's $\theta$ function,    
\begin{equation} \label{rtheta} \theta(\bl,
\bw)=\sum_{\bn\in \bZ^g} e^{\bn\tau_\bl\bn^\tr + \bn \bw^\tr} \text{ \cite[p.~1]{Fay} on
$\sT_{g,m}\times
\bC^g$},\end{equation} is invariant under $\bw\mapsto -\bw$ (it is {\sl even}), and  under translating  $\bw$ by its $\balpha$
periods. 

Significantly, \eqref{rtheta} depends on  $\bl$; even its  divisor of zeros -- denote this $\Theta_\bl$ -- depends on 
 $\bl \!\!\mod 2\bl$.  Giving
$m\in \sM_g$ does not canonically give $\bl$, so I comment on that dependence now: How can we compare info on $X_m$ varying in a
{\sl given\/} family (a Hurwitz space, or $\sM_g$)  with variation of
$\theta(\bl,
\bw)$. The quotient of $X_m^k$ by the symmetric group $S_k$, permutating its coordinates, gives the degree $k$
divisors. 

First, we compare the zero set $\Theta_\bl$ on $\bl\times \bC^g$ with the positive, degree
$g-1$, divisor classes on $X_m$:  $W_{g-1}\eqdef W_{g-1,m}$.  Fundamental Fact: $W_{g-1}$ is birational to $X_m^{g-1}/S_{g-1}$, degree $g\nm
1$ divisors. If $\delta$ is a  
\hc\ divisor class on $X_m$, then the translate $W_{g-1,m}-\delta$  by
$\delta$ is a
$g\nm 1$-dimensional set in $\Pic^{(0)}(X_m)$. 

\begin{thm} \label{eventhetanulls} Pullback of $W_{g-1,m}-\delta$ from $\Pic^{(0)}(X_m)$ to its universal covering space $\bl\times
 \bC^g$ gives a divisor of the form $\Theta_\bl+\mu$ with $\mu=\mu_m$ representing a point of order 2 on $\bC^g/\lrang{\Pi(\bl)}$ (2-division
point). 

As $\mu$ runs over 2-division representatives, translates, $\theta(\bl,\bw+\mu)$, of $\theta(\bl,\bw)$ run over a collection of $2^{2g}$
functions, each  either even or odd. Each has zero divisor of form $W_{g-1,m}-\delta$, for some half-canonical class
representative $\delta$, and $\theta(\bl,\bw+\mu)$ is zero at $\bw=\pmb 0$ if and only if  the class of $\delta$ contains a positive
divisor. Both 
$\bl_0$ and 
$\theta(\bl_0,\bw+\mu_0)$ (uniquely) analytically continue along any path $P$ in $\sT_{g}$ based at $\bl_0$.   
\end{thm}

\begin{proof}[Comments on the proof] \cite[p.~7, Thm.~1.1]{Fay} states the characterization of $W_{g-1,m}-\delta$ as being some translate of
the
$\theta$ divisor. It comes from Riemann's precise solution of the Jacobi Inversion Problem. The minimal expression of that says 
the map 
\eqref{jacInv} from $X_m^g$ to $\bC^g/\lrang{\Pi(\bl)}$ is onto. The characterization of the $\theta$ divisor is that these are the
(degree 0) divisor classes of form $[D-\delta]$ where this map fails to be one-one, having as fiber a copy of projective space of dimension
one less than that of the linear system of $D$. He quotes \cite{lewittes} or \cite{mayer} for a proof. 

From Riemann-Roch:  $W_{g-1,m}-\delta$ is closed under multiplication by $-1$; and it determines the $\theta$ function with it as divisor up
to a   holomorphic exponential  (in $\bw$), which we can take to be even in $\bw$. So,  
$\theta(\bl, \bw+\mu)$, $2\mu\equiv \pmb 0 \mod \lrang{\bl}$, is either: 
\begin{equation} \label{even-odd} \begin{array}{rl} \text{even: } \theta(\bl, -\bw+\mu)&=\theta(\bl, \bw+\mu); \text{ or}\\ 
\text{odd: }\theta(\bl,
-\bw+\mu)&=-\theta(\bl,
\bw+\mu).\end{array}  \end{equation}

Suppose  a path $P: [0,1]\to \sT_g$ starts at $\bl_0$. Points on $\bC^g$ representing 2-division points on
$\bC^g/\lrang{\Pi(\bl)}$ form a discrete set. So, you can uniquely assign  
$t\in [0,1] \mapsto \mu_t
\in
\bC^g$ representing a 2-division point on $\bC^g/\lrang{\Pi(P(t))\eqdef \Pi(\bl_t)}$ to be continuous in $t$. Then, 
\begin{equation} \label{analconttheta} \theta(\bl_t,\bw+\mu_t)\text{ on }\bl_t\times \bC^g \end{equation} analytically continues 
$\theta(\bl_0,\bw+\mu_0)$ along $P$. \end{proof} 

\newcommand{\hct}{\text{\rm hc}}

\subsubsection{Even and odd $\theta\,$s and the spaces $\sM_{g,\pm}$} \label{tor2} A function 
$\theta(\bl_0,\bw+\mu_0)$ in Thm.~\ref{eventhetanulls} is called a $\theta$ with (a 2-division) characteristic. Even if $P$ is a closed
path, we don't expect $\bl$ or $\mu$ at the beginning and end  of $P$ to be the same.  
  
Now return to  the family $\sH$ in
\S\ref{hcfamilies}, $\bp_0\in \sH$ and let  $\delta_0\in \sM_{g,\pm}$ be a \hc\ class on $X_{\bp_0}$. Any $\bl_0$  in Torelli space
over
$\bp_0$ determines the unique theta $\theta(\bl_0,\bw+\mu_0)$ whose  divisor
$\Theta_{\bp_0}$ is the pullback of 
$W_{g-1,\bp_0}-\delta_0$. Let $\delta_t$ be the value at $t$ of the unique lift of the path
$P$ in $\sM_{g,\pm}$ starting at $\delta_0$. Then,  the divisor of
expression \eqref{analconttheta} is the pullback of $W_{g-1,P(t)}-\delta_t$. 

\begin{prop} \label{thetaChoices} If $\theta(\bl_0,\bw+\mu_0)$ is odd (resp.~even), then so is  $\theta(\bl_t,\bw+\mu_t)$   for all
 $t\in [0,1]$. Suppose $\sH$ is a family of covers with odd order branching, $\delta_0$ is the \hc\ class defined by Lem.~\ref{lnth.2} at
$\bp_0\in \sH$, and the divisor of $\theta(\bl_0,\bw+\mu_0)$ is the pullback of $W_{g-1,\bp_0}-\delta_0$. Then, if 
$P:[0,1]\to \sH$ is a closed path,  the divisors of  
$\theta(\bl_0,\bw+\mu_0)$ and $\theta(\bl_1,\bw+\mu_1)$ are the same. 
\end{prop} 

\begin{proof} Apply   Thm.~\ref{eventhetanulls} to the path $t\mapsto
(\Psi_{\sH,\sM_g}\circ P)(t)$ to get the analytic continuation. We know for each  $t$,  $\theta(\bl_t,\bw+\mu_t)$ is one of $\pm
\theta(\bl_t,-\bw+\mu_t)$. Their ratio is  continuous  in all variables, avoiding  $(t,\bw)$ that make the
denominator 0. So, the value is either $+1$ or
$-1$ giving the first conclusion. 
The 2nd conclusion follows from  Lem.~\ref{lnth.2}, saying $P(t)$ determines $\delta_t$, and $P(0)=P(1)$. 
\end{proof} 

Now consider in Prop.~\ref{hcspaces} the respective degrees of $\sM_{g,\pm}$ over $\sM_g$. Prop.~\ref{thetaChoices} says  analytic
continuation of \hc\ classes from one point on the connected space $\sM_g$ to another moves even (resp.~odd) classes to even (resp.~odd)
classes. So, for each
$g$, use connectedness of $\sM_g$ to analytically continue to where we can count these classes. Conclude this count using  
hyperelliptic curves in  Cor.~\ref{hypcomps}.

Next, consider for $\sH=\sM_g$ why the monodromy action is transitive on even (resp.~odd) \hc\
classes. That means, for paths  $P$ running over lifts to $\sT_g$ of closed paths in $\sH$, the
action is  transitive on both even and odd
$\theta\,$s.   

Assume along  $P$ based at   
$\bl_0\in \sT_{g,m}$,  the endpoint is 
$(\bl_0)\Psi_P\eqdef\bl_1\in \sT_{g,m}$, and
$\mu_0$ continues to
$\mu_1$. A  formula  explains how $\theta$ functions transform with an 
application 
of $\psi_P$ to  $\Pi(\bl_0)=(2\pi i I_g, \tau_{\bl_0})$ (\cite[p.~7]{Fay}, from \cite[p.~84]{igusa}). Write  $$\tau_{\bl_1}\eqdef 2\pi
i(A_P\tau_{\bl_0}+ 2\pi  iB_P)( C_P\tau_{\bl_0}+ 2\pi iD_P)^{-1}.$$ Denote  $C_P\tau_{\bl}+ 2\pi iD_P$ -- its entries are
functions in the entries of $\tau_{\bl}$ -- by
$M_\bl$.  If $M^*$ is any matrix with entries that are functions of the entries of $\tau_\bl$, then   
$\nabla_\Pi(M^*)$ is the matrix whose
$(i,j)$ entry is the partial with respect to the variable in the $(i,j)$ 
position of $\tau_\bl$. Then,
$\Pi(\bl_1)=\smatrix D C B A \Pi(\bl_0)$ and
$2\pi i\bw =\tilde
\bw M$. \cite[p.~8]{Fay} notation is compatible with  \S\ref{famCovers} by writing 
$\Pi(\bl)$ as $\bigl(\frac{\balpha}{\bbeta}\bigr)$ and having $\Psi(P)$ from the left as 
$U_P=\smatrix {A_P} {B_P} {C_P} {D_P}$ in
\eqref{SP2g}. The result is $\bigl(\frac{\tilde \balpha}{\tilde \bbeta}\bigr)$, with the following two provisos. First: Normalize by
multiplying on the right by the inverse of the matrix $\tilde \balpha$, and then multiply by $2\pi i$. Second: 
Fay  changes  $U_P$ to $\smatrix {D_P}  {C_P} {B_P} {A_P} $ (this is also in
$\Sp_{2g}(\bZ)$): He thought the result in more standard notation?      

For a $g\times g$ matrix $N$, use bracket notation $\{N\}$ for 
the vector of diagonal elements. The transformation formula:  
\begin{equation} \label{transform} \begin{array}{rl} \text{with}
& \mu_1^\tr= \smatrix D {-C} {-B} A
\mu^\tr +\frac 12 \Bigl[\begin{array}{c}\raise-2pt\hbox{{$\{C 
D^\tr\}$}}\\ \raise2pt\hbox{{$\{A B^\tr\}$}}
\end{array}\Bigr], \text{ then }\theta(\bl_1,\tilde \bw +\mu_1) = \\ \qquad  
&K_{\mu_0,P}\det(M)^{\frac 1 2}_{|\tau_{\bl_0}}  \exp(\frac 1 2 (\bw \nabla_\Pi \ln 
(\det(M))_{|\tau_{\bl_0}}  \bw^\tr)\theta(\bl_0,\bw+\mu_0). \\\end{array}\end{equation} 

Transitivity on even and odd $\theta\,$s is very old;  a corollary of \eqref{transform} by applying elements of 
$\Sp_{2g}(\bZ)$. That all even
$\theta\,$s are nonzero at the origin on a general surface is explained at 
\cite[p.~7]{Fay} by alluding to \cite[Cor.~3.2]{Fay}. Hershel Farkas  -- in the late
60's when we were Stony Brook colleagues -- attributed this argument to
\cite{poincare}.  Fay  expands along these lines: Deforming from period matrices of
\lq\lq products of elliptic curves\rq\rq\ to discern  objects on a general Riemann surface.  

Denote the function of $(t,\bw)$ in  Prop.~\ref{thetaChoices} corresponding to a path $P$ by $\theta_P(\bw)$.  As in \S\ref{hcfamilies}, for
  $\Psi:\sT\to \sH$ a family of surfaces, form $\sH_{\pm}$. 
\begin{defn}\label{nontrivthetanull}  We say $\Psi$ supports an even (resp.~odd)  \hc\ class 
if a component of $\sH_{+}$ (resp.~$\sH_-$) maps one-one from the fiber product to $\sH$. We say  an even class has a nontrivial
$\theta$-null if there is a path $P$ on this component (based at $(\bp_0,\delta_0)$ on the component), and some corresponding $\bl_0$, so
that for some $P$,  
$\theta_{P}(\bw)$ in Prop.~\ref{thetaChoices},  $\theta_{P}(\pmb 0)\not \equiv 0$ (as a function of $t$).  
\end{defn} 

 If $\Psi$ supports a \hc\ class, but the component is in $\sH_-$,
then
$\theta_{P}(\bw)$ is  odd in $\bw$ (Thm.~\ref{thetaNullRes}). Conclude: $\theta_{P}(\pmb 0)$ is an identically 0 function
of $t$. Still, we can ask if such a component defines a {\sl nondegenerate} odd $\theta$. That is, there is a path $P$ so that for most
values of $t$, the gradient of
$\theta_{P}(\bw)$ in $\bw$ doesn't vanish at $\bw=\pmb 0$.   

Thm.~\ref{thetaNullRes} shows  the $\oplus$ families of Figure 1 (Thm.~\ref{thmB}) have a nontrivial $\theta$-null for many
values of 
$(n,g)$. For $g=1$ (and $\sH=\sM_1$) there is a unique  nondegenerate odd $\theta$, but for $g\ge 2$, there is more than one (as in
Prop.~\ref{specreps}). See
\S\ref{comphc}. 

\begin{rem}[Stacks of compact Riemann surfaces] \label{stack} Significantly, $\sM_g$ has no total family $\sX\to \sM_g$
representing its points.  {\sl Stacks\/} arise in such situations to produce the {\sl algebraic\/} 
map
$\Psi_{\sH,\sM_g}:
\sH\to\sM_g$ (before Prop.~\ref{thetaChoices}). \cite[\S4]{Fr1} shows, for an ordinary (say, absolute or inner) Hurwitz space
(\S\ref{hurSpaces})  
$\sH=\sH(G,\bfC)$, the  space
has (a finite number of explicit) Zariski open subsets $\{U_i\}_{i=1}^s$, each having an \'etale
cover
$\Phi_i: W_i\to U_i$ with this property. There is a total family
$\sT_i\to W_i\times \prP^1_z$ so that for $w_i\in W_i$, the fiber $\sT_{i,w_i}\to w_i\times \prP^1_z$ over $w_i$
represents the equivalence class of covers of the image of $w_i$ in $\sH$. That is, the stack exists in the \'etale
topology. This produces  
a {\sl stacky\/} definition $\Psi_{\sH,\sM_g}$. There is always a unique global total family if the Hurwitz space has fine
moduli: as in \S\ref{intSNorm}, self-normalizing  in the absolute case; $G$ has no center in the
inner case. Example: In Lem.~\ref{ratPtSpin}, $\Spin_n$ has a nontrivial center, so  $\bp\in \sH_+(\Spin_n,\bfC)^\inn(K)$ may not
guarantee a
$K$ cover. 

\cite{wew} uses stack language; behind it is our construction. Pull reduced Hurwitz spaces
(\S\ref{rnz.4}) back to ordinary Hurwitz spaces; define the map from this. \end{rem}  
 
\begin{rem}[Warning!] Don't confuse 2-division points with \hc\ classes. Any \hc\ class  translates 2-division
points into \hc\ classes: The former is a homogeneous space for the latter. Still, the monodromy action on 2-division
points over $\sM_g$ has a length 1 orbit   (the origin) and another of length $2^{2g}-1$; different from the two orbit lengths on \hc\
classes, in  Prop.~\ref{hcspaces}.
\end{rem}

\subsection{$\theta$-nulls and roots of 1} \label{automorphy} \cite[Thm.~27.7]{Sh98} (in greater generality, quoting
\cite{KP1892}, but \cite[\S 28]{Sh98} has a detailed proof) says 
$K_{\mu_0,P}$ in  \eqref{transform} is a root of 1 that depends on which branch of 
$\det(M)^{\frac 1 2}$ we have chosen. 

\subsubsection{Hurwitz-Torelli automorphic functions and roots of 1} \label{HTautomorphic}  Consider a braid orbit $O$ on a reduced Nielsen
class
$\ni(G,\bfC)^{*,\rd}$, $*=\abs$ or $\inn$, and denote its corresponding reduced Hurwitz space component by $\sH_O$. Regard $\bp_0\in \sH_O$
as a base point for analytic continuation. Suppose $f$ is a meromorphic function around  $\bp_0$ and it continues  
along {\sl any\/} path on $\sH_O$ based at $\bp_0$. Denote the subgroup,  of
$H_r$ that fixes a particular element, $\bg_O\in O$, by $H_{r,\bg_O}$.  Then, the monodromy action on  $f$ explicitly  interprets as an
action of $H_{r,\bg_O}$: $q\in H_{r,\bg_O}$ takes $f$ to $f_q$.  

Choose a classical generators
(\S\ref{gpNC}) of
$\pi_1(U_\bz,z_0)$ to determine how $H_r$ acts in \S\ref{suppLems}. If in this identification 
$\bg_0$ is a branch cycle description of  $\phi_0: X_0\to \prP^1_z$, then by restriction, $H_{r,\bg_O}$ 
acts on the fundamental group of
$X_0\setminus\{\phi^{-1}(\bz)\}$. 

\begin{defn} \label{HT} Refer to $f$ as a $\pi_1$-H(urwitz)-T(orelli)  (resp.~$H_1$-H-T) function if the  $H_{r,\bg_O}$ action determined
by it factors through
$\pi_1(X_{\bp_0})$ (resp.~$H_1(X_{\bp_0,\bZ})$). For an $H_1$-H-T function, in the notation of \eqref{transform}, $q\in H_{r,\bg_O}$ acts
through $\smatrix {D_q}{C_q}{B_q}{A_q}$ with an associated $g\times g$ matrix $M_q(\bl)=C_q\tau_{\bl}+2\pi iD_q$. With $P$ a path based at
$\bp_0$, for  
$\bp=P(t)$ denote the analytic continuation of $\bl_0$ over $P(t)$ by  $\bl=\tilde P(t)$. An
$H_1$-H-T function $f$ is {\sl (H-T) automorphic}, of weight $m$,  if: $$ \text{for each $q\in H_{r,\bg_O}$, 
and path $P$ based at  $\bp_0$},  f_q(P(t))=M_q(\tilde P(t))^mf(P(t)).$$  \end{defn}  

The automorphic definition matches the form of \cite[\S25.1]{Sh98}. We now give examples of
$H_1$-H-T functions. Start with $\bp_0\in \sH(G,\bfC)^{*,\rd}$, and
let
$\delta_{\bp_0}$ be any
\hc\ class on $X_{\bp_0}$ defining the $\Theta_{\bp_0}$ as in Thm.~\ref{eventhetanulls}.

Assume
$\bl_0$ and $\mu_0$ in Prop.~\ref{eventhetanulls}, over $\bp_0$,  give an even $\theta$ with zero divisor $\Theta_{\bp_0}$. For $P:[0,1]\to
\sH(G,\bfC)^{*,\rd}$ a path based at $\bp_0$, denote  analytic continuation of the
$\theta$ (resp.~$\delta_{\bp_0}$) along $P$ by 
$\theta_{P}(\bl_0,\bw+\mu_0)$ (resp.~$\delta_P$). It's value at $t\in [0,1]$ is then 
$\theta(\bl_t,\bw+\mu_t)$ (resp.~$\delta_{P(t)}$) compatible with the notation of \S\ref{tor2}.  Denote the
connected component of
$\sH_O\times_{\sM_g}\sM_{g,+}$ containing
$(\bp_0,\delta_{\bp_0})$ by
$\sH_{O,\delta_{\bp_0}}$. 

\begin{defn} \label{thetanull} The result, $\theta_{P}(\bl_0,\mu_0)$, of setting   $\bw=\pmb 0$ and letting $P$ vary over all paths $P$
based at
$\bp_0$ is a $\theta$-{\sl null} on $\sH(G,\bfC)^{*,\rd}$.\end{defn} 
\smallskip
Let $P'$ be a closed path on $\sH(G,\bfC)^{*,\rd}$ representing $q\in H_{r,\bg_O}$ with $\bl_1$ and $\mu_1$ the 
analytic continuations of $\bl_0$ and $\mu_0$ to the end of $P'$. With
$P$ running over all paths based at $\bp_0$, 
\eqref{transform} compares  
$\theta_{P'\cdot P}(\bl_0,\mu_0)=\theta_{P}(\bl_1,\mu_1)$ and 
$\theta_{P}(\bl_0,\mu_0)$. Analytic continuations only depend on the homotopy class of a path with the homotopy keeping the  endpoints
fixed. So, we can replace the fixed closed path  $P'$ by $q$. 
\smallskip

\begin{prop} \label{thetanullprop} \!\!The $theta$-null $\theta_{P}(\bl_0,\mu_0)$ \!is an \!$H_1$-H-T function. \!It is  \!iden\-ti\-cally
zero if and only if $\delta_{\bp}$ contains a positive divisor for each 
$(\bp,\delta_{\bp})\in \sH_{O,\delta_{\bp_0}}$. 

For $\ni(G,\bfC)$ a Nielsen class of odd-order branching, and  $q=[P']$ as above, 
$$\theta_{q\cdot P}(\bl_0,\mu_0)=\theta_{P}(\bl_1,\mu_1)=K_{\mu_0,q}\det(M_q(\tilde P))^{\frac 1 2} \theta_{P}(\bl_0,\mu_0).$$ For some
minimal positive integer $m$,   $(\theta_{P}(\bl_0,\mu_0))^m$ is automorphic under $H_{r,\bg_O}$.  
\end{prop} 

\begin{proof} Most everything follows from the definition of the $\theta$-null and \eqref{transform}, except these two. The criterion for
nonzeroness at the origin of a $\theta$ is from Thm.~\ref{eventhetanulls}, and we apply Prop.~\ref{thetaChoices} to a Nielsen class of
odd order branching. The minimal integer
$m$ in the last paragraph is one for which
$K_{\mu_0,q}^m=1$ for $q$ running over any finite set of generators of $H_{r,\bg_0}$.  See Rem.~\ref{compExp}. 
\end{proof} 
 
When 
$g=1$, and $H_{r,\bg_O}$ identifies with $\PSL_2(\bZ)$, \cite[p.~101]{FKraTheta} shows explicitly $K_{\mu_0,q}\in \lrang{\sqrt{-i}}$: an 8th
root  of 1, where the  serious calculation is $q$ representing $\smatrix 1101\in \SL_2(\bZ)$. \cite[p. 176]{FKraTheta} repeats the
classical definition of automorphic form in
\lq\lq dimension 1,\rq\rq\ though as with everything in this book the functions are supported on congruence subgroups. In
particular, they explicitly compute that root of 1 in \cite[Prop. 2.1]{FKraTheta} for a $\theta$-null with other rational characteristics --
notating these as $\chi$ -- where the denominator in them is  an integer $k$. Such functions are supported on the
upper-half plane quotient of the congruence subgroup 
$\Gamma(k)$. The precise result is \cite[Prop. 2.1]{FKraTheta} where the $K_{\chi,q}$ are $8k$-th roots of
1. 

\begin{rem}[Effective monodromy] \label{compExp} Since $H_r$ is a
presented group, and $H_{r,\bg_0}$ is a stabilizer subgroup, an algorithm for computing generators of
$H_{r,\bg_O}$ comes from the Schreier algorithm for generators of a subgroup of finite index in a finitely generated free group 
\cite[p.~351]{FrJ}. Suppose $X_{\bp_0}\to\prP^1_z$ has branch cycles $\bg_0$, and $\hat X_{\bp_0}\to \prP^1_z$ is its Galois
closure. Let $X^{\text{un}}_{\bp_0}$ be the maximal unramified cover of $X_{\bp_0}$ lying between $\hat X_{\bp_0}$ and $X_{\bp_0}$.
\cite[Ex.~3.4]{combcomp} gives examples where $X^{\text{un}}_{\bp_0}\ne X_{\bp_0}$. Yet, typically they are equal, as for 
the Nielsen classes $\ni(A_n,\bfC_{3^r})^\abs$.  Given explicit
$q\in H_{r,\bg_O}$, \cite[\S3.2]{combcomp} displays   $\pi_1(X^{\text{un}}_{\bp_0})$ (and so $H_1$) {\sl in terms
of the branch cycles $\bg_O$}. That is, the presentation supports an explicit action of
$q$.   
\end{rem}

\begin{rem}[Modular Towers] \label{ModTowers} Each space $\sH_+(A_n,\bfC_{3^r})^\abs$ in Thm.~\ref{thmB} has above it 
an infinite collection of Hurwitz spaces for which Prop.~\ref{thetanullprop} produces an even $\theta$-null. The case $n=3$ is
included, where  for each prime $p$ (excluding
$p=3$) and for each nonnegative prime power $p^\ell$, the system has a Hurwitz space attached to a centerless group $(\bZ/p^{\ell})^2\xs
\bZ/3$, and the four conjugacy classes $\bfC_{\pm 3^2}$ (Ex.~\ref{pm3^2}). \cite[\S6]{lum} discusses the modular curve-like properties of
this system. Similarly, starting with $n=5$ and $r=4$ in Thm.~\ref{thmA}, though less obvious what are the ingredients for a system.   We can
compare the case of Prop.~\ref{A43-2} with 
\cite{FKraTheta}: Like a modular curve it is a family of genus 1 covers. Indeed, all these spaces are  1-dimensional
quotients of the upper half-plane, but they aren't modular curves. \cite[Main Thm.]{twoorbit} puts these families in a bigger context, and
extends the modular curve-like properties. As in the
 Farkas-Kra discussion above, if you go \lq\lq up\rq\rq\ in these systems, higher levels support higher characteristic 
$\theta$-nulls.   
 \end{rem} 

\subsubsection{$\theta$-nulls on  $\sH_{\pm}(A_n,\Ct 
r)^{*,\rd}$} \label{oddord} Consider $\sH$, a reduced Hurwitz space component of odd order branched  covers. Prop.~\ref{thetaChoices} canonically gives on it an analytic continuation of an even (resp.~odd)  $\theta$, if  the  \hc\ classes on the component are even (resp.~odd; Def.~\ref{nontrivthetanull}).  \cite{Se3} determines which from the Nielsen class and the component lifting invariant. Prop.~\ref{thetanullprop} gives the transformation formula for the corresponding $\theta$-null. It can only be nonzero if the $\theta$ is even. Assume $n\ge 4$. 

\begin{thm} \label{thetaNullRes} Assume $g=r-n+1$, the genus of the degree $n$ covers parametrized by $\sH_{\pm}(A_n,\bfC_{3^r})^{\abs,\rd}$ (Thm.~\ref{thmB}) is at least 1. For $r$ even, the  $\theta$ is even (resp.~odd) when $\sH=\sH_+(A_n,\bfC_{3^r})^{\abs,\rd}$ (resp.~$\sH_-(A_n,\bfC_{3^r})^{\abs,\rd}$). For $r$ odd, the results are switched. For the inner case, the result is independent of the parity of $r$: the  $\theta$ is even (resp.~odd) when $\sH=\sH_+(A_n,\bfC_{3^r})^{\inn,\rd}$ (resp.~$\sH_-(A_n,\bfC_{3^r})^{\inn,\rd}$). 

For $g=1$ or  for $n\ge 12\cdot g+4$,  the natural map $\Psi_{\sH_{\pm}, \sM_g}: \sH_{\pm}(A_n,\bfC_{3^r})^{\abs,\rd} \to \sM_g$ restricted 
to each component is dominant.  If  also, $r$ is even (resp.~odd), then
$\sH_+(A_n,\bfC_{3^r})^{\abs,\rd}$ (resp.~$\sH_-(A_n,\bfC_{3^r})^{\abs,\rd}$) supports a nonzero $\theta$-null.   \end{thm} 

\begin{proof} In \cite[Thm.~2]{Se3} we take the special case $X=\prP^1_z$, so  \cite[exp. (17)]{Se3} applies. Serre notes this reverts to \cite[exp. (10)]{Se2} (which references an early version of Cor.~\ref{cnf.3}): For $\bg=(\row g r)$ in the braid orbit of a Nielsen class $\ni(G,\bfC)$, we get even for the \hc\ class exactly when the product of the lift invariant \eqref{spinlift} and   $(-1)^{\sum_{i=1}^r w(g_i)}$  is 1: Serre's formula written multiplicatively. This only depends on the Nielsen class and the lifting invariant. In the 3-cycle absolute case each $w(g_i)$ is 1. For, however, the inner case, each $w(g_i)$ is $n!/6$, which is even (for $n\ge 4$). 

We review ingredients from 
\cite[Thm.~1]{AP}. Let $X$ be a compact
Riemann surface of genus $g$. Consider  $n
\ge  12g+4$ an integer, $\pmb k$ a \hc\ class on $X$ and $x_1,x_2,x_3$  any distinct points on $X$. Assume also:  
\begin{edesc} \label{exists} \item \label{existsa} There exists a meromorphic \hc\
differential (expression \eqref{half-values}) 
$\mu$ whose divisor of poles is $\le D_{X,{\pmb n}}\eqdef n_1x_1+n_2x_2+n_3x_3$; and 
\item \label{existsb} the square $\mu\otimes \mu$ of $\mu$ is the differential $df$ of a meromorphic $f$ on $X$. \end{edesc}
From the Riemann-Roch Theorem (see \S\ref{comphc}), \eql{exists}{existsa} follows if  a \hc\  class has sufficient polar degree to
guarantee a global section. Denote the linear system of sections of $\pmb k$ with polar divisor $\le$  a divisor $D$ by $H^0(X,O(\pmb
k,D))$.

The argument for \eql{exists}{existsb} crucially gives a differential satisfying  Square Hypothesis
\eqref{merSquare}. Consider the $\bC$ bilinear pairing  $$\Gamma: H^0(X,O(\pmb k,D))\times H^0(X,O(\pmb k,D)) \to H^1(X\setminus D, \bC)$$ by
 $(s_1,s_2)\mapsto s_1\otimes s_2$, a differential with pole divisor supported in $D$. Let $\Delta_{\pmb k}$ be restriction of $\Gamma$ to the
diagonal. Then,  $\Delta_{\pmb k}^{-1}(0)$ is  sections whose square is exact: 

\begin{triv} \label{conclus}  $df$ whose divisor satisfies the square hypothesis \eqref{merSquare}, so  \\ $f: X\to
\prP^1_z$ has odd ramification. \end{triv} 
\noindent If  $\dim(H^0(X,O(\pmb
k,D)))>  \dim(H^1(X\setminus D,
\bC))$, by putting   
$\Gamma$ in standard form see that  $\Delta_{\pmb k}^{-1}(0)\setminus \{0\}$ is nonempty.  We can change $\pmb n$ to have this happen because
$D$ has just three support points, so the target dimension doesn't change  with
$\pmb n$. We want to know in increasing order for which
$(n,r)$ the following hold: 

\begin{edesc} \label{goal} \item  \label{goala} Whether (and which) $\Psi_{\sH_{+}, \sM_g}$ and/or $\Psi_{\sH_{-}, \sM_g}$ is dominant; and 
\item  \label{goalb} if some cover $\phi_\bp: X_\bp\to \prP^1_z$ in such a component produces a $\theta$ whose value at $\pmb 0$
is nonzero. 
\end{edesc}  All even $\theta\,$s on a  general curve of genus $g$ are nonzero at $\pmb 0$ (above Def.~\ref{nontrivthetanull}). 
Suppose $\sH$ is a component which supports an even $\theta$-null. 
Conclude:  
\eql{goal}{goala} holding for  $\Psi_{\sH, \sM_g}$ implies \eql{goal}{goalb}. When $g=1$, \eql{goal}{goala} holds \cite{FrKK}.   

Since, however, the divisor of poles has just three points of support, it doesn't tell us which $X\,$s are in the image of points in
$\sH_{\pm}(A_n,\bfC_{3^r})$. It is easy that a dense set of  $m\in \sM_g$ represent $X_m$ with an odd order branched cover 
not in any Nielsen classes $\ni(A_n,\bfC_{3^r})$. Yet,   $X_m$ may {\sl also\/} be a cover in such a Nielsen class. 

It
is also easy to show that if $m\in \sM_g$ has complex coordinates that over $\bQ$ generate the function field of $\sM_g$ ($m$ a generic
point), then any odd order branched cover $f_m: X_m\to\prP^1_z$ (degree $n$) is primitive \cite[p.~26]{Fr1}. 

It is much harder that if $g\ge 3$ (and $m$ generic) then the
monodromy group of $f_m$ must  be $A_u$ for $u$ given in Rem.~\ref{genAnmon}. To show \eql{goal}{goala} requires 
knowing when we can find an alternate $f^*_m$ to  $f_m$, so  $f_m^*: X_m\to\prP^1_z$ has 3-cycle branching. 

The answer is
always. Start by writing each of the branch cycles $\bg=(\row g r)$ for $f_m$ as a product of 3-cycles, to give a
(possibly) larger value $r^*$ and branch cycles $\bg^*\in \ni(A_n,\bfC_{3^{r^*}})$ with $r^*-n+1=g$. Then, apply
Riemann's existence theorem to produce a cover  that must also represent a generic surface of genus $g$, and so by
specialization represents
$X_m$ as a 3-cycle cover of $\prP^1_z$ (\cite[\S4]{Fr1} or \cite[\S5]{AP}). 

Notice, however, this argument \lq\lq deforms\rq\rq\ between different Nielsen
classes. To be certain both components of $\sH_{\pm}$ (for the appropriate $(n,r)$) contain a \lq\lq generic\rq\rq\ Riemann surface, we need
to know this deformation preserves information about the evenness and oddness of the \hc\ classes. 

Both \cite[following
Thm.~2]{Se3} and
\cite[\S5]{AP} have arguments that  handle this: The former uses  topology to characterize \hc\  parity. The latter would adjust to use that (in our 3-cycle case) the even \hc\ classes correspond to an unramified spin cover
of the Galois closure of $\phi: X\to \prP^1_z$ (Cor.~\ref{cnf.3}).      
\end{proof}

\begin{quest} \label{charthetanull} 
For those absolute and inner Hurwitz space components carrying an even $\theta$-null in Thm.~\ref{thetaNullRes} (or as in Rem.~\ref{ModTowers}), generalize the theorem to find those for which the $\theta$-null is nonzero. 
\end{quest}

\begin{rem}[Generic alternating group monodromy] \label{genAnmon} Let $f_m: X_m\to\prP^1_z$ present the generic compact Riemann surface of
genus $g$ as an odd branched cover following the proof of Thm.~\ref{thetaNullRes}. If $g\ge 3$, then the monodromy group is a copy of
$A_{\deg(f_m)}$ according to   \cite{GN95}, \cite{GM98} and  \cite{GS07}. The case $g=0$ is a source of considerable
applications. Then,  excluding a finite number of significant special cases,  either
the monodromy group is $A_{\deg(f_m)}$, or it is $A_l$ with 
$\deg(f_m)=\frac{l(l-1)}2$ ($A_{l}$ acting on unordered pairs of integers from $\{1,\dots,l\}$). \cite[p.~76]{thomp-gen0} lists all cases,
not just odd order branching. \end{rem}

\subsection{The Hyperelliptic Locus} \label{hyplocus} Suppose the affine part of a
hyperelliptic curve $X$  is $\{(z,w)\mid
w^2=h(z)\}$. \S\ref{sqhypth} lists differentials
$\omega$ satisfying  square hypothesis \eqref{merSquare}. \S\ref{sqdivisors} uses these to list
representative divisors for all \hc\ classes, and Prop.~\ref{specreps}  computes the monodromy orbit
lengths.

\subsubsection{Half-canonical divisors} \label{halfcan} Suppose $\omega$ is a meromorphic 
differential on a Riemann surface $X$, written locally as $f_\alpha(z_\alpha)dz_\alpha$, as in \S\ref{cocycs} (or in detail in
\cite[Chap.~2, \S2.4]{FrBook}), on simply connected   domains $U_\alpha$. On $U_\alpha$ its divisor is the divisor
$(f_\alpha(z_\alpha))$ of the function. Assume also the {\sl square hypothesis}: 
\begin{triv} \label{merSquare} 
$(f_\alpha(z_\alpha))$ has the form
$2D_\alpha$ for $U_\alpha$ running over a subchart covering $X$.\end{triv}
\noindent Then, there is a branch $h_\alpha(z_\alpha)$ of square
root (of $f_\alpha(z_\alpha)$) on $U_\alpha$ \cite[Chap.~2, \S6.1]{FrBook}. Of course, there are two of
these; our notation means we have chosen  one. Call the symbol collection $\{\tau_\alpha=
h_\alpha(z_\alpha)\sqrt{dz_\alpha}\}_{\alpha\in I}$, a {\sl half-canonical 1-chain\/} on $U_\alpha$. In bundle
language, this is a 1-chain with values in the square-root of a canonical bundle. Still, it is more than
that, for  the squares of these form a global differential on
$X$. So, we refer to 
$\{h_\alpha(z_\alpha)\}_{\alpha\in I}$ by $\bh$ and call it a square-root of $\omega$. 

\begin{lem}[Half-canonical divisor]  The 1-chain 
$\{(h_\alpha(z_\alpha))\}_{\alpha\in I}$ from a square root of $\omega$ 
give a well-defined divisor: a half-canonical divisor on
$X$. 
\end{lem}

\begin{proof} Let $D=(\omega)$ be the divisor of $\omega$. Since, $h_{\alpha}^2=f_{\alpha}$, the support
multiplicities of $D$ are all even integers. So, a square-root of $\omega$ defines $D_{1/2}=(\omega)/2$, a
divisor uniquely given  by the zeros and poles of the $h_{\alpha}\,$s. 
\end{proof}

\subsubsection{Square-hypothesis for hyperelliptic curves} \label{sqhypth}  App.~\ref{cocycs} tells precisely what it means -- expression
\eqref{half-values} -- for there to be a
\hc\ differential representing a \hc\ divisor class. This subsection presents the much easier case of hyperelliptic curves.  With no
loss, assume an  
odd degree polynomial $h$ with distinct zeros $\row z {r-1}$ (use $z_r=\infty$). Denote the point on $X\eqdef X_\bz$ over $z_i$ by
$x_i$, with
$x_\infty$ lying over $z=\infty$. As in \cite[p.~7]{Mu76}, form  
$$\omega_i=\frac{(z-z_i)^{\frac 12}}{(\prod_{j\ne i}z-z_j)^{\frac 1 2}}\,dz,\ i=1,\dots,r-1.$$ 
Since $w=\sqrt{h(z)}$, the factor in front of the $dz$ (a meromorphic differential on $X$) in $\omega_i$ is just
$\frac{z-z_i}{w}$, a meromorphic function on $X$. 

Express $\omega_i$ in a local parameter $t_{z'}$ -- a fractional power of
$(z-z')$ -- over each \\ $z'\in \prP^1_z$. Take $t_{z'}$ to be $z-z'$ for $z'\in U_\bz=\prP^1_z\setminus \{\bz\}$,
$(z-z')^{\frac 1 2}$ for $z'\in \{z_1,\dots,z_{r-1}\}$ and 
$1/z^{\frac 1 2}$ for $z'=\infty$. The multiplicity  of
$\omega_i$  at $t_{z'}=t=0$ is 0 for $z'\in U_\bz$; the multiplicity, 2, of $t^2/h(t)\,d(t^2)$  (including 
 $t$ in the denominator) for $z'=z_i$; 0 for $z'=z_j$, with $j\ne i$; and the multiplicity, $r-6$, of
$((1/t^2-z_i)/h(1/t^2))^{\frac 1 2}\,d(1/t^2)$ at $t=0$ at $z'=\infty$. Conclude: The total divisor of
$\omega_i$ is \begin{equation} \label{defDi} 2x_i+(r-6)x_\infty=2\cdot D_i;  x_i+(\frac{r-6}2)x_\infty\text{ is a \hc\
divisor}.\end{equation} The following divisors are linearly equivalent to 0:
\begin{edesc} \label{extraEquiv} \label{linequiv0} \item \label{linequiv0a} $2(x_i-x_j)$ (divisor of  $\frac{z-z_i}{z-z_j}$); and 
\item \label{linequiv0b} $\sum_{i=1}^{r-1}
(x_i-x_\infty)$ (divisor of
$\sqrt{h(z)}$). \end{edesc}

\subsubsection{$\sM_{g,\pm}$ over $\hyp_g$} \label{sqdivisors}  Use 
$\bz$ as the basepoint of
$U_r$. Start with
$\omega_1$ and multiply it by powers of the $(z\nm z_i)\,$s to get differentials satisfying square
hypothesis \eqref{merSquare}. Divide by 2 to get representatives of all \hc\ classes on $X_\bz$ from $$\sD\eqdef
\{D_1+\sum_{1\le i\le r-1}\epsilon_{i}(x_i -x_\infty)\}_{\epsilon_{i}\in \{0,1,-1\}}.$$ The expression {\sl orbit\/} refers to the
action of
$\pi_1(U_r,\bz)$ (or of $\pi_1(U_r/\PGL_2(\bC),\bz_0)$; \S\ref{gpNC}) on a divisor class. For convenience, we sometimes write
$x_\infty=x_r$. 

The divisor $mx_i$, $m$ odd (resp.~\!even), is equivalent to $x_i+(m\nm1)x_\infty$  (resp.~\!$mx_\infty$). The  
collection  $\sD$ 
modulo linear equivalence, represent all \hc\ divisors on $X_\bz$. 
From \eql{linequiv0}{linequiv0a}, you can  replace $x_i-x_j$ by $x_j-x_i$: each   $D\in \sD$ is equivalent to 
\begin{equation} \label{sstratify} \sD_s\eqdef \{\sum_{u=1}^sx_{i_u}+mx_r\}_{1\le i_1<\cdots<i_s<r}\text{ with } s+m=g\nm1.\end{equation} 
Add 
\eql{linequiv0}{linequiv0b} to any $D$ in \eqref{sstratify}to conclude it is equivalent to one in
$\cup_{s=0}^{\frac{r\nm 2}2}
\sD_{s}$. 

\begin{prop} \label{specreps} The set $\cup_{s=0}^{\frac{r\nm 2}2}
\sD_{s}$   consists of inequivalent divisors. 
With $r=2\np 4\ell$ (resp.~$4\ell$), here are the $\pi_1(U_r,\bz)$ orbits. For  $s=0,2,\dots,g\nm2$
(resp. \\ $s=1,3,\dots,g\nm2$),  $\sD_{s}\cup \sD_{s+1}$ forms 1 orbit of of length 
${{r\nm 1}\choose s}+{{r\nm 1}\choose {s\np 1}}$,  leaving one (resp.~two) orbit(s)  $\sD_{g}$ (resp.~$\sD_0$ and $\sD_{g}$)
of length(s)
${{r\nm 1}\choose {\frac{r\nm 2}2}}$ (resp.~${{r\nm 1}\choose {\frac{r\nm 2}2}}$ and ${{r\nm
1}\choose 0}$). There are $1\np\ell$ orbits on the 
$\sum_{s=0}^{\frac{r\nm2}2}{{r\nm1}\choose s}=2^{2g}$ distinct classes in $\sD$. 
\end{prop} 

\begin{proof}  Suppose  $\cup_{s=0}^{\frac{r\nm 2}2}
\sD_{s}$ contains equivalent divisors.  Then, \eql{linequiv0}{linequiv0a} implies, for some $m\le r\nm 1$, $mx_\infty$ is
equivalent to $\sum_{u=1}^{s'}x_{i_u}$, contrary to $h$ is the lowest degree polynomial defining  $X$ as a hyperelliptic
curve.  

We start with the cases
$r=6$ ($g=2$) and
$r=8$ ($g=3$). For
$r=6$: 
\begin{edesc} \item $\sD_0\cup \sD_1=\{x_{i_1}, 1\le i\le 6\})$:  $6={6\choose 1}={5\choose 0}+{5\choose 1}$
elements,  and
\item $\sD_{2}=\{x_{i_1}+x_{i_2}-x_\infty,\ 2\le i<j\le 5\}$: $10={5\choose2}$ elements. 
\end{edesc} 

For all $r$, it is easy to compute the monodromy action. Each element of $\pi_1(U_r,\bz_0)$  is represented by a permutation of  entries
of $\bz$ (\S\ref{BHMon}),  inducing the same action on $\row x {r}$. Then, $\sD_0\cup \sD_1$ (resp.~$\sD_{2}$), an orbit of length 6
(resp.~10), consists of the odd (resp.~even) \hc\ classes. 

When $s=0$ and $m$ is even (say, when $r=8$), then $\sD_0$ contains $mx_\infty$ whose orbit has length 1. Here are the rest of the orbits
for 
$r=8$:   
\begin{equation} \label{r=8} \begin{array}{rl}
\sD_1\cup \sD_2:&\{ x_i+x_\infty, 1\le i\le 7\}\cup  \{x_i+x_j, 1\le i<j\le 7\}.\\
\sD_3:& \{x_i+x_j+x_k-x_\infty, 1\le i <j<k\le 7\}.
\end{array} \end{equation} 
This case has three orbits, with respective representatives (and orbit lengths) $2x_\infty\ (1), x_1+x_2\ (28={ 8\choose
2}={7\choose 1}+{7\choose 2}), x_1+x_2+x_3-x_\infty\ (35={7\choose 3})$. 

From this format, an easy induction counts the orbits as in the last paragraph of the Proposition's statement. The quantity
$\sum_{s=0}^{\frac{r\nm2}2}{{r\nm1}\choose s}$ is 
\begin{equation} \label{sumbinom} \sum_{s=0}^{r\nm 1}{{r\nm1}\choose s}/2=(1+1)^{r-1}/2=2^{2g},\end{equation} giving the count of the classes
represented in
$\sD$. 
\end{proof} 

In Prop.~\ref{specreps} notation, refer to $\sD_s\cup\sD_{s\np1}$, as appropriate, by $\sD_{s}'$. A  
\hc\ class is, respectively, even or odd  if the  dimension of its linear system (\S\ref{halfdiffs}) is even or odd.
{\sl Nondegenerate\/} if this dimension is, respectively, 0 or 1.  

\begin{cor}[Even and odd $\theta\,$s] \label{hypcomps} The set $\sD_g$ (resp. $\sD_{g\nm 2}'$) consists of nondegenerate even (resp.~odd)
\hc\ classes. Similarly, for even (resp.~odd) $g$, the divisors in $\sD_{g\nm 2k}'$, $k=1,\dots,\frac{g}2$ (resp.~$k=1,\dots,\frac{g\nm
1}2$) are odd (resp.~even) if $k$ is odd (resp.~even).  When $g$ is odd,  $\sD_{0}$ has the same parity as
$\frac{g\np1} 2$. 

The total number of even
(resp.~odd) classes is $2^{2g-1}+2^{g-1}$ (resp.~$2^{2g-1}-2^{g-1}$).
\end{cor} 

\begin{proof} There is one pole of order 1 for a divisor in $\sD_g$, so there is no divisor of a function that adds to one of these to give
a positive divisor. All divisors in $\sD_{g\nm1}$ have support consisting of points of multiplicity 1, so only constant functions can make
them positive divisors.  

Given $r$, let $a_r$ and $b_r$ be the respective count of even and odd \hc\ classes. From Prop.~\ref{specreps}, $a_r+b_r=2^{2g}$.
If we show $a_r-b_r=2^{g}$. Then, solving for $a_r$ and $b_r$ gives the 2nd paragraph statement. 
The two cases are similar, so we do just the tougher, $r=4\ell$. For each odd $s$ use  ${{r\nm1}\choose s}\np {{r\nm1}\choose
{s\np 1}}={r\choose {s\np 1}}$. Write
$a_r-b_r$ as 
\begin{equation} \label{oddeven} {{r\nm1}\choose g}-{{r}\choose {g\nm 1}}+{{r}\choose {g\nm 3}} \np \cdots \np (-1)^{\ell\nm 1}{{r}\choose
{2}}
\np (-1)^{\ell}{{r\nm1}\choose {0}}.\end{equation} Use  ${{r\nm1}\choose g}= {{r\nm1}\choose {g+1}}$ and
${{r\nm1}\choose {0}}={{r}\choose {0}}$ to write \eqref{oddeven} -- akin to \eqref{sumbinom} -- as the real part of
$(-1)^\ell(1+i)^r/2=(\sqrt{2}e^{i\pi/4})^{4\ell}/2$. The result is $2^{2\ell-1}=2^g$, and we are done.   \end{proof} 

\begin{rem}[Nondegenerate even and odd $\theta\,$s]  Riemann showed every Riemann
surface of genus $g$ has nongenerate odd
$\theta\,$s (see \S\ref{tor2}). Apply Rem.~\ref{hypcomps} to label these $\theta\,$s explicitly using the corresponding \hc\ class. Example:
For odd $g$, the  component of $\sM_{g,\pm}$ over $\hyp_g$ of degree 1  is degenerate for $g\ge 3$. \end{rem} 

\begin{appendix} 
\renewcommand{\labelenumi}{{{\rm (\teql \alph{enumi})}}} 

\section{Interpreting points on Hurwitz spaces} I review how Hurwitz spaces interpret
Inverse Galois aspects. \cite{Fr08b} has examples. 
For
$z$ transcendental over a field $K$,  a field extension
$L/K(z)$ is  {\sl regular\/} if
$L\cap \hat K=K$. As in \S\ref{gpNC}, we abbreviate Riemann's Existence Theorem by RET. As in the rest of the paper we denote the degree of $L/K(z)$ by $n$. 

\subsection{RIGP and AIGP} \label{secrigp-aigp} The Regular  Inverse Galois
Problem over a field
$K$ ($\supset
\bQ$ for  simplicity) asks when does a group $G$ ($\le S_n$) appear as the group of a Galois extension $\hat
L/K(z)$ with  $L\cap \bar K=K$. We use the acronym RIGP for this. 

More often we seek regular extensions $L/K(z)$ with Galois closure $\hat L/K(z)$
where, with $\hat K$ the constant field of $\hat L$,  $G(\hat L/\hat K(z))$ is $G$.  Then,  
$G$ is the {\sl geometric\/} monodromy group (of $L/K(z)$), a normal
subgroup of the {\sl arithmetic\/} monodromy
$\hat G=G(\hat L/K(z))$ (also, $\le S_n$). This is a $(G,\hat G)$ {\sl realization}. Finding such a
realization for some $\hat G$ is the A(bsolute)IGP (we often care which $\hat G\,$s are
achieved).

Attached to a cover of projective curves
$\phi:X \to
\prP^1_z$ over $K$ is a collection of branch points $\bz=\{\row z r\}$
invariant under
$G_K$.  Each $z_i$ corresponds to a conjugacy class $\C_i$ in the  geometric
monodromy group. The collection $\row \C r$ is a deformation invariant of
the cover. We say all equivalence classes of covers with the datum $(G,\bfC)$ are 
in the Nielsen class
$\ni(G,\bfC)$, using the notation of \S\ref{nz.1} including for the braid and Hurwitz monodromy
groups $B_r$ and $H_r$. 

This section reviews  equivalences on covers,  and the corresponding spaces associated to  $(G,
\bfC)$ whose points interpret a solution to the RIGP or AIGP (\cite[\S2]{BaFr}, much from \cite{Fr1}
or
\cite{FrV2}). That material identifies 
$H_r$ with  the 
fundamental group of  projective $r$-space, $\prP^r$ minus the discriminant 
locus 
$D_r$: $U_r=\prP^r\setminus D_r$. The natural map 
$(\prP^1)^r\to \prP^r$, modulo the action of
$S_r$, takes the fat diagonal $\Delta_r$ to  $D_r$. This interprets 
$U_r$ as the space of  $r$ distinct unordered 
points in $\prP^1$. 

\subsection{Absolute equivalence} \label{absEquiv} With 
$\bfC$  a collection of  conjugacy classes of $G$,  
denote the automorphisms of $G$, from conjugations in $S_n$, that permute entries  
in $\bfC$ by $\Aut_\bfC(G)$. Conjugation by  $G$ induces inner 
automorphisms of $G$. Denote the subgroup of $\Aut_\bfC(G)$ induced by 
inner 
automorphisms by ${\rm Inn}(G)$. With $G'$ any group between $\Aut_\bfC(G)$ 
and $\text{\rm Inn}(G)$ (allow both extremes), denote the quotient of the conjugation  
action 
of $G'$ on $\ni(G,\bfC)$ by $\ni(G,\bfC)/G'$. When $G'=G$, use 
$\ni(G,\bfC)^\inn$ for $\ni(G,\bfC)/G'$. 

\begin{defn} \label{selfnormalizing} For $H \le G$, call  $H$ self-normalizing in $G$ if $N_G(H)=H$. 
This is equivalent to no element of $S_{n(H)}\setminus H$ centralizes $H$
\cite[Lem.~2.1]{Fr1}.  Primitivity of $T_H: G\to S_{n(H)}$ is sufficient but not necessary 
(Ex.~\ref{dihk}). 
\end{defn}

\subsubsection{Absolute Nielsen classes} 
For $H\le G$, of index $n(H)=(G:H)$, denote  
elements in $S_{n(H)}$ that normalize $G$ and 
permute 
the conjugacy classes $\bfC$ by $S_{n(H)}(G,\bfC)$. Denote the  quotient
$\ni(G,\bfC)/S_{n(H)}(G,\bfC)$ by  
$\ni(G,\bfC)^{\abs(H)}$. 

\begin{lem}[Equivalence Lemma] \label{lnz.1} The action of  $B_r$ on 
$\ni(G,\bfC)/G$ induces an action of $B_r$ (or $H_r$) on 
$\ni(G,\bfC)/G'$. 
Suppose $H\le G$ is self-normalizing in $G$ and $H$ contains no nontrivial normal subgroup of 
$G$. Let  
$T_H:G\to S_{n(H)}$ be the faithful representation on the $n(H)$ 
cosets of $H$. Then, $T_H$ extends to $G'$ if and only 
if 
$G'$ maps the conjugates of $H$ in $G$ among themselves. This case induces an $H_r$ 
equivariant map from $\ni(G,\bfC)/G'$ to 
$\ni(G,\bfC)^{\abs(H)}$. \end{lem}

\begin{proof} Two transitive permutation representations of $G$ are 
(permutation) equivalent if and only they have the same point stabilizers. Let 
$\row 
H t$ (with $H_1=H$) be the conjugates of $H$ in $G$. An element of $G$ acts 
by conjugation on $\row H t$. Since $H$ is self-normalizing, the stabilizer of 
$H$ 
in this action is just $H$. Let $T: G\to S_{n(H)}$  be the corresponding 
permutation 
representation. As $T$ and $T_H$ have the same point stabilizers, these 
representations are permutation equivalent. To extend $T_H$ to $G'$, it 
suffices to extend $T$ to $G'$. This is now 
automatic. 
Everything else in the lemma follows immediately from the definitions. 
\end{proof}

As above, $G'$ is any group between 
$\Aut_\bfC(G)$ and $G$. From fundamental group theory, as in \cite{Fr1} or \cite{FrV2},  
 $H_r$ acting on $\ni(G,\bfC)/G'$ produces an unramified cover $\Phi^{G'}: \sH(G, 
\bfC)^{G'}\to U_r$. So,  $\sH(G, \bfC)^{G'}$ is a manifold. 
Connected 
components of  $\sH(G, \bfC)^{G'}$ correspond one-one with $B_r$ 
(or $H_r$) orbits on $\ni(G,\bfC)/G'$. 

There are two natural equivalences of covers. 
\begin{edesc} \label{enz.2}
\item \label{enz.2a} $\phi_i: X_i\to \prP^1$, $i=1,2$, are
equivalent if some continuous $\alpha:
X_1\to X_2$ satisfies $\phi_2\circ \alpha=\phi_1$.
\item \label{enz.2b} As in \eql{enz.2}{enz.2a}, except there is 
$\beta\in \PGL_2(\bC)$ (\S\ref{Ancond}) with 
$\phi_2\circ\alpha=\beta\circ\phi_1$.
\end{edesc}

Call equivalence \eql{enz.2}{enz.2a} (resp.~\eql{enz.2}{enz.2b}) {\sl absolute \/} (resp. 
 {\sl reduced\/} absolute) equivalence. Any (degree $n$)  cover produces a 
 (degree $n$) permutation 
representation of 
the geometric monodromy group of the cover. 

\subsubsection{Hurwitz spaces} \label{hurSpaces}  
RET says 
elements of $\ni(G,\bfC)^{\abs(H)}$  correspond one-one to  covers of $\prP^1$ with unordered branch
points
$\bz$ (modulo absolute  equivalence) in the Nielsen class 
$\ni(G,\bfC)$. The complete collection of these equivalence classes of covers as $\bz$
varies in
$U_r$ forms a covering space, $\sH(G,\bfC)^{\abs(H)}$ of $U_r$. We drop the $H$ in $\abs(H)$ if we know
it from the context. 

There is a 
natural compactification $\bar \sH(G,\bfC)^{\abs(H)}$ of
$\sH(G,\bfC)^{\abs(H)}$ over
$\prP^r$. If $\sH'$ is any component of $\sH(G,\bfC)^{\abs(H)}$, let $\bar \sH'$ be
the normalization of $\prP^r$ in the function field of $\sH'$, and take the disjoint
union of these $\bar \sH'$, a projective variety by Grauert-Remmert 
(\cite{Grauert}, as used in \cite[\S4]{Fr1}). 

\cite[p.~56]{Se1} succinctly goes through the literature to show an analytic cover of a Zariski
subspace of a projective variety extends to a general cover of projective varieties. It notes
that 
\cite{Grauert}, extending the original cover to a cover of complete analytic spaces, is delicate.
By contrast, \fval{p.~788} notes:  When
$G$ has a self-normalizing representation, we see the affine structure on $\sH(G,\bfC)^{\abs(H)}$ from
the 1-dimensional RET at generic branch points and normalization in a function field. So, this gives
the quasiprojective structure for all absolute and inner spaces directly.  

This also works for reduced
Hurwitz spaces, equivalence \eql{enz.2}{enz.2b}. 
\cite[\S10.2]{VB} overlooked this. When
$r>4$ there are other compactifications, related to admissible covers \cite{wew}, 
but this one is useful. 

The index of $g\in S_n$ is $n-m$ where $m$ is the number of orbits of $g$. Covers in an absolute
Nielsen class
$\ni(G,\bfC)^{\abs(H)}$ have a genus $g_{G,\bfC,H}=g_\bg$ defined by 
\begin{equation} \label{RH} \text{Riemann-Hurwitz: } 2(n+g_\bg-1)=\sum_{i=1}^r \ind(g_i).\end{equation} 
This is  the genus of a cover of
$\prP^1_z$ with branch cycles given by $\bg$.  As above, equivalence
\eql{enz.2}{enz.2a} presents a space  
$\sH$ of covers with $r$ branch 
points as 
a cover of  $U_r$. Equivalence \eql{enz.2}{enz.2b} gives
a different target, 
$J_r$: the bi-quotient of  $(\prP^1)^r\setminus
\Delta_r$ by
$\PGL_2(\bC)$ (linear fractional transformations)   
and $S_r$. Here, $\PGL_2(\bC)$ acts diagonally on $(\prP^1)^r$ and $S_r$ acts by permuting these
coordinates.  These actions commute. Example: $J_4$ is the 
traditional $j$-line minus the cusp at $\infty$.  
 
\subsubsection{Interpreting  self-normalizing Def.~\ref{selfnormalizing}} \label{intSNorm} Self-normalizing is 
equivalent to  $\phi_\bp: X_\bp\to \prP^1$ -- representing  
$\bp\in \sH(G,\bfC)^{\abs(H)}$ -- has no automorphisms commuting with $\phi_\bp$. 
{\sl Fine moduli\/} for  
$\sH(G,\bfC)^{\abs(H)}$ means there is a  {\sl unique\/}  family 
of representing covers (as in \S\ref{covs-reps}) and 
self-normalizing  is equivalent  \cite[Lem.~2.1]{Fr1}.

Let $\Psi^{\abs(H)}:\sT^{\abs(H)}\to \sH(G,\bfC)^{\abs(H)}\times\prP^1$ be 
the corresponding family: For $\bp\in 
\sH(G,\bfC)^{\abs(H)}$, 
restricting $\sT^{\abs(H)}$ to the fiber of $\sT^{\abs(H)}$ over 
$\bp\times \prP^1$ is a cover representing $\bp$. \cite[\S4]{FrV2}  
shows  $\sH(G,\bfC)^\inn$ has a unique representing family 
if 
$G$ has no center. When there is a
self-normalizing
$H$,
\S\ref{innerHurw} constructs $\sH(G,\bfC)^\inn$ directly using geometric Galois closure.

\begin{exmp} \label{exA43-2} Use  $(\gamma_0,\gamma_1,\gamma_\infty)$ from the $\sh$-incidence
calculation in Prop.~\ref{A43-2}. Denote their restrictions to lifting
invariant $+1$  (resp.~-1) orbit  
 by $(\gamma_0^+,\gamma_1^+,\gamma_\infty^+)$ (resp.~$(\gamma_0^-,\gamma_1^-,\gamma_\infty^-)$).
We read indices of the $+$ (resp.~$-$) elements from the $\ni_0^+$ (resp.~$\ni_0^-$) matrix block: 
Cusp  widths over
$\infty$ add to the degree $4+2+3=9$ (resp.~$4+1+1=6$) to give
$\ind(\gamma_\infty^+)=6$ (resp.~$\ind(\gamma_\infty^+=3$); since $\gamma_1^+$ (resp.~$\gamma_1^-$) has 1
(resp.~no) fixed point and $\gamma_0^\pm$ have no fixed points, $\ind(\gamma_1^+)=4$
(resp.~$\ind(\gamma_1^-)=3$) and
$\ind(\gamma_0^+)=6$ (resp.~$\ind(\gamma_0^+)=4$). The genus of $\bar \sH_{0,\pm}$ is $g_{\pm}=0$:  
 $$2(9+g_+-1)=6+4+6=16\text{  and }2(6+{g_-}-1)=3+3+4=10.$$
\end{exmp}

\begin{exmp} \label{dihk} Consider
the  dihedral 
group, $D_{p^{k+1}}$, of order $2p^{k+1}$,  with $p$ an odd prime and $k>0$. The 
standard permutation representation of $D_{p^{k+1}}$ on $\bZ/p^{k+1}$ (an involution 
generates $H$) is imprimitive though $H$ is self-normalizing. \end{exmp}

\subsubsection{Inner Hurwitz spaces} \label{innerHurw} \cite[\S3.1.3]{BaFr} gives a Galois closure process
(more direct than in \fva) for inner Hurwitz spaces $\sH(G,
\bfC)^\inn$ from the absolute spaces. To form the Galois closure of a separable, finite flat morphism 
$\Phi: \sT\to W$ of normal varieties 
over a field $K$, form the fiber product of $\Phi$,  $n=\deg(\Phi)$ times: $$\{(\row t n)\in \sT^n\mid
\Phi(t_i)=\Phi(t_j)\}.$$ The Galois closure (of $\Phi$ over $K$) identifies with the normalization of a
non-diagonal absolutely irreducible component
$\hat \sT$ of this algebraic set. Prop.~\ref{habs-hin} shows how to go from $\sH^\abs$ 
to
$\sH^\inn$ giving the self-centralizing fine moduli condition. 

The space we seek is an unramified cover
$\sH^\inn\to \sH^\abs$: Its points
$\bp^\inn\in
\sH^\inn$ over
$\bp^\abs\in \sH^\abs$ representing an
absolute Nielsen class, represent the class of pairs $$(\hat X\to \prP^1, h: G\to \Aut(\hat 
X/\prP^1))$$ in the inner Nielsen class $\ni(G,\bfC)^\inn$. Then,  $\hat X\to \prP^1$ is a geometrically Galois 
cover with group $G$ having
branch points $\bz$; and $h$ is an isomorphism between $G$ and the 
automorphism 
group of the cover. Mapping between inner and absolute spaces takes $\bp^\inn$ to 
$\bp^\abs=\Phi^\inn_{\abs(H)}(\bp^\inn)$ and 
$\bz=\Phi^\abs\circ\Phi^\inn_{\abs(H)}(\bp^\inn)$.

Since $\hat X$ is a subset of $X^n$, identify $G$ as the subgroup of $S_n$ 
mapping $\hat X$ into itself. Given self-normalizing, Prop.~\ref{habs-hin} goes from a
family of absolute covers over any parameter space to a family of Galois closures of these covers.  

\begin{prop} \label{habs-hin} If $H\le G$ is self-normalizing (\S\ref{intSNorm}) a $K$ component of the Galois
closure of a  total family $\Phi^\abs: \sT^\abs \to \sH^\abs\times \prP^1_z$ gives $\Phi^\inn: \sT^\inn \to
\sH^\inn\times
\prP^1_z$, with
$\sH^\inn$ the normalization of $\sH^\abs$ in $\sT^\inn$. Suppose $L\supset K$. Then, an $L$ point $\bp^\inn\in
\sH(G,\bfC)^\inn$ corresponds to an $L$ component of the Galois closure of $\phi_{\bp^\abs}$. 

We have 
a sequence of covers $$\sH(G, \bfC)^\inn\mapright 
{\Phi^\inn_{\abs(H)}}\sH(G,\bfC)^{\abs(H)}\mapright{\Phi^\abs} U_r$$ 
from  inner to absolute Hurwitz space. 
For $\sH$  a component of $\sH(G, 
\bfC)^\inn$ and $\sH'$ its image by ${\Phi^\inn_{\abs(H)}}$,  $\sH\to\sH'$ is Galois and unramified.  
Its group is the subgroup of $N_{S_{n(H)}}(G,\bfC)/G$ 
that 
stabilizes the braid orbit in $\ni(G,\bfC)$ associated to $\sH'$. \end{prop}

\begin{proof}[Comments]  Let $\Phi: \sT\to \sH\times \prP^1_z$ be any $r$ branch point, degree $n$, family of
$\prP^1_z$ covers.  Form the $n$-fold fiber product of $\Phi$. Normalize the result, and take
a (non-diagonal) connected component $\hat \sT$. This induces  $\Phi^{(n)}:\hat \sT \to \sH\times \prP^1_z$. 
Then, take $\hat \sH$ to be the normalization of $\sH$ in the function field of
$\hat \sT$. 

We can (geometrically) compare the  Galois closure of the sphere covers at any fiber $\Phi_\bp:\sT_\bp\to
\bp\times
\prP^1_z$, for $\bp\in \sH$, to the components of the fiber
$\hat
\sT_\bp \to \bp\times \prP^1_z$. We expect several geometric copies ($u$, say, locally constant in $\bp$) of the
Galois closure of
$\sT_\bp\to
\bp\times \prP^1_z$ for each $\bp$. Then, $u$ is the degree of  $\hat \sH$ over $\sH$. Local constancy 
of $u$ results from deforming classical generators of the $r$-punctured sphere, as in 
\cite[\S4]{Fr1}. This works in positive characteristic only for tamely ramified covers. 
\end{proof}
 
\begin{exmp}[Ex.~\ref{mod8exs} cont.] Even for $G=A_n$, the value of $u=\deg{\Phi^\inn_{\abs(H)}}$ in
Prop.~\ref{habs-hin} is nonobvious. It is 1 for $K=\bar \bQ$ in Ex.~\ref{mod8exs}, $n\equiv 1 \mod 8$, while
it is 2 for
$K=\bQ$. It is always 2 in Ex.~\ref{mod8exs}, $n\equiv 5 \mod 8$, or in 
Thms.~\ref{thmA} or \ref{thmB}.  
\end{exmp} 

Here is the RIGP-AIGP interpretation, a variant on  \cite[Main Thm.]{FrV2}. 

\begin{prop} \label{rigp-aigp} A regular extension $L/K(z)$  over $K$ 
in the Nielsen class
$\ni(G,\bfC)^{\abs(H)}$ corresponds to $\bp\in \sH(G,\bfC)^{\abs(H)}(K)$. If $H$ is self-normalizing,
then conversely, any  $\bp\in
\sH(G,\bfC)^{\abs(H)}(K)$ corresponds to a regular extension over $K$ in the Nielsen class. This gives a
$(G,\hat G)$ ($K$) realization (as in \S\ref{secrigp-aigp}) with $\hat G/G=G(K(\hat \bp)/K)$ for $\hat \bp\in
\sH(G,\bfC)^{\inn}$ over
$\bp$. 

Regular Galois  $\hat L/K(z)$ in 
$\ni(G,\bfC)^\inn$ corresponds to $\bp\in \sH(G,\bfC)^{\inn}(K)$: A solution of
the RIGP over $K$. Conversely, for $H$ centerless, any  $\hat \bp \in  
\sH(G,\bfC)^{\inn}(K)$ gives a regular Galois extension in the Nielsen class.
\end{prop}

\subsection{Reduced versions of moduli 
spaces} \label{rnz.4}  Points on nonreduced  absolute moduli spaces correspond 
to sphere covers in a given Nielsen class. Suppose $G$ is centerless and $H\le G$ is self-normalizing. We see
the Inverse Galois problem structure from the relation  between $\sH(G,\bfC)^{\abs(H)}$ and $\sH(G,\bfC)^\inn$
(\S \ref{nth.1}).  Even with self-normalizing, 
conveniently interpreting  
points for the reduced moduli spaces 
$\sH(G,\bfC)^{\abs(H)}/\PGL_2(\bC)$ and 
$\sH(G,\bfC)^\inn/\PGL_2(\bC)$ depends on the situation.  
 
To see this, consider the case tied to modular curves: 
$G=D_p$ ($p$ an odd prime) with $\bfC=\bfC_{2^4}$,  four repetitions 
of the involution conjugacy class (\cite{Fr2} and \cite[\S5.1--\S5.2]{DFr94}) Then, 
$\sH(D_p,\bfC_{3^4})^{\inn}/\PGL_2(\bC)$  identifies with the classical  $Y_1(p)$ (modular curve without
cusps). Points of $Y_1(p)$ correspond to equivalence classes of pairs $(E,\be)$ with $E$ 
an elliptic 
curve and $\be$ an order $p$ torsion point on $E$. 

If multiplication by $-1$ is the only automorphism of $E$, then $(E,-
\be)$ is also in this class. Also, if $(E,\be')$ is here, 
then $\be'=\pm \be$. Let $O$ be 
the $\PGL_2(\bC)$ orbit in $U_r$ that maps to the $j$-line value of $E$. Choose any  $\bx\in O$. 
To show   
$Y_1(p)$ is  $\sH(D_p,\bfC_{3^4})^{\inn}/\PGL_2(\bC)$ requires 
 recovering $(\hat X\to \prP^1, h: D_p\to \Aut(\hat 
X/\prP^1))$ up to conjugation by $D_p$ from the triple $(E,\be,\bx)$ up to 
equivalence. This works; $\bx$ determines a degree 2 map from 
$E/\lrang{\be}\to \prP^1$. The sequence $E\to E/\lrang{\be}\to \prP^1$ is 
$E\to\prP^1$, a (geometric) Galois cover with group $D_p$. The collection of 
points $\pm \be$ determines an isomorphism of this group with $D_p$ up to 
conjugacy by $D_p$. 

Finally, the reduced absolute space in this case interprets naturally as
$Y_0(p)$ in a diagram coming from the map $\sH(G,\bfC)^\inn\to
\sH(G,\bfC)^{\abs(H)}$. That is, map the equivalence class of $(E,\be)$
to the equivalence class $(E,E/\lrang{\be})$. 

Now we quote the literature for the algebraic
structure on the reduced  spaces. 

\begin{prop}[Reduction Proposition]\label{pnz.5} With  $G\le G'\le \Aut_\bfC(G)$ as
above, the quotients 
$U_r^\rd=J_r$ (\S\ref{hurSpaces}), and $\sH(G, \bfC)^{G'}/\PGL_2(\bC)\eqdef \sH(G, \bfC)^{G',\rd} $ are affine
algebraic varieties with an induced finite map $\Psi^\rd: \sH(G, 
\bfC)^{G',\rd}\to J_r$. \end{prop}

\begin{proof} The discriminant locus in $\prP^r$ is a hypersurface. Therefore, 
its complement, $U_r$, is affine. 
Since 
$\sH(G, \bfC)^{G'}$ doesn't ramify over $U_r$, it is 
normal. Further,  the Grauert-Remmert version of RET  
applies \cite[p.~442]{Hart}. There is a unique normal projective 
variety $\overline{\sH(G, \bfC)^{G'}}$ with a finite covering  $\bar\Psi: 
\overline{\sH(G, \bfC)^{G'}}\to \prP^r$ whose restriction over $U_r$ gives 
the natural covering map $\Psi: \sH(G, \bfC)^{G'}\to U_r$.  
Since  $\Psi$ is a finite cover, it is an affine cover and $\sH(G, 
\bfC)^{G'}$ is 
an affine variety.  

Apply \cite[Thm.~1.1, p.~27]{Mu4}  to 
the affine scheme $\sH(G, \bfC)^{G'}$ and the  reductive group $\PGL_2(\bC)$. If  
the action of $\PGL_2(\bC)$ is closed, $\sH(G, 
\bfC)^{G',\rd}$ with the natural map $\Gamma: \sH(G, 
\bfC)^{G'}\to \sH(G, \bfC)^{G',\rd}$ is a {\sl universal geometric 
quotient\/} and  affine. We see the $\PGL_2(\bC)$ orbit of  
$\bp\in \sH(G, \bfC)^{G'}(\bC)\leftrightarrow \phi$ (with branch set $\bz$) is closed by considering any limit 
$\alpha_n\circ\phi$ with $\alpha_n\in \PGL_2(\bC)$. This comes to showing any limit of 
$\bz$  under $\{\alpha_n\}_{n=0}^\infty$ is in $U_r$; or the
analog for  $(\row z r)\in U^r\eqdef (\prP^1)^r\setminus \Delta_r$
(\S\ref{gpNC}) replacing $\bz$. What 
$\alpha_n$  does to $(\row z r)$ determines it, so this determines the limit of the
$\alpha_n\,$ in $\PGL_n(\bC)$. When
$r=4$, $\sH(G, \bfC)^{G',\rd}$ is a curve; its completion ramifies over $\bar J_4=\prP^1_j$ with
connected components  one-one to those of 
$\sH(G, \bfC)^{G'}$. For $r\ge 5$, $\sH(G, \bfC)^{G',\rd}$ may have singularities. \end{proof}

For $r=4$, recall  $\sQ''$ of \eql{cuspgp}{cuspgpa} and the elements $\gamma_0$ and $\gamma_1$ of
\eql{cuspgp}{cuspgpc}.  
\begin{rem}[Fine reduced moduli] \S\ref{intSNorm} discussed fine
moduli for inner (no center in $G$) and absolute ($H$ self-normalizing) Hurwitz spaces. These
respective conditions must hold for reduced fine moduli. Yet, you need more \cite[Prop.~4.7]{BaFr}. 
For
$r\ge 5$, it is that the reduced space $\sH^*/\PGL_2(\bC)$ ($*=\inn$ or
$\abs$) is nonsingular. The analog  for $r=4$ is that $\gamma_0$ and $\gamma_1$ on {\sl reduced Nielsen
classes}, $\ni(G,\bfC)^*/\sQ''$,  have no fixed points. For
$r=4$ there is also one more condition: $\sQ''$ has  {\sl only maximal length\/}  (that would be 4; action
through a Klein 4-group) orbits on 
$\ni(G,\bfC)^*$. \cite[Ex.~8.5]{BaFr} notes that neither condition holds for $\ni(A_5,\bfC_{3^4})^\inn$ (in
Thm.~\ref{thmA}). The two components $\sH_{\pm}$ attached to $\ni(A_4,\bfC_{\pm 3^2})^{\inn,\rd}$ both fail the
$\sQ''$ condition (its orbit lengths on both are 2), while neither $\gamma_0$ or $\gamma_1$ have fixed
point on $\bar \sH_{-}$, but $\gamma_1$ does on  $\bar \sH_{+}$ \cite[\S6.6.3]{lum}. 
\end{rem} 

\section{Producing \hc\ differentials} \label{cocycleNot}

Let $X$ be a compact Riemann surface. Riemann used
{\sl certain\/} theta functions on $X$ to give  a constructive approach to  all its  functions and
differential forms. We collects observations on the key ingredient, \hc\ differentials. 

\subsection{$\Theta$ data} \label{famCovers} Suppose
$X$ appears in a smooth family 
$\Psi: \sX\to\sP$ of Riemann surfaces as the fiber $X_{\bp}$ over 
$\bp\in\sP$.  Let $(\balpha\eqdef (\row \alpha g),\bbeta\eqdef (\row \beta g)$ be a {\sl canonical homology
basis\/} (of
$H_1(X,\bZ)$) for $X$: The cup-product image of $(\alpha_i,\beta_j)$ (resp. $(\alpha_i,\alpha_j)$ and
$(\beta_i,\beta_j)$) in $H_2(X,\bZ)\equiv \bZ$ is $\delta_{i,j}$ (resp.~0) for all $1\le i,j\le g$. 

If 
so, we can represent the cup product as a skew-symmetric
$2g\times 2g$ matrix
$E$  which extends to an $\bR$-bilinear form on $H_1(X,\bR)$. We can write
$$\Bigl(\frac{\balpha}
{\bbeta}\Bigr)E(\balpha^\tr|\bbeta^\tr)  =J_{2g} \text{ (on the left the $i$th row is $\alpha_i$, $1\le i\le g$,
etc.)}.$$ Then,  an element
$U=\smatrix A B C D \in \Sp_{2g}(\bZ)$  as in \S\ref{tor1}, acting on the left, transforms $\bigl(\frac{\balpha}
{\bbeta}\bigr)$ so as to 
preserve the pairing. \cite[\S3.1]{Sh98} calls  $f(\bw)$ a $\theta$ function  if $f(\bw
+\gamma)=f(\bw)e^{l_\gamma(\bw)+c_\gamma}$ for $\gamma\in H_1(X,\bZ)$, $l_\gamma(\bw)$ linear in $\bw$ and
$c_\gamma\in \bC$. 

This treatment expresses that any $\theta$ comes from a slight
generalization of $E$. Conversely, there is a $\theta$ on a complex torus if and only if there is an
associated $E$. Riemann's
$\theta$ is that attached to this particular $E$.  

Denote positive divisors of degree $g-1$ on
$X$ by $W_{g-1}$.  Note that $W_{g\nm1}$ is independent of $(\balpha,\bbeta)$, but the \hc\ class and
$\Theta_{X_\bp}$ is not.  

\cite[App.~B]{BaFr} explains that for his generalization of Abel's Theorem, 
Riemann wanted 
$\theta_X$   {\sl odd\/} and nondegenerate: 
 $\theta$  has nonzero gradient  (see \S\ref{comphc}) at the 
origin of $\widetilde \Pic^{(0)}=\bfC^g$, as in \S\ref{hcfamilies}. This
gave his generalization of Abel's Theorem. So, instead of the even  $\theta(\bl,\bw)$ in \eqref{rtheta}, he
needed a different  $\theta$ with 2-division characteristic. We want $\theta_X$ even, but also nondegenerate: 
not zero at the origin.  

In Thm.~\ref{thetaNullRes}, we needed a point on the Hurwitz space so that {\sl no\/} 
$\theta$ with 2-division characteristic  was zero at the origin. Riemann showed some $\theta\,$s with
2-division characteristic will work at each Riemann surface. Still, for a general family of covers 
$\sH$ with odd order branching (as in Prop.~\ref{thetaChoices}) it may be that none of those match the
\hc\ class  
$\theta$-null defined by Lem.~\ref{lnth.2}. Here is the exact criterion. 

\begin{prop} The corresponding $\theta$-null is nontrivial exactly when there is  
$\bp\in\sH$ so the
\hc\ class $\delta_\bp$ attached to $\bp$ contains no positive divisor. \end{prop}   

\newcommand{\hH}{\text{\rm Hol}}

\subsection{\hc\ differentials from coordinate charts}  \label{cocycs} 
Suppose $X$ is an $n$-dimensional complex manifold. Let  
$\{U_\alpha,\phi_\alpha\}_{\alpha\in I}$ be the coordinate chart, with 
$\{\psi_{\beta,\alpha}=\phi_\beta\circ \phi_\alpha^{-1}\}_{\alpha,\beta\in I}$ the
corresponding  transition functions. Each
$\psi_{\beta,\alpha}$ is one-one and analytic on an open subset of $\bC^n$ whose
coordinates we label $z_{\alpha,1},\dots,z_{\alpha ,n}$. 

\subsubsection{Reminder on cocycles} \label{transfunctions} Denote the $n\times n$ {\sl complex Jacobian
matrix\/}  for $\psi_{\beta,\alpha}$ by $J(\psi_{\beta,\alpha})$. Call the matrices
$\{J(\psi_{\beta,\alpha})\}_{\alpha,\beta\in I}$ the (transformation)  {\sl cocycle\/} attached to
meromorphic differentials. 

Similarly
$\{J(\psi_{\beta,\alpha})^{-1}\}_{\alpha,\beta\in I}$ is the cocycle attached to meromorphic
tangent vectors. Recall the notation for $n\times n$ matrices, $\bM_n(R)$ with entries in an integral
domain $R$ and for the {\sl invertible\/} matrices $\GL_n(R)$ with entries in $R$ under
multiplication.  Cramer's rule says for each $A\in \bM_n(R)$ there is an adjoint matrix $A^*$ so
that $AA^*$ is the scalar matrix $\det(A)I_n$ given by the determinant of $A$. This shows the
invertibility of
$A\in
\bM_n(R)$ is equivalent to $\det(A)$ being a {\sl unit\/} (in the multiplicatively
invertible elements 
$R^*$) of $R$. Denote the $n\times n$ identity matrix (resp.~zero matrix) in $\GL_n(R)$ by $I_n$
(resp.~${\pmb 0}_n$). If $U\subset X$ is an open set, denote the {\sl holomorphic\/} functions on $U$ by
$\hH(U)$.    

\newcommand{\sGL}{\tsp{G}{\tsp{L}}} \newcommand{\sbM}{\tsp{M}}
\begin{defn}[1-cycocle] Suppose $g_{\beta,\alpha}\in \GL_n(\hH(U_\alpha\cap U_\beta))$,
$\alpha,\beta\in I$. Assume  $g_{\gamma,\beta}g_{\beta,\alpha}=g_{\gamma,\alpha}$ for all
$\alpha,\beta,\gamma\in I$ on $U_\alpha\cap U_\beta\cap U_\gamma$ (if nonempty).
Then, $\{g_{\beta,\alpha}\}_{\alpha,\beta\in I}$ is a multiplicative {\sl 1-cocycle with values in
$\sGL_{n,X}$}.  Similarly, suppose  $g_{\beta,\alpha}\in \bM_n(\sH(U_\alpha\cap U_\beta))$,
$\alpha,\beta\in I$. Suppose $g_{\gamma,\beta}+g_{\beta,\alpha}=g_{\gamma,\alpha}$ for all
$\alpha,\beta,\gamma\in I$ on $U_\alpha\cap U_\beta\cap U_\gamma$. Then, 
 $\{g_{\beta,\alpha}\}_{\alpha,\beta\in I}$ is an additive {\sl 1-cocycle with values in
$\bM_{n,X}$}. 
\end{defn}

When there are k-cocycles, there are also (k-1)-chains and their associated k-boundaries.
We write the definition for
$\GL_n$, recognizing there are analogous versions for all other types of cocycles. 

\begin{defn}[1-boundary] Consider $u_\alpha\in\GL_n(\hH(U_\alpha))$,
$\alpha\in I$. If $$g_{\beta,\alpha}=u_\beta(u_\alpha)^{-1}\text{ in   }U_\alpha\cap U_\beta$$ (if nonempty) for
all
$\alpha,\beta\in I$), then $\{g_{\beta,\alpha}\}_{\alpha,\beta\in I}$ is a 1-cocycle, called a {\sl
1-boundary\/} with values in
$\sGL_{n,X}$.  Call the set
$\{u_\alpha\}_{\alpha\in I}$ a {\sl 0-chain\/} with values in $\sGL_{n,X}$. 
\end{defn}

\subsubsection{Square hypothesis versus sections of \hc\ cocycles} \label{comphc}  
Let $\omega$ -- represented by $\{f_\alpha(z_\alpha)dz_\alpha\}_{\alpha\in I}$ -- be a differential on a
compact Riemann surface
$X$ satisfying Square Hypothesis
\eqref{merSquare}. Example: One  produced by the differential of a function with odd order branching, as
in Lem.~\ref{lnth.2}. If each $U_\alpha$ is simply connected, then $f_\alpha(z_\alpha)$ has two
meromorphic square-roots
$\pm h_\alpha(z_\alpha)$ on
$\phi_\alpha(U_\alpha)$.   

Use transition function notation $\psi_{\beta,\alpha}$ from
\S\ref{transfunctions} to consider existence of a well-defined  {\sl
\hc\ differential\/}  whose divisor (on $X$) is
$\pmb k=(\omega)/2$.  For that, we must choose signs on the $h_\alpha\,$s so as, on
$\phi_\alpha(U_\alpha\cap U_\beta)$, to  assert equality: \begin{equation} \label{half-values}
\tau_\alpha(z_\alpha)= h_\alpha(z_\alpha)\sqrt{dz_\alpha} =
\tau_\beta(\psi_{\beta,\alpha}(z_\alpha))=h_\beta(\psi_{\beta,\alpha}(z_\alpha))\sqrt{d\psi_{\beta,\alpha}(z_\alpha)}.\end{equation}
If so, call such a collection
$\{h_\alpha\}_{\alpha\in I}$ a (meromorphic) section of (the {\sl bundle\/} of) $\pmb k$. 

\newcommand{\sg}{\bg}

\begin{prop} \label{half-coboundary} Assume  $U_\alpha\cap
U_\beta$,
$(\alpha,\beta)\in I\times I$ is simply connected and you have chosen  
$\sqrt{J(\psi_{\beta,\alpha})}=k_{\beta,\alpha}$ on
$U_\alpha\cap U_\beta$.  
If you can sign the $h_\alpha\,$s to give equality in 
\eqref{half-values} for all $(\alpha,\beta)$, then  
$\{k_{\beta,\alpha}\}_{(\alpha,\beta)\in I\times I}=\pmb k$ is a (\hc\!\!) cocycle. 

Suppose $\pmb k$ and $\pmb k'$ are two \hc\ cocycles  that differ by a coboundary. Then, one has a section if
and only if the other does. So, there are  $2^{2g}$ such \hc\ 
cocycles modulo coboundaries corresponding to $h_\alpha$ sign changes. 
\end{prop}  

\begin{proof} Assume there is a choice of signs that gives equality in \eqref{half-values}. Then, on a triple
$\alpha,\beta,\gamma$ with $U_\alpha\cap U_\beta\cap U_\gamma\ne \emptyset$, dropping the extra
evaluation notation, check the compatibility of the equations $h_\gamma=h_\alpha
k_{\gamma,\alpha}$, $h_\gamma=h_\beta
k_{\gamma,\beta}$ and $h_\beta=h_\alpha
k_{\beta,\alpha}$. Substitute the 3rd in the 2nd, then equate to the first to see $\pmb k$ is a co-cycle.  

If $\pmb k$ and $\pmb k'$ differ by a coboundary from $\pmb m=\{m_\alpha (\in \{\pm 1\})\}_{\alpha\in I}$, then
we can multiply the section $\{h_\alpha\}_{\alpha \in I}$ by $\pmb m$ to get $\{h_\alpha\cdot
m_\alpha\}_{\alpha
\in I}$ to get a section of $\pmb k'$. More generally, all choices of $\bm$ as above, modulo coboundaries, that give allowable sign changes in the $h_\alpha\,$s correspond to homomorphisms of the fundamental group of $X$ into $\{\pm 1\}$. There are $2^{2g}$ such homomorphisms.   
\end{proof} 

The square of a \hc\ differential (\eqref{half-values} holds) gives  a differential satisfying the Square
Hypothesis.  Yet, the converse may not hold. Consider those \hc\ classes on a compact surface $X$ for
which there is an $h_\delta$ cocycle section  with no poles.  
\cite[Thm.~4.21, due to Mumford]{Fay} says the collection of their 
squares  generates the space of holomorphic differentials. \cite[p.~16]{Fay}  says this
produces a nondegenerate (after Def.~\ref{nontrivthetanull}) odd
\hc\ class $\delta$. Since Riemann used this result, it is disconcerting the \lq\lq proof\rq\rq\
quotes only modern papers. 

Given such a $\delta$, let $\mu$ be the 2-division point giving $\theta(\bl,\bw+\mu)$ with 
$W_{g\nm1,X}-\delta$ (compatible with \S\ref{tor1} notation) as zero divisor. This produces the {\sl prime
form\/}, $$\frac{\theta(\psi(x-x')+\mu)}{h_\delta(x)h_\delta(x')}\text{ with }x,x'\in X, \text{ and $\psi$ an
embedding of
$X$ in
$\Pic^{(0)}(X)$}.$$ 
From the prime form,  \cite[Chap.~II]{Fay} constructs 
all the important objects on
$X$, including functions as expressing Riemann's generalization of Abel's Theorem. 
 Unattributed, but calling
it classical, \cite[p.~17]{Fay} directly presents the square of the prime form for a surface
$X$ presented as a $\prP^1_z$ cover. 

\begin{prob} Work out  ingredients of the prime form along Hurwitz spaces of odd order branching, 
in Thms.~\ref{thmA} and \ref{thmB} and the many more spaces in Rem.~\ref{ModTowers}. \end{prob}

Thm.~\ref{thetaNullRes} uses Riemann-Roch  to produce sections of any \hc\ cocycle with degree
bounds on the polar divisor $D$ of the meromorphic section. Further, the existence of \hc\ sections is what
makes the argument using  the map $H^0(X,O(\pmb k,D)) \to H^1(X\setminus D, \bC)$ work. Still, there may be a
difference between even and odd $\pmb k\,$s in the fiber dimensions of this \lq\lq quadratic\rq\rq\ map. 

\begin{prob} With $\phi: X\to\prP^1_z$ having odd order branching, characterize when 
some section of a \hc\ bundle (\S\ref{halfcan}) has divisor half of $(d\phi)$. 
\end{prob} 
 
\end{appendix} 



\providecommand{\bysame}{\leavevmode\hbox to3em{\hrulefill}\thinspace}

\end{document}